\documentclass{amsart}
\usepackage{amssymb}
\usepackage[all]{xy}

\DeclareMathOperator{\clos}{clos}
\DeclareMathOperator{\id}{id}
\DeclareMathOperator{\Ran}{Ran}
\DeclareMathOperator{\Ker}{Ker}
\DeclareMathOperator{\Ext}{Ext}
\DeclareMathOperator{\Mor}{Mor}
\DeclareMathOperator{\Ob}{Ob}
\DeclareMathOperator{\Sys}{Sys}
\DeclareMathOperator{\Mod}{Mod}
\DeclareMathOperator{\Tfn}{Tfn}
\DeclareMathOperator{\Cfn}{Cfn}
\DeclareMathOperator{\Inv}{Inv}

\begin{document}

\newtheorem{theorem}{Theorem}[section]
\newtheorem*{theorema}{Theorem A}
\newtheorem*{theoremb}{Theorem B}
\newtheorem{lemma}[theorem]{Lemma}
\newtheorem*{lem}{Lemma}
\newtheorem{lemm}{Lemma}
\renewcommand{\thelemm}{{\it(\roman{lemm})}}
\newtheorem{prop}[theorem]{Proposition}
\newtheorem*{prp}{Proposition}
\newtheorem{cor}{Corollary}
\renewcommand{\thecor}{}
\newtheorem*{defin}{Definition}
\newcommand{\myproof}{{\noindent\it Proof. }}
\newcommand{\proofend}{$\;\square$\\[3pt]}
\newcommand{\exmpl}{\noindent\textbf{Example.}}
\newcommand{\myrem}{\noindent\textbf{Remark.}}

\title[Factorizations and invariant subspaces]
{Factorizations and invariant subspaces\\ for weighted Schur classes}

\author[Alexey Tikhonov]{Alexey Tikhonov}

\address{
Department of Mathematics\\
Taurida National University\\
Yaltinskaya str., 4\\
Simferopol 95007 Crimea\\
Ukraine }

\email{tikhonov@club.cris.net}



\date{October 2, 2005}

\begin{abstract}
We study factorizations of operator valued functions of weighted Schur classes over
multiply-connected domains. There is a correspondence between functions from weighted
Schur classes and so-called ``conservative curved'' systems introduced in the paper.
We show that the fundamental relationship between invariant subspaces of the main
operator of a conservative system and factorizations of the corresponding operator
valued function of Schur class, which is well known in the case of the unit disk, can
be extended to our case. We develop new notions and constructions and discuss changes
that should be made to the standard theory to obtain desired generalization.
\end{abstract}

\maketitle \thispagestyle{empty}

\setcounter{section}{-1}

\section{Introduction}

It is well known~\cite{NF,Br} that there is an one-to-one correspondence between (simple) unitary
colligations
\[
{\mathfrak A}=\left(\begin{array}{cc}T&N\\M&L\end{array}\right) \in{\mathcal L}(H\oplus{\mathfrak
N}, H\oplus{\mathfrak M}),\quad {\mathfrak A}^*{\mathfrak A}^{}=I,\;\,{\mathfrak A}^{}{\mathfrak
A}^*=I
\]
and operator valued functions $\Theta(z)$ of the Schur class
\[
S=\{\Theta\in H^\infty({\mathbb D},{\mathcal L}({\mathfrak N},{\mathfrak M})) :
||\Theta||_\infty\le 1\}\,.
\]
Here $\,H,{\mathfrak N},{\mathfrak M}\,$ are separable Hilbert spaces and $\,{\mathcal
L}({\mathfrak N},{\mathfrak M})\,$ is the space of all bounded linear operators acting from
$\,{\mathfrak N}\,$ to  $\,{\mathfrak M}\,$. The mapping defined by the formula
$\;\Theta(z)=L^*+zN^*(I-zT^*)^{-1}M^*\,,\,|z|<1$ is one of the directions of the above mentioned
correspondence. The operator valued function $\,\Theta(z)\,$ is called the characteristic function
of the unitary  colligation ${\mathfrak A}$ and its property $\,||\Theta||_\infty\le 1\,$ is the
consequence of the unitarity property of the colligation ${\mathfrak A}$.


The reverse direction of the correspondence is realized via functional model~\cite{NF,Br}, whose
essential ingredients are Hardy spaces $H^2$ and $H^2_-$ (see~\cite{Du}). These two sides of the
theory (unitary colligations and Schur class functions) are equipollent: both have simple, clear
and independent descriptions and we can easily change a point of view from unitary colligations to
Schur class functions and back. This context gives a nice opportunity to connect operator theory
and function theory in a very deep and fruitful manner~\cite{Ni}.

One of the cornerstones of the theory is the link (see~\cite{NF,Br}) between factorizations of
characteristic function $\Theta(z)$ and invariant subspaces of operator $T$, which goes back
to~\cite{Li} and~\cite{B}. The most simple way to explain this connection is to look at it from the
point of view of systems theory and to employ the well-known correspondence between unitary
colligations ${\mathfrak A}$ and conservative linear discrete time-invariant systems
$\,\Sigma=(T,M,N,L;H,\mathfrak{N},\mathfrak{M})\,$ (see, e.g.,~\cite{BGK})
\[
\left\{%
\begin{array}{rccl}
    x(n+1)&=&T x(n)+Nu(n)\,, &\quad x(n)\in H,\;\,u(n)\in \mathfrak{N}\,, \\
    y(n)&=&Mx(n)+Lu(n)\,, &\quad  y(n)\in \mathfrak{M},\;\,n\ge 0\,.
\end{array}%
\right.
\]
The conservative property (the property of energy balance) of $\,\Sigma\,$ corresponds
to the unitarity property of the colligation ${\mathfrak A}$. If we send the sequence
$\,u(n)\,$ into the system $\,\Sigma\,$ with the initial state $\,x(0)=0\,$, we get
the identity $\,\hat{y}(z)=S(z)\hat{u}(z)\,$, where $\,\hat{u}(z)=\sum_{n=0}^{\infty}
z^nu(n)\,$, $\,\hat{y}(z)=\sum_{n=0}^{\infty} z^ny(n)\,$, and
$\;S(z)=L+zM(I-zT)^{-1}N\,$. Note that the transfer function $\;S(z)\,$ of the system
$\Sigma$ is equal to $\,\Theta^{\sim}(z)\,$, where
$\,\Theta^{\sim}(z):=\Theta^*(\bar{z})\,$ is the dual function to the characteristic
function $\,\Theta(z)\,$ of the unitary colligation $\,{\mathfrak A}\,$.

Sending the output of a system $\Sigma_2=(T_{2},M_{2},N_{2},L_{2};H_{2},\mathfrak{N},\mathfrak{L})$
into the input of a system $\Sigma_1=(T_{1},M_{1},N_{1},L_{1};H_{1},\mathfrak{L},\mathfrak{M})$, we
obtain the cascade system
$\Sigma_{21}:=\Sigma_{2}\cdot\Sigma_{1}=(T_{21},M_{21},N_{21},L_{21};H_{21},\mathfrak{N},\mathfrak{M})$.
It is clear that the transfer function $\,S_{21}(z)\,$ of the system $\Sigma_{21}$ is the product
of the transfer functions of systems $\Sigma_1$, $\Sigma_2$ and it is easily shown that
\[
\Sigma_{21}=
\left(\left(%
\begin{array}{cc}
  T_1 & N_1M_2 \\
  0 & T_2 \\
\end{array}%
\right),\; (M_1,L_1M_2)\,,\;
\left(%
\begin{array}{c}
  N_1L_2 \\
  N_2 \\
\end{array}%
\right)\,,\; L_2 L_1\,\right),
\]
where $H_{21}=H_{1}\oplus H_{2}$. The subspace $H_1$ is invariant under the operator $T_{21}$ and
therefore, if we fix the characteristic function $\Theta_{21}(z)$, one may hope to study invariant
subspaces of the operator $T_{21}$ using this approach. Unfortunately, there are some pitfalls for
this: the operator $T_{21}$ can vary when we run over all factorizations of $\Theta_{21}(z)$. More
precisely, the variable part is the unitary component $T_{21u}$ from the decomposition
$T_{21}=T_{21s}\oplus T_{21u}$ into completely non-unitary and unitary parts~\cite{NF}. In this
connection, recall that any conservative system $\Sigma$ can be uniquely represented in the form
$\Sigma=\Sigma_s\oplus\Sigma_u$, where $\,H_s=H_c\vee H_o,\,H_u=H\ominus H_s\,$, $\,H_c=\vee_{n\ge
0}T^n N(\mathfrak{N}),\,H_o=\vee_{n\ge 0}(T^*)^nM^*(\mathfrak{M})\,$. Here $\Sigma_s$ and
$\Sigma_u$ are the simple  and  ``unitary'' parts of the system $\Sigma$, respectively. A system
$\,\Sigma\,$ is called simple if $\,H=H_s\,$. A system $(T,0,0,0;H,\{0\},\{0\})$ is called ``purely
unitary'' system if the operator $T$ is unitary.

B.~Sz.-Nagy and C.~Foia\c{s} established the following criterion (see~\cite{NF,Br}): \textit{the
product of consecrative systems $\Sigma_{21}=\Sigma_{2}\cdot\Sigma_{1}$ is simple if and only if
the corresponding factorization $\Theta_{21}(z)=\Theta_{2}(z)\Theta_{1}(z)$ is regular}. The
product $\Theta_{21}(z)=\Theta_{2}(z)\Theta_{1}(z)$ of Schur class functions is called
regular~\cite{Br} if
\[
\Ran\,(I-\Theta^*_{2}(z)\Theta_{2}(z))^{1/2}\cap
\Ran\,(I-\Theta_{1}(z)\Theta^*_{1}(z))^{1/2}=\{0\}\,,\quad\mbox{a.e.}\; z\in\mathbb{T}\,.
\]
This definition of regularity is equivalent to standard one from~\cite{NF}.

Moreover, B.~Sz.-Nagy and C.~Foia\c{s} described (Theorems VII.1.1 and VII.4.3
in~\cite{NF}) an order preserving one-to-one correspondence between regular
factorizations of a characteristic function and invariant subspaces of the
corresponding model operator, where the order relation for invariant subspaces is the
ordinary inclusion and for factorization the order relation is
$\,\Theta_{2}\Theta_{1}\prec\Theta_{2}'\Theta_{1}'\,$, where we write
$\,\Theta_{2}\Theta_{1}\prec\Theta_{2}'\Theta_{1}'\,$ if there exists $\,\theta\in
S\,$ such that $\,\Theta_{2}=\Theta_{2}'\theta\,$ and
$\,\Theta_{1}'=\theta\Theta_{1}\,$. Extension of this correspondence between
factorizations  and invariant subspaces to the case of weighted Schur classes is the
main aim of the present paper.

We shall consider operator valued functions (or rather, sets of operator valued functions) of
\textit{weighted Schur classes} $S_{\Xi}$  :
\[
\begin{array}{ll}
S_{\Xi}\;:=&\{\;(\Theta^{+},\Xi_+,\Xi_-)\; :\; \Theta^{+}\in H^\infty(G_+,{\mathcal
L}({\mathfrak
N}_+,{\mathfrak N}_-))\,,\\[2pt]
&\quad\forall\;\zeta\in C\;\;\forall\;n\in {\mathfrak N}_+\quad
||\Theta^{+}(\zeta)n||_{-,\zeta}\le ||n||_{+,\zeta}\}\,,
\end{array}\eqno{\rm{(Cfn)}}
\]
where ${\mathfrak N}_{\pm}$ are separable Hilbert spaces; $G_+$ is a finite-connected
domain of the complex plane $\,\mathbb{C}\,$ bounded by a rectifiable Carleson curve
$C$, $\,G_-=\mathbb{C}\setminus \clos\,G_+\,$ and $\infty\in G_-$;\, $\,\Xi_{\pm}\,$
are operator valued weights such that $\Xi_{\pm},\Xi_{\pm}^{-1}\in
L^\infty(C,{\mathcal L}({\mathfrak N}_{\pm}))$, $\Xi_{\pm}(\zeta)\ge 0,\;\zeta\in
C\;$, and $\,||n||_{\pm,\zeta}:=(\Xi_{\pm}(\zeta)n,n)^{1/2},\;n\in {\mathfrak
N}_{\pm}$. We shall also use the parallel notation $\,\Theta\in \Cfn\,$ whenever
$\,\Theta\in S_{\Xi}\,$.

First, we recall the construction of free \textit{functional model} of Sz.-Nagy-Foia\c{s} type
(see~\cite{T1,T2,NV}). Let $\,\Pi=(\pi_+,\pi_-)\,$ be a pair of operators $\pi_{\pm}\in {\mathcal
L}(L^2(C,{\mathfrak N_{\pm}}),{\mathcal H})$ such that
\[
\begin{array}{clcl}
(i)_1&(\pi_{\pm}^*\pi_{\pm})z=z(\pi_{\pm}^*\pi_{\pm});&
(i)_2&\pi_{\pm}^*\pi_{\pm}>>0;\\
(ii)_1&(\pi_-^{\dag}\pi_+)z=z(\pi_-^{\dag}\pi_+);&
(ii)_2&P_{-}(\pi_{-}^{\dagger}\pi_{+})P_+=0;\\
(iii)&\Ran\pi_+\vee \Ran\pi_-={\mathcal H}\,,
\end{array}  \eqno{\rm{(Mod)}}
\]
where ${\mathfrak N}_{\pm},{\mathcal H}$ are separable Hilbert spaces; the notation
$A>>0$ means that $\exists\; c>0$ such that $\forall u\;(Au,u)\ge c(u,u)\,$; the
(nonorthogonal) projections $\,P_{\pm}\,$ are uniquely determined by conditions $\Ran
P_{\pm}=E^2(G_{\pm},{\mathfrak N}_{\pm})$ and $\Ker P_{\pm}=E^2(G_{\mp},{\mathfrak
N}_{\mp})$ ; the spaces $E^2(G_{\pm},{\mathfrak N}_{\pm})$ are Smirnov
spaces~\cite{Du} of vector valued functions with values in $\,{\mathfrak N}_{\pm}\,$
(since the curve $C$ is a Carleson curve, such projections exist); the operators
$\,\pi_{\pm}^{\dag}\,$ are adjoint to $\,\pi_{\pm}\,$ if we regard $\;\pi_{\pm}\colon
L^2(C,\Xi_{\pm})\to {\mathcal H}\;$ as operators acting from weighted spaces $L^2$
with operator valued weights $\,\Xi_{\pm}=\pi_{\pm}^*\pi_{\pm}\,$. In this
interpretation $\pi_{\pm}$ are isometries.

Note that, in our model, we strive to retain analyticity in both the domains $G_+$ and
$G_-$ with the aim to reserve possibility to exploit techniques typical for boundary
values problems (singular integral operators, the Riemman-Hilbert problem, the
stationary scattering theory, including the smooth methods of T.Kato). Thus we will
use both the Smirnov spaces $E^2(G_{\pm})$, which are analogues of the Hardy spaces
$H^2$ and $H^2_-$. The requirement of analyticity in both the domains conflicts with
orthogonality: in general, the decomposition $L^2(C)=E^2(G_{+})\dot+E^2(G_{-})$ is not
orthogonal. Note that the combination ``analyticity only in $G_+$ and orthogonality''
is a mainstream of development in the multiply-connected case starting from~\cite{AD}.
In this paper we sacrifice the orthogonality and therefore at this point we fork with
traditional way of generalization of Sz.-Nagy-Foia\c{s} theory~\cite{Ba, PF, F}.
Nevertheless, our requirements are also substantial and descends from applications
(see~\cite{T1,Ya,VYa}): in~\cite{T1} we studied the duality of spectral components for
trace class perturbations of a normal operator with spectrum on a curve; the
functional model from~\cite{Ya} goes back to the paper~\cite{VYa}, which is devoted to
spectral analysis of linear neutral functional differential equations.

The operator $\pi_{-}^{\dagger}\pi_{+}$ can be regarded as an analytic operator valued function
$(\pi_{-}^{\dagger}\pi_{+})(z),\;z\in G_+$. In this connection, we shall say that the set of
operator valued functions
\[
\Theta=(\pi_{-}^{\dagger}\pi_{+},\pi_{+}^*\pi_{+},\pi_{-}^*\pi_{-})\in S_{\Xi}\,.\eqno{\rm{(MtoC)}}
\]
is the \textit{characteristic function} for a model $\Pi$. Note also that the relation (MtoC)
defines the transformation $\Theta=\mathcal{F}_{cm}(\Pi)$.  Conversely, for a given $\Theta\in
S_{\Xi}\,$, it is possible to construct (up to unitary equivalence) a functional model $\Pi\in
\Mod$ such that $\Theta=(\pi_{-}^{\dagger}\pi_{+},\pi_{+}^*\pi_{+},\pi_{-}^*\pi_{-})\,$, i.e.,
there exists the inverse transformation $\mathcal{F}_{mc}:=\mathcal{F}_{cm}^{-1}$ (see
Prop.\ref{CtoM}).

At this moment we should look for a suitable generalization of conservative systems (=unitary
colligations). We define curved conservative systems in terms of the functional model. Let
$\Pi\in\Mod$. Define the model system
$\,\widehat{\Sigma}=\mathcal{F}_{sm}(\Pi):=(\widehat{T},\widehat{M},\widehat{N},\widehat{\Theta}_u,\widehat{\Xi};
\mathcal{K}_{\Theta},\mathfrak{N}_+,\mathfrak{N}_-)\,$,
where\\
\[
\begin{array}{ll}
{\widehat T}\in\mathcal{L}(\mathcal{K}_{\Theta})\,,\;\, & {\widehat T}f:=\mathcal{U}f-\pi_{+}{\widehat M}f\,;\\[2pt]
{\widehat M}\in\mathcal{L}(\mathcal{K}_{\Theta},{\mathfrak N_+})\,,\;\, & {\widehat
M}f:=\displaystyle \frac{1}{2\pi
i}\int_{C} (\pi_{+}^{\dag}f)(z)\,dz \,;\\
{\widehat N}\in\mathcal{L}({\mathfrak N_-},\mathcal{K}_{\Theta})\,,\;\, & {\widehat
N}n:=P_{\Theta}\pi_{-}n \,;\\[6pt]
\widehat{\Theta}_u \;\;\mbox{is the unitary "part" of} &
\widehat{\Theta}=(\pi_-^{\dag}\pi_+,\,\widehat{\Xi})\,;
\\[7pt]
\widehat{\Xi}:=(\pi_+^*\pi_+,\,\pi_-^*\pi_-)\,;
\end{array}
\eqno{\rm{(MtoS)}}
\]
$f\in \mathcal{K}_{\Theta}:=\Ran P_{\Theta}\,$,
$\;P_{\Theta}:=(I-\pi_{+}P_+\pi_{+}^{\dag})(I-\pi_{-}P_-\pi_{-}^{\dag})\,,\;n\in {\mathfrak
N_-}\,$, and the normal operator $\mathcal{U}$ with absolutely continuous spectrum lying on $C$ is
uniquely determined  by conditions $\,\mathcal{U}\pi_{\pm}=\pi_{\pm}z\,$. In the sequel, we shall
refer the operator $\,\widehat{T}\,$ as the model operator. The unitary "part" $\widehat{\Theta}_u$
is determined by the unitary constant part $\Theta_{u}^0$ from pure-unitary
decomposition~\cite{NF,Br} of Schur class function
$\Theta^{0}(w)=\Theta_{p}^0(w)\oplus\Theta_{u}^0,\;w\in\mathbb{D}$, where $\Theta^{0}(w)$ is the
lift of the (multiple valued character-automorphic) operator valued function
$\,(\chi_-\Theta^+\chi_+^{-1})(z)\,$ to the universal cover space~\cite{AD}; $\,\chi_{\pm}\,$ are
outer (character-automorphic) operator valued functions such that
$\,\chi_{\pm}^*\chi_{\pm}=\Xi_{\pm}\,$. Note also that the formulas (MtoS) define the
transformation $\,\widehat{\Sigma}=\mathcal{F}_{sm}(\Pi)\,$.

\vskip 2pt A coupling of operators and Hilbert spaces
$\Sigma=(T,M,N,\Theta_u,\Xi;H,\mathfrak{N},\mathfrak{M})$ is called a \textit{conservative curved
system} if there exists a functional model $\Pi$ with $\mathfrak{N}_+=\mathfrak{N}$ and
$\mathfrak{N}_-=\mathfrak{M}$, a Hilbert space $\mathcal{K}_u$, a normal operator $\widehat{T}_u\in
\mathcal{L}(\mathcal{K}_u)\,$, and an operator $X\in
\mathcal{L}(H,\mathcal{K}_{\Theta}\oplus\mathcal{K}_u)$ such that $\,\sigma(\widehat{T}_u)\subset
C$, $X^{-1}\in \mathcal{L}(\mathcal{K}_{\Theta}\oplus\mathcal{K}_u,H)$, and
\[
\Sigma\,=\,(T,M,N,\Theta_u,\Xi;H,\mathfrak{N},\mathfrak{M})
\;{\stackrel{X}{\sim}}\;(\widehat{\Sigma}\oplus\widehat{\Sigma}_u)\,,\eqno{\rm{(Sys)}}
\]
where $\,\widehat{\Sigma}=\mathcal{F}_{sm}(\Pi)\,$ and
$\,\widehat{\Sigma}_u=(\widehat{T}_{u},0,0,0;\mathcal{K}_u,\{0\},\{0\})\,$. We write
$\,\Sigma_1\,{\stackrel{X}{\sim}}\,\Sigma_2\,$ if
\[
XT_1=T_2 X\,,\quad M_1=M_2 X\,,\quad N_1 X=N_2\,,\quad
\Theta_{1u}=\Theta_{2u}\,,\quad\Xi_{1}=\Xi_{2}\,.
\]
The spaces $\mathcal{K}_{\Theta}$ and $\mathcal{K}_u$ play roles of the simple and ``unitary''
subspaces of the system $\widehat{\Sigma}\oplus\widehat{\Sigma}_u$, respectively. A curved
conservative system $\,\Sigma\,$ is called simple if
\[
\rho(T)\cap G_+\ne\emptyset\qquad \mbox{and} \qquad\bigcap_{z\in\rho(T)}\,\Ker M(T-z)^{-1}=\{0\}\,.
\]
In the case of unitary colligations this definition is equivalent to standard one~\cite{Br}
whenever $\,\rho(T)\cap \mathbb{D}\ne\emptyset\;$. Note that there appear some troubles if we
attempt to extend the standard definition (simple subspace = controllable subspace  $\vee$
observable subspace) straightforwardly.

In the case when $G_+=\mathbb{D}$ and $\Xi_{\pm}\equiv I$, for a conservative curved
system $\,\Sigma=(T,M,N,\Theta_u,\Xi;H,\mathfrak{N},\mathfrak{M})\,$, we can consider
the block-matrix ${\mathfrak A}=\left(\begin{array}{cc}T&N\\M&L\end{array}\right)$\,,
where $\,L=\Theta^+(0)^*\,$. It is readily shown that ${\mathfrak A}$ is a unitary
colligation and $\,\Theta^+(z)=L^*+zN^*(I-zT^*)^{-1}M^*\,,\,|z|<1\,$, i.e.
$\,\Sigma\,$ is a conservative system and $\,\Theta^+\,$ is the Sz.-Nagy-Foia\c{s}
characteristic function. In the  case of simple-connected domains we lose the property
of unitarity for the matrix ${\mathfrak A}$ but we can regard a system
$\,\Sigma=(T,M,N)\,$ as the result of certain transformation (deformation) of a
unitary colligation (=conservative system)~\cite{T1,T2}. Another reason to call our
systems ``curved conservative'' is the fact that the characteristic function of such a
system is a weighted Schur function.

Thus, we have defined the notion of conservative curved system. Note that linear similarity
(instead of unitary equivalence for unitary colligations) is a natural kind of equivalence for
conservative curved systems and duality is a substitute for orthogonality. The following diagram
shows relationships between models, characteristic functions, and conservative curved systems
\[
\xymatrix{\Cfn \ar@<0.4ex>[r]^{\mathcal{F}_{mc}} & \Mod \ar@<0.4ex>[l]^{\mathcal{F}_{cm}}
\ar[r]^{\mathcal{F}_{sm}} & \Sys }\,.\eqno{\rm{(dgr)}}
\]
As we can now see, characteristic functions and conservative curved systems are not on
equal terms: first of them plays leading role because the definition of conservative
curved system depends on the functional model, which, in turn, is uniquely determined
by the characteristic function. But, surprisingly, the conservative curved systems is
a comparatively autonomous notion (i.e., though we define such systems in terms of the
functional model, many properties and operations with conservative curved systems can
be formulated intrinsically and  do not refer explicitly to the functional model) and
one of the aims of this paper is to ``measure'' a degree of this autonomy with the
point of view of the correspondence ``factorizations of characteristic function
$\leftrightarrow$ invariant subspaces''.

If we are going to follow the way described above for conservative systems, we need to
introduce transfer functions. For a curved conservative system $\Sigma$, we define the
\textit{transfer function}
\[
\Upsilon=(\Upsilon(z),\,\Theta_u,\,\Xi) \;\mbox{,
where}\quad\Upsilon(z):=M(T-z)^{-1}N\,.\eqno{\rm{(Tfn)+(StoT)}}
\]
The formula (StoT) defines also the transformation
$\Upsilon=\mathcal{F}_{ts}(\Sigma)$. Then, using the functional model, the
transformation $\,\mathcal{F}_{tc}=\mathcal{F}_{ts}\circ\mathcal{F}_{sc}\,$ can be
computed as
\[
\Upsilon(z)=\left\{%
\begin{array}{ll}
    \Theta^-_+(z)-\Theta^+(z)^{-1}\,, & z\in G_+\cap\rho(T)\,; \\
    -\Theta^-_-(z)\,, & z\in G_-\,.
\end{array}%
\right.\eqno{\rm{(CtoT)}}
\]
In this connection, note that the spectrum of a model operator coincides with the
spectrum of a characteristic function, i.e., $\,z\in
G_+\cap\rho(T)\,\Leftrightarrow\,\exists\;\Theta^+(z)^{-1}\,$. The operator valued
functions $\Theta^-_{\pm}(z)$ are defined by the formulas
\[
\begin{array}{l}
    \Theta^-_{\pm}(z)\,n:=(P_{\pm}\Theta^-n)(z),\;\;z\in G_{\pm},\;\;n\in \mathfrak{N}_-\,; \\[3pt]
    \Theta^-(\zeta):=(\pi_+^{\dag}\pi_-)(\zeta)=\Xi_+(\zeta)^{-1}\Theta^+(\zeta)^*\Xi_-(\zeta),\;\;\zeta\in C\,.
\end{array}
\]
In the case when $\,G_+=\mathbb{D}\,$ and $\,\Xi_{\pm}\equiv I\,$ we get
$\,\Theta^-(\zeta)=\Theta^+(1/{\bar\zeta})^*\,,|\zeta|=1$ and therefore,\,
$\,\Theta_{+}^-(z)=\Theta^+(0)^*\,,\;|z|<1\,$;\,\, $\,\Theta_{-}^-(z)=\Theta^+(1/{\bar
z})^*-\Theta^+(0)^*,\;|z|>1\,$.

\vskip 2pt Thus we arrive at the complete diagram
\[
\xymatrix{\Mod \ar@<0.4ex>[r]^{\mathcal{F}_{cm}} \ar[d]_{\mathcal{F}_{sm}}
& \Cfn \ar@<0.4ex>[l]^{\mathcal{F}_{mc}} \ar[d]^{\mathcal{F}_{tc}}\\
\Sys \ar[r]_{\mathcal{F}_{ts}} & \Tfn}   \eqno{\raisebox{-20pt}{\rm{(Dgr)}}}
\]
Unfortunately, we have obtained almost nothing for our purpose: to study the
correspondence ``factorizations $\leftrightarrow$ invariant subspaces''. The main
difficulty is to invert the arrows $\mathcal{F}_{tc}$ and $\mathcal{F}_{ts}$. In the
case of the unit circle the transfer function  can be calculated as
$\Upsilon(z)=\Theta^+(0)^*-\Theta^+(1/{\bar z})^*,\;|z|>1$ and, conversely, one can
easily recover the characteristic function $\Theta^+(z)$ from the transfer function
$\Upsilon(z)$ (see~\cite{T3} for this case and for the case of simple connected
domains). But, in general, especially for multiply-connected domains, the latter is a
considerable problem. Note that the condition $\,\Upsilon(z)=M(T-z)^{-1}N\in N(G_+\cup
G_-)\,$ (that is, the transfer function $\,\Upsilon(z)\,$ is an operator valued
function of Nevanlinna class: $\Upsilon(z)=1/\delta(z)\,\Omega(z)\,$, where $\delta\in
H^{\infty}(G_+\cup G_-)$ and $\Omega\in H^{\infty}(G_+\cup
G_-,\mathcal{L}(\mathfrak{N}_{-},\mathfrak{N}_{+}))$) is sufficient for uniqueness of
characteristic function and there is a procedure recovering the characteristic
function from a given transfer function. Moreover, under this assumption it is
possible to give intrinsic description for conservative curved systems. Note that we
reap the benefit of functional model when we are able to determine that some set of
operators $(T,M,N)$ is a conservative curved system~\cite{T1,Ya,T2}. The author plans
to address these problems elsewhere.

Thus we distinguish notions of characteristic and transfer function and there are no simple enough
(and suitable in the study of factorizations) relationships between them. These circumstances
dictate that we have to use only the partial diagram (dgr) and to ignore other objects and
transformations related to transfer functions from the complete diagram (Dgr). Note also that we
study the correspondence ``factorizations of \textit{characteristic} function $\leftrightarrow$
invariant subspaces of operator $T$'' in contrast to the correspondence studied in~\cite{BGK}:
``factorizations of \textit{transfer} function $\leftrightarrow$ invariant subspaces''. At this
point we fork with~\cite{BGK}.

The paper is organized as follows. In Section~1 we deal with the fragment
$\xymatrix{\Cfn \ar@<0.4ex>[r]^{\mathcal{F}_{mc}} & \Mod
\ar@<0.4ex>[l]^{\mathcal{F}_{cm}}}$: in the context of the functional model we develop
the constructions corresponding to factorizations of characteristic functions. If we
restrict ourselves to regular factorizations, we can keep on to exploit the functional
model $\Mod$. But to handle arbitrary factorizations and to obtain a pertinent
definition of the product of conservative curved systems we need some generalization
of $\Mod$. Moreover, the order relation
$\,\Theta_{2}\Theta_{1}\prec\Theta_{2}'\Theta_{1}'\,$ implies the factorizations like
$\,\Theta_{2}'\theta\Theta_{1}\,$ and therefore we need a functional model suited to
handle factorizations of characteristic function with three or more multipliers. With
this aim we introduce the notion of n-model $\Mod_n$ and extend the transformations
$\mathcal{F}_{mc}$ and $\mathcal{F}_{cm}$ to this context. In the rest part of the
section we study geometric properties of n-models in depth and do this mainly because
they form a solid foundation for our definition of the product of curved conservative
systems in the next section.

At this moment it is unclear how to define the product of conservative curved systems.
As a first approximation we can consider the following construction. Let
$\,\Sigma_1\sim\widehat{\Sigma}_1=\mathcal{F}_{sm}(\mathcal{F}_{mc}(\Theta_1))\,$ and
$\,\Sigma_2\sim\widehat{\Sigma}_2=\mathcal{F}_{sm}(\mathcal{F}_{mc}(\Theta_2))\,$.
Then the candidate for their product is
$\,\widehat{\Sigma}_{21}=\mathcal{F}_{sm}(\mathcal{F}_{mc}(\Theta_2\Theta_1))\,$,
where $\,\mathcal{F}_{mc}(\Theta_2\Theta_1)\,$ is 3-model corresponding to the
factorization $\,\Theta_{21}=\Theta_2\cdot\Theta_1\,$. Our aim is to define the
product $\,\Sigma_{2}\cdot\Sigma_{1}\,$ by explicit formulas without referring to the
functional model. In Section~2 we suggest such a definition and study basic properties
of it. The main one among those properties is the property that the product of systems
$\,\Sigma_{2}\cdot\Sigma_{1}\,$ is a conservative curved system too (Theorem A). The
geometry properties of n-model established in Section 1 play crucial role in our
reasoning.

In Section~3 we establish the correspondence between two notions of regularity. The first of them
is the regularity of the product of conservative curved systems $\,\Sigma_{2}\cdot\Sigma_{1}\,$,
the second one is the notion of regular factorization of operator valued functions~\cite{NF,Br},
which we extend to the weighted Schur classes. We obtain this correspondence indirectly: introduce
the notion of regularity for models and establish separately the correspondences
$\,\Cfn\leftrightarrow\Mod\,$ and $\,\Sys\leftrightarrow\Mod\,$.

In Section~4 we study the transformation $\,\mathcal{F}_{ic}$ defined therein, which
takes a factorization $\Theta_2\Theta_1$ of characteristic function to the invariant
subspace of the model operator of the system
$\,\widehat{\Sigma}_{21}=\mathcal{F}_{sm}(\mathcal{F}_{mc}(\Theta_2\Theta_1))\,$. We
show that this mapping is surjective. Combining this property of $\,\mathcal{F}_{ic}$
with the criterion of regularity from Section 3, we establish the main result of the
paper: \textit{there is an order preserving one-to-one correspondence between regular
factorizations of a characteristic function and invariant subspaces of the resolvent
of the corresponding model operator}. In conclusion we translate results obtained for
model operators into the language of conservative curved systems.

Note that the multiply connected domain specific appears essentially only in the proof
of Prop.~\ref{u_chain}. So, at first a reader can study the paper assuming that the
domain $G_+$ is simple connected. On the other hand, the multiply connected specific
influences on the choice of other our proofs: note that, for simple connected domains,
some of them can be reduced to the case of the unit disk (see, e.g.,~\cite{T2,T3}).

\section{Geometric properties of n-model}

We start with the definition of an \textit{n-characteristic function}, which formalizes products of
weighted Schur class functions like the following
$\,\theta_{n-1}\cdot\ldots\cdot\theta_{2}\theta_{1}\,$ : in fact, we merely rearrange them
$\,\Theta_{ij}:=\theta_{i-1}\cdot\ldots\cdot\theta_{j}\,$.
\begin{defin}
Let $\,\Xi_k,\;k=\overline{1,n}\;$ be operator valued weights such that
$\Xi_{k},\Xi_{k}^{-1}\in L^\infty(C,{\mathcal L}({\mathfrak N}_{k}))$,
$\Xi_{k}(\zeta)\ge 0,\;\zeta\in C\;$. A set of analytic in $G_+$ operator valued
functions $\,\Theta=\{\Theta_{ij} : i\ge j\}$ is called an n-characteristic function
if $\;\Theta_{ij}\in S_{\Xi}$ with weights $\Xi=(\Xi_{i},\Xi_{j})$ and $\;\forall \;
i\ge j \ge k \;\; \Theta_{ik}=\Theta_{ij}\Theta_{jk}\,$.
\end{defin}
\noindent We assume that $\,\Theta_{kk}:=I\,$ and denote by $\,\rm{Cfn}_n\,$ the class
of all n-characteristic functions. In the sequel, we shall usually identify
3-characteristic function with the factorization of Schur class function
$\,\theta=\theta_2\cdot\theta_1\,$, where $\,\theta=(\Theta_{31},\Xi_1,\Xi_3)\,$,
$\,\theta_1=(\Theta_{21},\Xi_1,\Xi_2)\,$, $\,\theta_2=(\Theta_{32},\Xi_2,\Xi_3)\,$. It
is clear how to define the product of n-characteristic functions
$\,\Theta=\Theta''\cdot\Theta'$ : assuming that $\,\Xi_{n'}'=\Xi_{1}''$, we need only
to renumber multipliers, for instance,
$\,\Theta_{ij}=\Theta_{i-n'+1,1}''\Theta_{n'j}'\,$, $\,i\ge n'\ge j\,$.

\vskip 3pt In the context of functional models a corresponding notion is the notion of
\textit{n-functional model}.
\begin{defin}
An n-tuple $\,\Pi=(\pi_1,\dots,\,\pi_n)$ of operators $\pi_{k}\in {\mathcal L}(L^2(C,{\mathfrak
N_{k}}),{\mathcal H})$ such that
\[
\begin{array}{clll}
(i) &   \forall\; k & \;(\pi_{k}^*\pi_{k})z=z(\pi_{k}^*\pi_{k}); &\; \pi_{k}^*\pi_{k}>>0;\\[2pt]
(ii) & \forall\; j \ge k & \;(\pi_j^{\dag}\pi_k)z=z(\pi_j^{\dag}\pi_k);&
\;P_{-}(\pi_{j}^{\dagger}\pi_{k})P_+=0;\\[2pt]
(iii) & \forall\; i \ge j \ge k & \;\;\;\pi_i^{\dag}\pi_k=\pi_i^{\dag}\pi_j\pi_j^{\dag}\pi_k;\\[2pt]
(iv) &  & \;\;\;\mathcal{H}_{\pi_n\vee\dots\vee\pi_1}=\mathcal{H}\\[2pt]
\end{array}  \eqno{\rm{(Mod_n)}}
\]
is called an n-model.
\end{defin}
\noindent Here $\,\mathcal{H}_{\pi_n\vee\dots\vee\pi_1}:=\vee_{k=1}^n\Ran\pi_k\,$. The
definition is an extension of the definition (Mod): namely, $\Mod=\Mod_2$. It is
readily seen that $\Theta=\{\pi_i^{\dag}\pi_j\}_{i \ge j}$ is an n-characteristic
function with weights $\Xi_k=\pi_k^*\pi_k$ and therefore we have defined the
transformation $\,\mathcal{F}_{cm} : \Mod_n \to \Cfn_n \,$. The existence of the
``inverse'' transformation $\,\mathcal{F}_{mc}\,$
 follows from
\begin{prop}\label{CtoM}
Suppose $\Theta\in\Cfn_n$. Then $\;\exists\,\Pi\in\Mod_n\,$ such that
$\,\Theta=\mathcal{F}_{cm}(\Pi)\,$. If also $\,\Theta=\mathcal{F}_{cm}(\Pi')\,$, then there exists
an unitary operator $\,X :
\mathcal{H}_{\pi_n\vee\dots\vee\pi_1}\to\mathcal{H}_{\pi_n\vee\dots\vee\pi_1}'\,$ such that
$\,\pi_k'=X\pi_k\,$.
\end{prop}
\myproof We put $\;\mathcal{H}=\oplus_{k=1}^{n}\mathcal{H}^{\Delta}_k\,$, where
$\;\mathcal{H}^{\Delta}_k=\clos\Delta_{kk+1k}L^2(C,{\mathfrak
N_{k}}),\;k=\overline{1,n-1},\;$\\$\mathcal{H}^{\Delta}_n=L^2(C,{\mathfrak N_{n}})\,$,
$\,\Delta_{kk+1k}:=(I-\Theta^{\dag}_{k+1k}\Theta_{k+1k})^{1/2}$, and $\Theta^{\dag}_{k+1k}$ is
adjoint to the operator $\;\Theta_{k+1k} : L^2(C,\Xi_{k})\to L^2(C,\Xi_{k+1})$. Let
$\,\nu_k,\;k=\overline{1,n}$ be the operators of embedding of $\mathcal{H}^{\Delta}_k$ into
$\mathcal{H}\,$ and
\[
\pi_n:=\nu_n\,,\quad \pi_k:=\pi_{k+1}\Theta_{k+1k}+\nu_k\Delta_{kk+1k}\,,\quad
k=\overline{1,n-1}\,.
\]
It can easily be calculated that
\[
\pi_k=\nu_n\Theta_{nk}+\nu_{n-1}\Delta_{n-1nn-1}\Theta_{n-1k}+
\ldots+\nu_{j}\Delta_{jj+1j}\Theta_{jk}+\ldots+\nu_k\Delta_{kk+1k}\,.
\]
From this identity we get $\,\pi_i^{\dag}\pi_j=\Theta_{ij}\,,\;i \ge j \,$.

\noindent The existence and unitary property of $\,X\,$ follows from the identity
\[
\begin{array}{lll}
||\,\pi_1 u_1+\ldots+\pi_n u_n\,||^2 &=& \sum\limits_{i,j=1}^n(\pi_j^{\dag}\pi_i u_i,\,
u_j)_{L^2(C,\Xi_{j})}\;=
\\[3pt]
\sum\limits_{i,j=1}^n({\pi'}_{j}^{\dag}\pi_i' u_i,\, u_j)_{L^2(C,\Xi_{j})}\ &=& ||\,\pi_1'
u_1+\ldots+\pi_n' u_n\,||^2\,.\qquad\qquad\square
\end{array}
\]
The construction of Prop.~\ref{CtoM} is simplified if all functions $\Theta_{ij}$ are
two-sided $\Xi$-inner. In this case $\,\mathcal{H}=L^2(C)\,$ and
$\,\pi_k=\Theta_{nk}\,$.

We can consider an equivalence relation $\,\sim\,$ in $\Mod_n\,$ : $\,\Pi\sim\Pi'\,$
if there exists an unitary operator $\,X :
\mathcal{H}_{\pi_n\vee\dots\vee\pi_1}\to\mathcal{H}_{\pi_n\vee\dots\vee\pi_1}'\,$ such
that $\,\pi_k'=X\pi_k\,$. It is clear that the transformation $\,\mathcal{F}_{cm}\,$
induces a transformation $\;\mathcal{F}_{cm}^{\sim} : \Mod_n^{\sim} \to \Cfn_n \;$
such that $\,\mathcal{F}_{cm}^{\sim}(\Pi^{\sim})=\mathcal{F}_{cm}(\Pi)\,$,
$\,\Pi\in\Pi^{\sim}\,$. By Prop.\ref{CtoM}, there exists the inverse transformation
$\;\mathcal{F}_{mc}^{\sim} : \Cfn_n \to \Mod_n^{\sim} \;$. But, in the sequel, we
shall usually ignore this equivalence relation and use merely the transformations
$\;\mathcal{F}_{cm}\,$ and $\,\mathcal{F}_{mc}\,$.

The product of n-models $\,\Pi',\Pi''\,$ with the only restriction
$\,\pi_{n}'^*\pi_{n}'=\pi_{1}''^*\pi_{1}''$ is defined (up to unitary equivalence) as
$\,\Pi=\Pi''\cdot\Pi':=\mathcal{F}_{mc}(\mathcal{F}_{cm}(\Pi'')\cdot\mathcal{F}_{cm}(\Pi'))\,$.

Using the construction of Prop.\ref{CtoM}, we can uniquely determine the normal
operator $\mathcal{U}=XzX^{-1}\in \mathcal{L}(\mathcal{H}_{\pi_n\vee\dots\vee\pi_1})$
with absolutely continuous spectrum $\sigma(\mathcal{U})\subset C\,$ such that
$\,\mathcal{U}\pi_k=\pi_kz\,$, where $\,X :
\hat{\mathcal{H}}_{\pi_n\vee\dots\vee\pi_1}\to\mathcal{H}_{\pi_n\vee\dots\vee\pi_1}\,$
is an unitary operator such that $\,\pi_k=X\hat{\pi}_k\,$; the operators
$\,\hat{\pi}_k\,$ are constructed for n-characteristic function
$\,\Theta=\mathcal{F}_{cm}(\Pi)\,$ as in Prop.\ref{CtoM}.

Taking into account the existence of a such operator $\mathcal{U}$, note that
$\,\mathcal{F}_{sc}(\Theta)=\mathcal{F}_{sc}(\Theta_{n1})\oplus\widehat{\Sigma}_u\,$, where the
system $\widehat{\Sigma}_u=(\widehat{T}_u,0,0,0)\,$ is a ``purely normal'' system with the normal
operator $\;\widehat{T}_u=\,\mathcal{U}\,|\,(\mathcal{H}_{\pi_n\vee\dots\vee\pi_1}\ominus
\mathcal{H}_{\pi_n\vee\pi_1})\,$, $\;\sigma(\widehat{T}_u)\subset C\,$.\\

Let $\Pi\in \Mod_n$. Now we define our building bricks:  orthoprojections
$\,P_{\pi_i\vee\dots\vee\pi_j}\,$ onto $\,\mathcal{H}_{\pi_i\vee\dots\vee\pi_j}$ and projections
$\,q_{i\pm}:=\pi_iP_{\pm}\pi_i^{\dag}\,$.
\begin{lemma}\label{lemm_mod1}
For $\;i\ge j\ge k\ge l\ge m\,$ \,\,\,1) $\;q_{i-}q_{j+}=0$;
\,\,\,2) $\;q_{i+}+q_{i-}=\pi_i\pi_i^{\dag}=P_{\pi_i}$;\\
\,3) $\;P_{\pi_i\vee\dots\vee\pi_j}(I-\pi_k\pi_k^{\dag})P_{\pi_l\vee\dots\vee\pi_m}=0$\,; \,4)
$\;P_{\pi_l\vee\dots\vee\pi_m}(I-\pi_k\pi_k^{\dag})P_{\pi_i\vee\dots\vee\pi_j}=0$\,.
\end{lemma}
\proof Statement 1) is a direct consequences of $(ii)$ in (Mod$_n$). Statement 2) is obvious.
Statement 3) is equivalent to the relation
\[
\forall f,g\in \mathcal{H}\qquad
((I-\pi_k\pi_k^{\dag})P_{\pi_l\vee\dots\vee\pi_m}f,P_{\pi_i\vee\dots\vee\pi_j}g)=0\,.
\]
The latter can be rewritten in the form
\[
((I-\pi_k\pi_k^{\dag})\pi_{l'} u,\pi_{i'} v)=0\;,\quad j\le i'\le i\,,\quad m\le l'\le l\,
\]
and is true because of $(iii)$ in (Mod$_n$). Statement 4) can be obtained from Statement 3) by
conjugation. \proofend

\noindent We also define the projections
\[
P_{(ij)}:=P_{\pi_i\vee\dots\vee\pi_j}(I-q_{j+})(I-q_{i-}),\;\; i\ge j\,\,.
\]
It is easily shown that $\,P_{(ii)}=0\,$ and
\[
P_{(ij)}=(I-q_{j+})P_{\pi_i\vee\dots\vee\pi_j}(I-q_{i-})=(I-q_{j+})(I-q_{i-})P_{\pi_i\vee\dots\vee\pi_j}\,.
\]
Indeed, $\,P_{\pi_i\vee\dots\vee\pi_j}-P_{\pi_i}\,$ is orthoprojection onto
$\,\mathcal{H}_{\pi_i\vee\dots\vee\pi_j}\ominus\mathcal{H}_{\pi_i}\,$,
$\,P_{\pi_i}q_{i\pm}=q_{i\pm} P_{\pi_i}\,$ and
$\,(P_{\pi_i\vee\dots\vee\pi_j}-P_{\pi_i})q_{i\pm}=q_{i\pm}(P_{\pi_i\vee\dots\vee\pi_j}-P_{\pi_i})=0\,$.
The same is hold for $\,\pi_j\,$. Then,
\[
\begin{array}{ll}
P_{(ij)}^2&=P_{\pi_i\vee\dots\vee\pi_j}(I-q_{j+})(I-q_{i-})(I-q_{j+})(I-q_{i-})P_{\pi_i\vee\dots\vee\pi_j}\\[2pt]
&=P_{\pi_i\vee\dots\vee\pi_j}[(I-q_{j+})(I-q_{j+})(I-q_{i-})-
(I-q_{j+})q_{i-}(I-q_{i-})]P_{\pi_i\vee\dots\vee\pi_j}\\[2pt]
&=P_{\pi_i\vee\dots\vee\pi_j}[(I-q_{j+})(I-q_{i-})]P_{\pi_i\vee\dots\vee\pi_j}=P_{(ij)}\,.
\end{array}
\]
Note also that $\;P_{(ij)}=\Theta_{nj}P_-\Theta_{ij}^{-1}P_+\Theta_{ni}^{-1}\,$
whenever all functions $\Theta_{ij}$ are two-sided $\Xi$-inner (recall that then we
can choose $\,\pi_k=\Theta_{nk}\,$).
\begin{lemma}\label{lemm_mod2}
For $\;i\ge j\ge k\ge l\,$, one has
\[
\begin{array}{lll}
1) \;P_{(ij)}q_{k+}=0; & 2) \;q_{i-}P_{(jk)}=0; & 3) \;P_{(ij)}P_{(kl)}=0; \\[4pt]
4) \;P_{(ik)}P_{(jk)}=P_{(jk)}; & 5) \;P_{(ij)}P_{(ik)}=P_{(ij)}; & 6)
\;P_{(jk)}P_{(ij)}=0\,.\qquad\qquad\qquad\qquad\qquad\qquad
\end{array}
\]
\end{lemma}
\proof Using Lemma~\ref{lemm_mod1}, we have\\[2pt]
1)
\[
\begin{array}{ll}
P_{(ij)}q_{k+}&=P_{\pi_i\vee\dots\vee\pi_j}(I-q_{j+})(I-q_{i-})q_{k+}=
P_{\pi_i\vee\dots\vee\pi_j}(I-q_{j+})q_{k+}\\[4pt]
&=P_{\pi_i\vee\dots\vee\pi_j}[(I-\pi_j\pi_j^\dag)+q_{j-}]q_{k+}
=P_{\pi_i\vee\dots\vee\pi_j}(I-\pi_j\pi_j^\dag)\pi_kP_{+}\pi_k^{\dag}=0\,;
\end{array}
\]
2)
\[
\begin{array}{ll}
q_{i-}P_{(jk)}&=q_{i-}(I-q_{k+})(I-q_{j-})P_{\pi_j\vee\dots\vee\pi_k}=
q_{i-}(I-q_{j-})P_{\pi_j\vee\dots\vee\pi_k}\\[4pt]
&=q_{i-}[(I-\pi_j\pi_j^\dag)+q_{j+}]P_{\pi_j\vee\dots\vee\pi_k} =
\pi_iP_{-}\pi_i^{\dag}(I-\pi_j\pi_j^\dag)P_{\pi_j\vee\dots\vee\pi_k}=0\,;
\end{array}
\]
3)
\[
\begin{array}{ll}
P_{(ij)}P_{(kl)}&=P_{(ij)}(I-q_{l+})(I-q_{k-})P_{\pi_k\vee\dots\vee\pi_l}
=P_{(ij)}(I-q_{k-})P_{\pi_k\vee\dots\vee\pi_l}\qquad\qquad\\[4pt]
&=P_{(ij)}[(I-\pi_k\pi_k^\dag)+q_{k+}]P_{\pi_k\vee\dots\vee\pi_l}\\[4pt]
&=(I-q_{j+})(I-q_{i-})P_{\pi_i\vee\dots\vee\pi_j}(I-\pi_k\pi_k^\dag)P_{\pi_k\vee\dots\vee\pi_l}=0\,;
\end{array}
\]
4)
\[
\begin{array}{ll}
P_{(ik)}P_{(jk)}&=P_{\pi_i\vee\dots\vee\pi_k}(I-q_{k+})(I-q_{i-})P_{(jk)}
=P_{\pi_i\vee\dots\vee\pi_k}(I-q_{k+})P_{(jk)}=P_{(jk)}\,;
\end{array}
\]
5)
\[
\begin{array}{ll}
P_{(ij)}P_{(ik)}&=P_{(ij)}(I-q_{k+})(I-q_{i-})P_{\pi_i\vee\dots\vee\pi_k}
=P_{(ij)}(I-q_{i-})P_{\pi_i\vee\dots\vee\pi_k}=P_{(ij)}\,;
\end{array}
\]
6)
\[
\begin{array}{ll}
P_{(jk)}P_{(ij)}&=
P_{\pi_j\vee\dots\vee\pi_k}(I-q_{k+})(I-q_{j-})(I-q_{j+})(I-q_{i-})P_{\pi_i\vee\dots\vee\pi_j}\qquad\qquad\\[4pt]
&=P_{\pi_j\vee\dots\vee\pi_k}(I-q_{k+})(I-\pi_j\pi_j^\dag)(I-q_{i-})P_{\pi_i\vee\dots\vee\pi_j}\\[4pt]
&=(I-q_{k+})P_{\pi_j\vee\dots\vee\pi_k}(I-\pi_j\pi_j^\dag)P_{\pi_i\vee\dots\vee\pi_j}(I-q_{i-})=0\,.
\end{array}
\]
\proofend

\noindent We also define the subspaces
\[
\begin{array}{ll}
\mathcal{K}_{(ij)}:=\Ran P_{(ij)}\,,\quad & \mathcal{H}_{ij}:=\mathcal{H}_{\pi_i\vee\dots\vee\pi_j}\\[6pt]
\mathcal{H}_{ij+}:=\mathcal{H}_{ij}\cap\Ker q_{i-}\,,\qquad & \mathcal{D}_{j+}:=\Ran
q_{j+}\,.\qquad\quad
\end{array}
\]
\vskip 4pt \noindent It is easy to prove that $\;\mathcal{H}_{ij+}\cap \Ker
P_{(ij)}=\mathcal{D}_{j+}\,$. Indeed, let $\,f\in \mathcal{H}_{ij+}\cap \Ker P_{(ij)}\,$. Then
\[
\begin{array}{ll}
f&=(I-P_{(ij)})f=f-(I-q_{j+})(I-q_{i-})P_{\pi_i\vee\dots\vee\pi_j}f\\[4pt]
&=f-(I-q_{j+})(I-q_{i-})f=f-(I-q_{j+})f=q_{j+}f\in \mathcal{D}_{j+}\,.
\end{array}
\]
Conversely, let $\,f\in \mathcal{D}_{j+}\,$. Then $\,f=q_{j+}f\in \mathcal{H}_{ij}\,$ and therefore
$\,q_{i-}f=q_{i-}q_{j+}f-0\,$, that is, $\,f\in \mathcal{H}_{ij+}\,$. Hence we have
\[
\begin{array}{ll}
P_{(ij)}f&=(I-q_{j+})(I-q_{i-})P_{\pi_i\vee\dots\vee\pi_j}f =(I-q_{j+})(I-q_{i-})f=(I-q_{j+})f=0
\end{array}
\]
and $\,f\in \mathcal{H}_{ij+}\cap \Ker P_{(ij)}\,$.

\vskip 4pt Translating the assertions of the above lemmas into the language of
geometry, we obtain
\[
\mathcal{K}_{(ij)}\subset\mathcal{H}_{ij+}\,,\qquad
\mathcal{K}_{(jk)}\subset\mathcal{K}_{(ik)}\,,\qquad \mathcal{H}_{jk+}\subset
\mathcal{H}_{il+}\,,\qquad i\ge j\ge k\ge l\,.
\\[2pt]
\]
Indeed, let $\,f\in \mathcal{K}_{(ij)}\,$. Then $\,f=P_{(ij)}f\in
\mathcal{H}_{ij}\,$ and
$\,g_{i-}f=g_{i-}(I-q_{j+})(I-q_{i-})f=g_{i-}(I-q_{i-})f=0\,$
$\,\Rightarrow\, f\in \Ker q_{i-}$. The inclusion
$\,\mathcal{K}_{(jk)}\subset\mathcal{K}_{(ik)}\,$ is a
straightforward consequence of Lemma~\ref{lemm_mod2}(4). Let
$\,f\in \mathcal{H}_{jk+}\,$. Then
\[
g_{i-}f=g_{i-}(I-g_{j-})P_{\pi_j\vee\dots\vee\pi_k}f=
g_{i-}[(I-\pi_j\pi_j^\dag)+g_{j+}]P_{\pi_j\vee\dots\vee\pi_k}f=0\,.
\]
and therefore $\,f\in \mathcal{H}_{il+}\,$.

\vskip 4pt Let $\;1=m_1\le\ldots\le m_i\le \ldots\le m_N=n\,$. We define the operators\\
\[
\begin{array}{lcl}
P_{[m_i m_j]} &:\,=& P_{(m_{j+1} m_{j})}(I-P_{(m_{j+2} m_{j+1})})\ldots(I-P_{(m_{i} m_{i-1})})\\[3pt]
 &+& P_{(m_{j+2} m_{j+1})}(I-P_{(m_{j+3} m_{j+2})})\ldots(I-P_{(m_{i} m_{i-1})})\\[3pt]
 &+& \ldots+P_{(m_{i-1} m_{i-2})}(I-P_{(m_{i} m_{i-1})})+P_{(m_{i} m_{i-1})}\,,\qquad\;\;i\ge j\,.
\end{array}\\[3pt]
\]
\vskip 2pt \noindent  Note that our notation is ambiguous: the
projection $\,P_{[m_i m_j]}\,$ depends on the whole chain
$\;m_j\le\ldots\le m_i\,$ but not only on two numbers $\,m_j\,$
and $\,m_i\,$. The following properties of operators $\,P_{[m_i
m_j]}\,$ are straightforward consequences of
Lemma~\ref{lemm_mod2}.
\begin{prop}\label{proj_prop}
For $\;i\ge j\ge k\ge l\,$, \,\,1) $\;P_{[m_i m_j]}q_{m_k+}=0$; \,\,2) $\;q_{m_i-}P_{[m_j
m_k]}=0$;\\
\,\,3) $\,P_{[m_i m_j]}P_{[m_k m_l]}=0$\,.
\end{prop}
\noindent Further, since $\, I-P_{[m_i m_j]}= (I-P_{(m_{j+1} m_{j})})(I-P_{(m_{j+2}
m_{j+1})})\ldots(I-P_{(m_{i} m_{i-1})}) \,$, we get the following recursion relation
\[
P_{[m_im_k]}=P_{[m_jm_k]}(I-P_{[m_im_j]})+P_{[m_im_j]}\,,\qquad i\ge j\ge k\,.
\]
Since $\,P_{[m_{j+1} m_j]}=P_{(m_{j+1} m_j)}$, we obtain by induction that the
operator $\,P_{[m_i m_j]}\,$ is a projection and
\[
\mathcal{K}_{[m_i m_j]}=\mathcal{K}_{(m_i
m_{i-1})}\dot+\ldots\dot+\mathcal{K}_{(m_{j+1} m_j)}\,,
\]
where $\;\mathcal{K}_{[m_i m_j]}:=\Ran P_{[m_i m_j]}\,$. We use the notation
$\,H=H'\dot+ H''\,$ if there exists a projection $P'$ such that $\,H'=\Ran
P',\,H''=\Ker P'\,$. Besides we have
\[
\mathcal{K}_{[m_i m_j]}\subset\mathcal{H}_{m_i m_j+}\,,\qquad
\mathcal{D}_{m_j+}\subset \Ker P_{[m_i m_j]}\,.\\[4pt]
\]
The first inclusion follows straightforwardly from Prop.~\ref{proj_prop}(2). The
second one is a consequence of Prop.~\ref{proj_prop}(1).\\[-8pt]

Though we systematically strive to deal only with nonorthogonal projections
$\,q_{i+}\,$, sometimes we have to employ their orthogonal counterparts $\,q_{i+}'\,$.
By~\cite{AD}, there exists an isometries $\,\pi_i'\in {\mathcal L}(L^2(C,{\mathfrak
N_{i}}),{\mathcal H}) \,$ such that $\,\Ran \pi_i'=\Ran \pi_i\,$ and $\,\Ran
q_{i+}=\pi_i' E_{\alpha}^2(G_+,{\mathfrak N_{i}})\,$, where
$\,E_{\alpha}^2(G_+,{\mathfrak N_{i}})\,$ is the Smirnov space of
character-automorphic functions (see~\cite{AD,Ba,F} for the definition). Let
$\,q_{i+}'=\pi_i'P_{+}'\pi_i'^{*}\,$. Then $\,\Ran q_{i+}'=\Ran q_{i+}\,$,
$\,P_{-}\pi_i^\dag\pi_i'P_{+}'=0\,$, and $\,P_{-}'\pi_i'^{*}\pi_iP_{+}=0\,$, where
$P_{+}'$ is the orthoprojection onto $\,E_{\alpha}^2(G_+,{\mathfrak N_{i}})\,$ and
$\,P_{-}'=I-P_{+}'\,$. Define also the projections
$\,q_{i-}'=\pi_i'P_{-}'\pi_i'^{*}\,$. Then we have
\[
q_{i-}q_{j+}'=0\qquad\quad\hbox{and}\quad\qquad q_{i-}'q_{j+}=0\,,\qquad\quad\,i\ge
j\,.
\]
Indeed,
\[
\begin{array}{ll}
q_{i-}q_{j+}'&=\pi_iP_{-}\pi_i^{\dag}\pi_j'P_{+}'\pi_j'^{*}
=\pi_iP_{-}\pi_i^{\dag}\pi_j\pi_j^{\dag}\pi_j'P_{+}'\pi_j'^{*}\\[4pt]
&=\pi_iP_{-}\pi_i^{\dag}\pi_jP_{+}\pi_j^{\dag}\pi_j'P_{+}'\pi_j'^{*}+
\pi_iP_{-}\pi_i^{\dag}\pi_jP_{-}\pi_j^{\dag}\pi_j'P_{+}'\pi_j'^{*}=0+0=0\,.
\end{array}
\]
By the same reason, $\,q_{i-}'q_{j+}=0\,$. Using these identities and repeating
mutatis mutandis proof of Lemma~\ref{lemm_mod2}, we obtain
\[
\begin{array}{l}
P_{(ij)}'q_{k+}=P_{(ij)}q_{k+}'=0\,;\qquad\quad  q_{i-}'P_{(jk)}=q_{i-}P_{(jk)}'=0\,;\\[7pt]
P_{(ij)}'P_{(kl)}=P_{(ij)}P_{(kl)}'=P_{(ij)}'P_{(kl)}'=0,\qquad i\ge j\ge k\ge l\,.
\end{array}
\]
Then, evidently,
\[
P_{[m_i m_j]}P_{[m_k m_l]}'=P_{[m_i m_j]}'P_{[m_k m_l]}=P_{[m_i m_j]}'P_{[m_k
m_l]}'=0\,.
\]
Since $\,\Ran q_{i+}'=\Ran q_{i+}\,$, we have
$\,\mathcal{D}_{i+}'=\mathcal{D}_{i+}\,$. Evidently,
$\,\mathcal{H}_{ij}'=\mathcal{H}_{ij}\,$. Further, let $\,f\in
\mathcal{H}_{ij+}=\mathcal{H}_{ij}\cap\Ker q_{i-}\,$. Then $\,q_{i-}f=0\,$, that is
$\,\pi_i{\pi_i}^{\dag}f=q_{i+}f\,$, and
\[
q_{i-}'f=\pi_i'{\pi_i'}^*f-q_{i+}'f=\pi_i{\pi_i}^{\dag}f-q_{i+}'f=q_{i+}f-q_{i+}'f \in
\mathcal{D}_{i+}'=\mathcal{D}_{i+}\,.
\]
Therefore, $\,q_{i-}'f=q_{i-}'^2f=q_{i-}'q_{i+}'g=0\,$ and
$\,\mathcal{H}_{ij+}\subset\mathcal{H}_{ij+}'\,$. For  the same reason,
$\,\mathcal{H}_{ij+}'\subset\mathcal{H}_{ij+}\,$. Thus, we have
\[
\mathcal{D}_{i+}'=\mathcal{D}_{i+}\,,\qquad
\mathcal{H}_{ij}'=\mathcal{H}_{ij}\,,\qquad \mathcal{H}_{ij+}'=\mathcal{H}_{ij+}\,
\]
and therefore
\[
\mathcal{K}_{[m_i m_j]}'\subset\mathcal{H}_{m_i m_j+}\,,\qquad
\mathcal{D}_{m_j+}\subset \Ker P_{[m_i m_j]}'\,.
\]

\vskip 10pt\noindent   The following Proposition affirms a more delicate property of
projections $\,P_{[m_i m_j]}\,$.
\begin{prop}\label{dir_sum}
One has $\;\mathcal{H}_{m_i m_j+}\cap \Ker P_{[m_i m_j]}=\mathcal{D}_{m_j+}\,,\;i\ge j\,$.
\end{prop}
\noindent But beforehand we need to prove the following elementary lemmas, which are
of interest in their own right.
\begin{lemm}\label{lem_i} Suppose $\,M,\,N_+,\,N_-\,$ are subspaces of a Hilbert
space and $\,N_+\bot N_-\,$. Then $\;\;(N\vee M)\ominus N_-=((N\vee M)\ominus N)\oplus
N_+\,$, where $\,N=N_+\oplus N_-\,$.
\end{lemm}
\proof Let $\,f\in (N\vee M)\ominus N_-\,$. Then $\,f\in (N\vee M),\;f\bot N_-\,$. We
have $\,f=f_N+f_N^\bot\,$, where $\,f_N\in N,\; f_N^\bot\in N^\bot\,$. Then
$\,f_N=f-f_N^\bot\,\bot\, N_-\,$ and $\,f_N\in N_+\,$. Hence, $\,f=f_N^\bot+f_N\in
((N\vee M)\ominus N)\oplus N_+\,$.

Conversely, let $\,f\in N_+\,$. Then $\,f\in N,\;f\bot N_-\,$ and therefore $\,f\in
(N\vee M)\ominus N_-\,$. Hence, $\,((N\vee M)\ominus N)\oplus N_+\subset (N\vee
M)\ominus N_-\,$. \proofend
\begin{lemm}\label{lem_ii} Suppose $P\in\mathcal{L}(\mathcal{H})$ is an projection;
$\,\mathcal{D}_+,\,\mathcal{H}_+\,$ are subspaces of $\,\mathcal{H}\,$ such that
$\,\mathcal{D}_+\subset\mathcal{H}_+\,$, $\,\mathcal{K}=\Ran P\subset\mathcal{H}_+\,$
and $\,\mathcal{D}_+\subset\Ker P\,$. Then the following conditions are
equivalent:\\[4pt]
\,\,1)\,\,$\,\mathcal{H}_+\cap\Ker
P=\mathcal{D}_+\,$;\,\,\,\,\,\,\,\,2)\,\,$\,\mathcal{H}_+=\mathcal{K}\dot+\mathcal{D}_+\,$;
\,\,\,\,\,\,\,\,3)\,\,$\,\Ker(P|\mathcal{H}_+)=\mathcal{D}_+\,$.
\end{lemm}
\proof \,\,$1)\,\Longrightarrow\,2)\,$ Let $\,f\in \mathcal{H}_+\,$. Then
$\,f=f_1+f_2\,$, where $\,f_1=Pf\in \mathcal{K}\,$ and $\,f_2=(I-P)f\in \Ker P\,$.
Since $\,\mathcal{K}\subset \mathcal{H}_+\,$, we get $\,f_2=f-Pf\in \mathcal{H}_+\,$
and therefore $\,f_2\in \mathcal{H}_+\cap\Ker P=\mathcal{D}_+\,$.\\[4pt]
\,\,$2)\,\Longrightarrow\,3)\,$ It is clear that $\,\mathcal{D}_+\subset
\Ker(P|\mathcal{H}_+)\,$. Let $\,f\in\Ker(P|\mathcal{H}_+)\subset\mathcal{H}_+\,$.
Then $\,f=f_1+f_2\,$, where $\,f_1\in \mathcal{K}\,$ and $\,f_2\in\mathcal{D}_+\,$.
Then $\,0=Pf=P(f_1+f_2)=f_1\,$ and therefore  $\,f=f_2\in\mathcal{D}_+\,$.\\[4pt]
\,\,$3)\,\Longrightarrow\,1)\,$ It is clear that
$\,\mathcal{D}_+\subset\mathcal{H}_+\cap\Ker P\,$. Let $\,f\in\mathcal{H}_+\cap\Ker
P\,$. Then $\,f\in\Ker(P|\mathcal{H}_+)=\mathcal{D}_+\,$. \proofend
\begin{lemm}\label{lem_iii}
Suppose $\,P_1\,$ and $\,P_2\,$ are projections such that $\,\Ker P_1=\Ker P_2\,$.
Then $\,P_1P_2=P_1\,$ and $\,P_2P_1=P_2\,$.
\end{lemm}
\proof Since $\,\Ran (I-P_2)=\Ker P_2\,$, we get $\,P_1(I-P_2)=0\,$. Hence,
$\,P_1P_2=P_1\,$.\proofend
\begin{cor}
Suppose $P_1,P_2\in\mathcal{L}(\mathcal{H})$ are projections;
$\,\mathcal{D}_+,\,\mathcal{H}_+\,$ are subspaces of $\,\mathcal{H}\,$ such that
$\,\mathcal{D}_+\subset\mathcal{H}_+\,$, $\,\Ran P_1\subset\mathcal{H}_+\,$, $\,\Ran
P_2\subset\mathcal{H}_+\,$, $\,\mathcal{D}_+\subset\Ker P_1\,$,
$\,\mathcal{D}_+\subset\Ker P_2\,$ and
$\,\Ker(P_1|\mathcal{H}_+)=\Ker(P_2|\mathcal{H}_+)=\mathcal{D}_+\,$. Then
$\,P_1P_2P_1=P_1\,$ and $\,P_2P_1P_2=P_2\,$.
\end{cor}
\proof It is clear that $\,P_1|\mathcal{H}_+,\,P_2|\mathcal{H}_+\,$ are projections.
By Lemma~\ref{lem_iii}, we have
$\,(P_1|\mathcal{H}_+)(P_2|\mathcal{H}_+)=P_1|\mathcal{H}_+\,$. Then
$\,P_1P_2P_1f=P_1P_2(P_1f)=P_1(P_1f)=P_1f\,$. \proofend

\noindent \textit{Proof} (of Proposition~\ref{dir_sum}). First, we prove our assertion
in the orthogonal context. Consider orthogonal projections
\[
P_{(ij)}'=P_{\pi_i\vee\dots\vee\pi_j}(I-q_{j+}')(I-q_{i-}')\,,
\]
Since operators $q_{j+}',q_{i-}'$ are selfadjoint, we have
$\,(q_{j+}'q_{i-}')^*=q_{i-}'^*q_{j+}'^*=q_{i-}'q_{j+}'=0\,$ and hence
\[
P_{(ij)}'=P_{\pi_i\vee\dots\vee\pi_j}(I-q_{j+}'-q_{i-}')
=P_{\pi_i\vee\dots\vee\pi_j}-q_{j+}'-q_{i-}'\,.
\]
Define subspaces $\,N_{k\pm}:=\Ran q_{k\pm}'\,$, $\,N_{k}:=N_{k+}\oplus
N_{k-}=\Ran\pi_k'\,\;k=\overline{1,n}\,$. Then we have
$\,P_{\pi_i\vee\dots\vee\pi_j}=q_{i-}'+P_{(ij)}'+q_{j+}'\,$ and
\[
\mathcal{H}_{ij}=N_{i-}\oplus \mathcal{K}_{(ij)}'\oplus N_{j+}\,,\qquad
\mathcal{H}_{ij+}=\mathcal{K}_{(ij)}'\oplus N_{j+}\,,\qquad \mathcal{D}_{j+}=
N_{j+}\,.
\]
In particular, we get $\,N_{k}\vee N_{k+1}=N_{k+1-}\oplus \mathcal{K}_{(k+1,k)}'\oplus
N_{k+}\,$ and therefore $\,N_{k+}\oplus\mathcal{K}_{(k+1,k)}'=(N_{k}\vee
N_{k+1})\ominus N_{k+1-}\,$. Applying the former identity and Lemma~\ref{lem_i} $i-j$
times, we have
\[
\begin{array}{c}
N_{j+}\oplus\mathcal{K}_{(j+1,j)}'\oplus\mathcal{K}_{(j+2,j+1)}'\oplus\ldots\oplus\mathcal{K}_{(i,i-1)}'\oplus
N_{i-}\\[8pt]
=[(N_{j}\vee N_{j+1})\ominus
N_{j+1-}]\oplus\mathcal{K}_{(j+2,j+1)}'\oplus\ldots\oplus\mathcal{K}_{(i,i-1)}'\oplus N_{i-}\\[8pt]
=[(N_{j}\vee N_{j+1})\ominus
N_{j+1}]\oplus N_{j+1+}\oplus\mathcal{K}_{(j+2,j+1)}'\oplus\ldots\oplus\mathcal{K}_{(i,i-1)}'\oplus N_{i-}\\[8pt]
=[(N_{j}\vee N_{j+1})\ominus
N_{j+1}]\oplus [N_{j+1+}\oplus\mathcal{K}_{(j+2,j+1)}']\oplus\ldots\oplus\mathcal{K}_{(i,i-1)}'\oplus N_{i-}\\[8pt]
=\ldots=[(N_{j}\vee N_{j+1})\ominus N_{j+1}]\oplus [(N_{j+1}\vee N_{j+2})\ominus
N_{j+2}]\oplus\ldots\\[8pt]
\oplus[(N_{i-1}\vee N_{i})\ominus N_{i}]\oplus N_{i+}\oplus N_{i-}=[(N_{j}\vee
N_{j+1})\ominus N_{j+1}]\oplus\\[8pt]
[(N_{j+1}\vee N_{j+2})\ominus N_{j+2}]\oplus\ldots\oplus[(N_{i-1}\vee N_{i})\ominus
N_{i}]\oplus N_{i}\\[8pt]
=[(N_{j}\vee N_{j+1})\ominus N_{j+1}]\oplus[(N_{j+1}\vee N_{j+2})\ominus
N_{j+2}]\oplus\ldots\oplus[(N_{i-1}\vee N_{i})]\\[8pt]
=\ldots=N_{j}\vee N_{j+1}\vee\ldots\vee N_{i-1}\vee N_{i}=\mathcal{H}_{ij}\,.
\end{array}
\]
On the other hand, we have already shown $\,\mathcal{H}_{ij}=N_{i-}\oplus
\mathcal{K}_{(ij)}'\oplus N_{j+}\,$. Therefore,
\[
\mathcal{K}_{(ij)}'=\mathcal{K}_{(ii-1)}'\oplus\ldots\oplus\mathcal{K}_{(j+1j)}'\,.
\]
Then we have $\,\mathcal{K}_{[m_im_j]}'=\mathcal{K}_{(m_im_j)}'\,$ and
$\,\mathcal{K}_{[m_i m_k]}'=\mathcal{K}_{[m_im_j]}'\oplus \mathcal{K}_{[m_jm_k]}'\,$.
It is easy to check that $\;\mathcal{H}_{m_i m_j+}\cap \Ker P_{[m_i
m_j]}'=\mathcal{D}_{m_j+}\,$.

For the nonorthogonal case, we shall use induction. In fact, we have already shown
that
\[
\mathcal{H}_{m_{j+1} m_j+}\cap \Ker P_{[m_{j+1} m_j]}=\mathcal{H}_{m_{j+1}
m_j+}\cap\Ker P_{(m_{j+1} m_j)}=\mathcal{D}_{m_j+}\,,\;j=\overline{1,n}\,.
\]
Let $\;i\ge j\ge k\,$. Assume that $\;\mathcal{H}_{m_i m_j+}\cap \Ker P_{[m_i
m_j]}=\mathcal{D}_{m_j+}\,$ and $\;\mathcal{H}_{m_j m_k+}\cap \Ker P_{[m_j
m_k]}=\mathcal{D}_{m_k+}\,$. Let $\,f\in \mathcal{H}_{m_i m_k+}\cap \Ker P_{[m_i
m_k]}\,$. Using the recursion relation
\[
P_{[m_im_k]}=P_{[m_jm_k]}(I-P_{[m_im_j]})+P_{[m_im_j]}\,
\]
and properties of projection $\,P_{[\cdot\,\cdot]}\,$, we have
\[
P_{[m_im_j]}f=P_{[m_im_j]}(P_{[m_jm_k]}(I-P_{[m_im_j]})+P_{[m_im_j]})f=P_{[m_im_j]}P_{[m_im_k]}f=0\,.
\]
Then, since $\,P_{[m_im_k]}f=0\,$, we also have $\,P_{[m_jm_k]}f=0\,$. On the other
hand, the vector $\,f\,$ can be decomposed $\,f=f_{ij}'+f_{jk}'+g\,$, where
$\,f_{ij}'\in\mathcal{K}_{[m_im_j]}',\,f_{ij}'\in\mathcal{K}_{[m_im_j]}'\,$ and
$\,g\in\mathcal{D}_{m_j+}\,$. Since $\,P_{[m_im_j]}P_{[m_jm_k]}'=0\,$, we have
\[
0=P_{[m_im_j]}f=P_{[m_im_j]}(f_{ij}'+f_{jk}'+g)=P_{[m_im_j]}f_{ij}'\,.
\]
By Lemma~\ref{lem_ii}, $\,\Ker(P_{[m_im_j]}|\mathcal{H}_{m_i
m_j+})=\mathcal{D}_{m_j+}\,$. Then, by Corollary of Lemma~\ref{lem_iii}, we obtain
$\,0=P_{[m_im_j]}'P_{[m_im_j]}f_{ij}'=P_{[m_im_j]}'P_{[m_im_j]}P_{[m_im_j]}'f_{ij}'=f_{ij}'\,$.
Further, $\,0=P_{[m_jm_k]}f=P_{[m_jm_k]}(f_{jk}'+g)=P_{[m_jm_k]}f_{jk}'\,$. As above,
we get $\,0=f_{jk}'\,$. Thus, we have $\,f=g\in\mathcal{D}_{m_j+}\,$ and therefore
$\,\mathcal{H}_{m_i m_k+}\cap \Ker P_{[m_i m_k]}\subset\mathcal{D}_{m_j+}\,$. The
inverse inclusion is obvious. \proofend

\myrem\,\, Since $\,P_{[31]}=P_{(32)}+P_{(21)}(I-P_{(32)})=P_{(32)}+P_{(21)}\,$, by
Prop.~\ref{dir_sum} and Corollary of Lemma~\ref{lem_iii}, we obtain the following
identities
\[
(P_{(32)}+P_{(21)})\,P_{(31)}\,(P_{(32)}+P_{(21)})\,=\,P_{(32)}+P_{(21)}
\]
and
\[
P_{(31)}\,(P_{(32)}+P_{(21)})\,P_{(31)}\,=\,P_{(31)}\,.
\]
This means that
\[
\;(P_{(31)}|\mathcal{K}_{[31]})^{-1}=(P_{(21)}+P_{(32)})|\mathcal{K}_{(31)}\,,\quad
((P_{(21)}+P_{(32)})|\mathcal{K}_{(31)})^{-1}=P_{(31)}|\mathcal{K}_{[31]}\,.
\]

\vskip 7pt \exmpl\,\, Let $\,w=\varphi(z)=z+\varepsilon
z^2\,,\;|\varepsilon|<1/2\,,\;G_+=\varphi(\mathbb{D})\,,\;C=\varphi(\mathbb{T})\,$. We put
\[
\theta(w)=\frac{2w}{1+\sqrt{1+4\varepsilon w}}\,,\quad\,w\in G_+\,,\qquad
\Theta_{ij}(w)=\theta(w)^{i-j}\,,\quad 1\le j\le i\le n\,,
\]
and $\,\Xi_i(w)=1\,,w\in C\;\,$. It can easily be checked that $\,|\theta(w)|=1\,$,
$\,w\in C\,$. Then
$\,P_{(ij)}=P_{(ij)}^{(n)}=\theta^{n-j}P_-\theta^{j-i}P_+\theta^{i-n}\,$. For the
functions
\[
f_k^{ij}(w)=\theta(w)^{n-j}\,w^{-k}\,,\qquad k=\overline{1,i-j}\,,\qquad
\]
we have $\,f_k^{ij}\in\mathcal{K}_{(ij)}^{(n)}=\Ran P_{(ij)}^{(n)}\,$. By~\cite{T2},
$\,\mathcal{K}_{(ij)}^{(n)}(\varepsilon)=P_{(ij)}^{(n)}(\varepsilon)\mathcal{K}_{(ij)}^{(n)}(0)\,$
and
$\,\mathcal{K}_{(ij)}^{(n)}(0)=P_{(ij)}^{(n)}(0)\mathcal{K}_{(ij)}^{(n)}(\varepsilon)\,$.
Since $\,\dim \mathcal{K}_{(ij)}^{(n)}(0)=i-j\,$, the functions $\,f_k^{ij}(w)\,$ form
a basis of the subspace $\,\mathcal{K}_{(ij)}^{(n)}\,$. Note also that
$\,P_{ij}^{(n)}=\theta^{n-i}P_{i-j+1,1}^{(i-j+1)}\theta^{i-n}\,$ and therefore
$\,\mathcal{K}_{(ij)}^{(n)}=\theta^{n-i}\mathcal{K}_{(i-j+1,1)}^{(i-j+1)}\,$.\\

\noindent Consider particular cases. In the case of $\,n=3\,$ we have
\[
f_1^{31}=\theta^2/w\,,\,f_2^{31}=\theta^2/w^2\qquad\mbox{and}\qquad
f_1^{21}=\theta^2/w\,,\;\,f_1^{32}=\theta/w\,.
\]
Hence, $\,\mathcal{K}_{(21)}\subset\mathcal{K}_{(31)}\,$,
$\,\mathcal{K}_{(32)}\nsubseteq\mathcal{K}_{(31)}\,$ and
$\;\mathcal{K}_{(32)}\dot+\mathcal{K}_{(21)}\ne\mathcal{K}_{(31)}\,$.\\

\noindent In the case of $\,n=5\,$ it can be calculated that
\[
P_{(21)}f_1^{53}=-\varepsilon^2 f_1^{21}\qquad\mbox{and}\qquad P_{(21)}f_2^{53}=2\varepsilon^3
f_1^{21}\,.
\]
Therefore, $\,P_{(21)}P_{(53)}\ne 0\,$. Our calculations are based on the formula
\[
P_{(kl)}f_p^{ij}=\theta^{n-l}w^{l-j-p}u_{l-j-p,j-l}(w)
-2^{l-k}\theta^{n-l}P_-\frac{w^{l-j-p}u_{k-j-p,j-k}(w)}{(1+\sqrt{1+4\varepsilon
w})^{l-k}}\,,
\]
where $\,1\le i,j,k,l\le n\,$, $\,p=\overline{1,i-j}\,$, $\,i\ge j\,,\;k\ge l\,$, and
\[
u_{q,r}(w):=2^{-r}w^{-q}P_-w^q (1+\sqrt{1+4\varepsilon w})^{r}\,.
\]
It can easily be checked that $\,u_{q,r}(w)\equiv 0\,,\;q\ge 0\,$. For $\,q<0\,$, we
make use of the Residue Theorem calculating $\,P_+w^q (1+\sqrt{1+4\varepsilon
w})^{r}\,$ and interpreting the projection $\,P_+\,$ as the boundary values of the
Cauchy integral operator. In particular, we get
\[
\begin{array}{ll}
u_{-1,r}(w)&=\;1\,;\qquad u_{-2,r}(w)\;=\;1+r\varepsilon w\,;\\[8pt]
u_{-3,r}(w)&=\;\frac{1}{2}(2+2r\varepsilon w+r(r-3)\varepsilon^2 w^2)\,;\\[8pt]
u_{-4,r}(w)&=\;\frac{1}{6}(2+6r\varepsilon w+3r(r-3)\varepsilon^2
w^2+3r(r-4)(r-5)\varepsilon^3 w^3)\,.
\end{array}
\]

\section{Product of conservative curved systems}

\begin{defin}
Let
$\,\Sigma_k=(T_{k},M_{k},N_{k},\Theta_{ku},\Xi_{k};H_{k},\mathfrak{N}_{k+},\mathfrak{N}_{k-})\,,\;k={1,2}\,$
be conservative curved systems, $\,G_{1+}=G_{2+}\,$, $\,\mathfrak{N}_{1-}=\mathfrak{N}_{2+}\,$ and
$\;\,\Xi_{1-}=\Xi_{2+}$. We define the product of them as\\
\[
\Sigma_{21}=\Sigma_2\cdot\Sigma_1:=(T_{21},M_{21},N_{21},\Theta_{21u},\Xi_{21};
H_{21},\mathfrak{N}_{1+},\mathfrak{N}_{2-})
\]\\[2pt]
with $\;\;\Theta_{21u}=\Theta_{2u}\Theta_{1u}\,$,\, $\;\Xi_{21}=(\Xi_{1+},\Xi_{2-})$,\,
$\,H_{21}=H_1\oplus H_2$\,,\\[4pt]
\[
\begin{array}{l}
T_{21}=\left(
\begin{array}{cc}
  T_1 & N_1M_2 \\
  0 & T_2 \\
\end{array}
\right),\quad M_{21}=(M_1,M_2^{21})\,,\quad N_{21}= \left(
\begin{array}{c}
  M_{*1}^{21*} \\
  N_2 \\
\end{array}
\right)\,,\\[14pt]
M_2^{21}f_2=\displaystyle-\frac{1}{2\pi
i}\int\limits_C\Theta_1^-(\zeta)\,[M_2(T_2-\cdot)^{-1}f_2]_-(\zeta)\,d\zeta,\quad f_2\in H_2\,,\\[10pt]
M_{*1}^{21}f_1=\displaystyle-\frac{1}{2\pi
i}\int\limits_{\overline{C}}\Theta_{*\,2}^-(\zeta)\,[N_1^*(T_1^*-\cdot)^{-1}f_1]_-(\zeta)\,d\zeta,\quad f_1\in H_1\,,\\[10pt]
\end{array}\eqno{\rm{(Prod)}}
\]
where $\,[M_2(T_2-\cdot)^{-1}f_2]_{-}\,$ and $\,[N_1^*(T_1^*-\cdot)^{-1}f_1]_-\,$ are the boundary
limits of $\;M_2(T_2-z)^{-1}f_2\,$ and $\,N_1^*(T_1^*-z)^{-1}f_1\,$ from the domains $\,G_-\,$ and
$\,\bar{G}_-:=\{\bar{z} \colon z\in G_-\}\,$, respectively\,;
$\;\Theta_{*\,2}^-=\Theta_{2}^{-\sim}\,$ (see {\rm(CtoT)} for the definition of $\,\Theta^-\,$).
\end{defin}
\noindent Note that we can consider the product $\,\Sigma_2\cdot\Sigma_1\,$ without the assumption
that $\,\Sigma_1,\Sigma_2\,$ are conservative curved systems. We need only to assume additionally
that $\,\forall\,f_2\in H_2\; \colon M_{2}(T_{2}-z)^{-1}f_2\in E^2(G_-)\,$ and $\,\forall\,f_1\in
H_1\; \colon N_{1}^*(T_{1}^*-z)^{-1}f_1\in E^2(\bar{G}_-)\;$. For conservative curved systems these
assumptions are always satisfied (it follows from the definition of conservative curved system).

\vskip 2pt We start to justify the definition with the observation that in case of
unitary colligations we get the standard algebraic definition~\cite{Br}:
$\,M_2^{21}=\Theta_1^+(0)^*M_2\,$ and $\,N_1^{21}:=M_{*1}^{21*}=N_1\Theta_2^+(0)^*\,$
(see the Introduction). Indeed, since in this case
$\,\Theta_{+}^-\equiv\Theta^+(0)^*=L\,$ and $\,M_{2}(T_{2}-z)^{-1}f_2\in E^2(G_-)\,$,
we obtain
\[
\begin{array}{ll}
M_2^{21}f_2&=\displaystyle-\frac{1}{2\pi
i}\int\limits_C\Theta_1^-(\zeta)\,[M_2(T_2-\cdot)^{-1}f_2]_-(\zeta)\,d\zeta\\[8pt]
&=\displaystyle-\frac{1}{2\pi
i}\int\limits_C\Theta_{1+}^-(\zeta)\,[M_2(T_2-\cdot)^{-1}f_2]_-(\zeta)\,d\zeta=L_1M_2\,.\\[8pt]
\end{array}
\]
By a similar computation, we get $\,N_1^{21}=N_1L_2\,$. Besides, we have
\begin{prop}\label{prod_sim}
\,1)
$\,\Sigma_1\sim\Sigma_1',\,\Sigma_2\sim\Sigma_2'\Rightarrow\Sigma_2\cdot\Sigma_1\sim\Sigma_2'\cdot\Sigma_1'$;
\,2) $\,(\Sigma_2\cdot\Sigma_1)^*=\Sigma_1^*\cdot\Sigma_2^*\,$.
\end{prop}
\noindent Here $\,\Sigma^*:=(T^*, N^*, M^*, \Theta_u^{\sim},
\Xi_*)\,,\;\Xi_{*\pm}=\Xi_{\mp}^{\sim-1}\,$. \proof\,\,1)\,\,\,\, Let
$\,\Sigma_1\,{\stackrel{X_1}{\sim}}\,\Sigma_1',\;\Sigma_2\,{\stackrel{X_2}{\sim}}\,\Sigma_2'\,$.
Then
\[
\begin{array}{ll}
M_2^{21}f_2&=\displaystyle-\frac{1}{2\pi
i}\int\limits_C\Theta_1^-(\zeta)\,[M_2(T_2-\cdot)^{-1}f_2]_-(\zeta)\,d\zeta\\[8pt]
&=\displaystyle-\frac{1}{2\pi
i}\int\limits_C\Theta_{1}^-(\zeta)\,[M_2'(T_2'-\cdot)^{-1}X_2f_2]_-(\zeta)\,d\zeta={M_2^{21}}'X_2f_2\,.\\[8pt]
\end{array}
\]
Hence, $\,M_2^{21}={M_2^{21}}'X_2\,$. Similarly, $\,X_1N_1^{21}={N_1^{21}}'\,$. Define
$\;X_{21}:=\left(
\begin{array}{cc}
  X_1 & 0 \\
  0 & X_2 \\
\end{array}
\right)\,$. Then, we get
\[
\begin{array}{rl}
T_{21}'X_{21}&= \left(
\begin{array}{cc}
  T_1' & N_1'M_2' \\
  0 & T_2' \\
\end{array}
\right) \left(
\begin{array}{cc}
  X_1 & 0 \\
  0 & X_2 \\
\end{array}
\right)= \left(
\begin{array}{cc}
  T_1'X_1 & N_1'M_2'X_2 \\
  0 & T_2'X_2 \\
\end{array}
\right)\\[14pt]
&=\left(
\begin{array}{cc}
  X_1T_1 & X_1N_1M_2 \\
  0 & X_2T_2 \\
\end{array}
\right) = \left(
\begin{array}{cc}
  X_1 & 0 \\
  0 & X_2 \\
\end{array}
\right) \left(
\begin{array}{cc}
  T_1 & N_1M_2 \\
  0 & T_2 \\
\end{array}
\right)=X_{21}T_{21}\,;\\[24pt]
M_{21}'X_{21}&=\left(
\begin{array}{cc}
  M_1'\,, & {M_2^{21}}'
\end{array}
\right) \left(
\begin{array}{cc}
  X_1 & 0 \\
  0 & X_2 \\
\end{array}
\right)=\left(
\begin{array}{cc}
  M_1'X_1\,, & {M_2^{21}}'X_2
\end{array}
\right)\\[14pt]
&=\left(
\begin{array}{cc}
  M_1\,, & {M_2^{21}}
\end{array}
\right)=M_{21}\,;\\[24pt]
N_{21}'&= \left(
\begin{array}{cc}
  {N_1^{21}}' \\ N_2'
\end{array}
\right)= \left(
\begin{array}{cc}
  X_1{N_1^{21}}' \\ X_2N_2'
\end{array}
\right) =\left(
\begin{array}{cc}
  X_1 & 0 \\
  0 & X_2 \\
\end{array}
\right) \left(
\begin{array}{cc}
  {N_1^{21}} \\ N_2
\end{array}
\right)=X_{21}N_{21}\,\\[4pt]
\end{array}
\]
and thus
$\,\Sigma_2\cdot\Sigma_1\;{\stackrel{X_{21}}{\sim}}\;\Sigma_2'\cdot\Sigma_1'\,$.\\

\noindent  2)\,\,\,\, Let
$\,\Sigma_{21*}=(\Sigma_2\cdot\Sigma_1)^*\,$. By straightforward calculation, we get\\
\[
T_{*21}=\left(
\begin{array}{cc}
  T_1^* & 0 \\
  M_2^*N_1^* & T_2^* \\
\end{array}
\right);\quad
M_{*21}=\left(
\begin{array}{cc}
  N_1^{21*}, & N_2^* \\
\end{array}
\right);\quad
N_{*21}=\left(
\begin{array}{cc}
  M_1^*  \\
  M_2^{21*} \\
\end{array}
\right).
\]
\vskip 4pt \noindent On the other hand, let
$\,\Sigma_{21*}'=\Sigma_1^*\cdot\Sigma_2^*\,$. Then we get\\
\[
T_{*21}'=\left(
\begin{array}{cc}
  T_{*1} & 0 \\
  N_{*2}M_{*1} & T_{*2} \\
\end{array}
\right)=\left(
\begin{array}{cc}
  T_1^* & 0 \\
  M_2^*N_1^* & T_2^* \\
\end{array}
\right)
=T_{*21}.
\]
\vskip 4pt \noindent Since
\[
\begin{array}{ll}
{M_{*1}^{21}}'f_1&=\displaystyle-\frac{1}{2\pi
i}\int\limits_{\overline{C}}\Theta_{*\,2}^-(\zeta)\,[M_{*1}(T_{*1}-\cdot)^{-1}f_1]_-(\zeta)\,d\zeta\\[10pt]
&=\displaystyle-\frac{1}{2\pi
i}\int\limits_{\overline{C}}\Theta_{*\,2}^-(\zeta)\,[N_1^*(T_1^*-\cdot)^{-1}f_1]_-(\zeta)\,d\zeta
=N_{1}^{21*}f_1,\quad f_1\in H_1\,,\\[10pt]
\end{array}
\]
we have
\[
M_{*21}'=\left(
\begin{array}{cc}
  {M_{*1}^{21}}'\,, & M_{*2} \\
\end{array}
\right)=\left(
\begin{array}{cc}
  {N_{1}^{21*}}'\,, & N_{2}^* \\
\end{array}
\right)
=M_{*21}\,.
\]
\vskip 4pt \noindent Since
\[
\begin{array}{ll}
{N_{*2}^{21}}'^*f_2&=\displaystyle-\frac{1}{2\pi
i}\int\limits_C\Theta_{**1}^-(\zeta)\,[N_{*2}^*(T_{*2}^*-\cdot)^{-1}f_2]_-(\zeta)\,d\zeta\\[10pt]
&=\displaystyle-\frac{1}{2\pi
i}\int\limits_C\Theta_1^-(\zeta)\,[M_2(T_2-\cdot)^{-1}f_2]_-(\zeta)\,d\zeta
=M_{2}^{21}f_2,\quad f_2\in H_2\,,\\[10pt]
\end{array}
\]
we have
\[
N_{*21}'=\left(
\begin{array}{cc}
  N_{*1}\\ {N_{*2}^{21}}' \\
\end{array}
\right)=\left(
\begin{array}{cc}
  M_{1}^*\\ {M_{1}^{21*}} \\
\end{array}
\right) =N_{*21}\,.
\]
\vskip 4pt \noindent Thus we get
$\,\Sigma_1^*\cdot\Sigma_2^*=\Sigma_{21*}'=\Sigma_{21*}=(\Sigma_2\cdot\Sigma_1)^*\,$.\proofend

\noindent  Further, we shall say that a triple of operators $\,(T,M,N)\,$ is a
realization of a transfer function $\,\Upsilon=\mathcal{F}_{tc}(\Theta)\,$ if
$\,\Upsilon(z)=M(T-z)^{-1}N\,$.
\begin{prop}
Suppose that triples $\,(T_{1},M_{1},N_{1})\,$ and $\,(T_{2},M_{2},N_{2})\,$ are
realizations of transfer functions $\;\Upsilon_{1}=\mathcal{F}_{tc}(\Theta_{1})\,$ and
$\;\Upsilon_{2}=\mathcal{F}_{tc}(\Theta_{2})\,$, respectively. Suppose also that
$\,\forall\,f_1\in H_1\; \colon N_{1}^*(T_{1}^*-z)^{-1}f_1\in E^2(\bar{G}_-)\;$ and
$\;\forall\,f_2\in H_2\; \colon M_{2}(T_{2}-z)^{-1}f_2\in E^2(G_-)\,$. Then the triple
$\,(T_{21},M_{21},N_{21})\,$ defined by {\rm (Prod)} is a realization of the transfer
function $\;\Upsilon_{21}=\mathcal{F}_{tc}(\Theta_{2}\Theta_{1})\,$.
\end{prop}
\proof  For the sake of simplicity, consider the case $\,\Theta_{k-}^-\in
H^{\infty}(G_-,{\mathcal L}({\mathfrak N}_{k-},{\mathfrak N}_{k+})))\,,$ $\;k=1,2\,$
(in the general case we need to use expressions like
$\,(M_{21}(T_{21}-\lambda)^{-1}N_{21}n,m)\,$). It can easily be shown that\\
\[
(T_{21}-\lambda)^{-1}=\left(
\begin{array}{cc}
  (T_1-\lambda)^{-1} & -(T_1-\lambda)^{-1}N_1M_2(T_2-\lambda)^{-1} \\
  0 & (T_2-\lambda)^{-1} \\
\end{array}
\right).
\]
\vskip 4pt \noindent Then, by straightforward computation, we obtain\\
\[
M_{21}(T_{21}-\lambda)^{-1}N_{21}=M_{1}(T_{1}-\lambda)^{-1}N_{1}^{21}
-\Upsilon_1(\lambda)\Upsilon_2(\lambda)+M_{2}^{21}(T_{2}-\lambda)^{-1}N_{2}\,.
\]
\vskip 4pt \noindent Here, we have\\
\[
\begin{array}{l}
M_{2}^{21}(T_{2}-\lambda)^{-1}N_{2}=\\[14pt]
\qquad\quad=\displaystyle-\frac{1}{2\pi
i}\int\limits_C\Theta_1^-(\zeta)\,[M_2(T_2-\cdot)^{-1}(T_{2}-\lambda)^{-1}N_{2}]_-(\zeta)\,d\zeta\\[14pt]
\qquad\quad=\displaystyle-\frac{1}{2\pi
i}\int\limits_C\frac{\Theta_1^-(\zeta)}{\zeta-\lambda}\,
[M_2(T_2-\cdot)^{-1}N_2-M_2(T_{2}-\lambda)^{-1}N_{2}]_-(\zeta)\,d\zeta\\[14pt]
\qquad\quad=\displaystyle-\frac{1}{2\pi
i}\int\limits_C\frac{\Theta_1^-(\zeta)(\Upsilon_2(\zeta)_--\Upsilon_2(\lambda))}{\zeta-\lambda}\,
\,d\zeta\\[14pt]
\qquad\quad=\displaystyle-\frac{1}{2\pi
i}\int\limits_C\frac{\Theta_{1+}^-(\zeta)(\Upsilon_2(\zeta)_--\Upsilon_2(\lambda))}{\zeta-\lambda}\,
\,d\zeta\\[14pt]
\qquad\quad=\displaystyle\frac{1}{2\pi
i}\int\limits_C\frac{\Theta_{1+}^-(\zeta)\Theta_{2-}^-(\zeta)}{\zeta-\lambda}\,
\,d\zeta\,+\,\frac{1}{2\pi
i}\int\limits_C\frac{\Theta_{1+}^-(\zeta)\Upsilon_2(\lambda)}{\zeta-\lambda}\,
\,d\zeta\\[14pt]
\qquad\quad=\displaystyle\frac{1}{2\pi
i}\int\limits_C\frac{\Theta_{1+}^-(\zeta)\Theta_{2-}^-(\zeta)}{\zeta-\lambda}\,
\,d\zeta\,+\left\{%
\begin{array}{ll}
    \Theta_{1+}^-(\lambda)(\Theta_{2+}^-(\lambda)-\Theta_{2}^+(\lambda)^{-1}), & \lambda\in G_+ \\
    0, & \lambda\in G_- \\
\end{array}%
\right.\,.    \\[14pt]
\end{array}
\]
\vskip 4pt \noindent Similarly, we have\\
\[
\begin{array}{l}
N_{1}^{21*}(T_{1}^*-\overline{\lambda})^{-1}M_{1}^*=\\[14pt]
\qquad\quad=\displaystyle\frac{1}{2\pi
i}\int\limits_{\overline{C}}\frac{\Theta_{*2+}^-(\zeta)\Theta_{*1-}^-(\zeta)}{\zeta-\overline{\lambda}}\,
\,d\zeta\,+\left\{%
\begin{array}{ll}
    \Theta_{*2+}^-(\overline{\lambda})(\Theta_{*1+}^-(\overline{\lambda})-\Theta_{*1}^+(\overline{\lambda})^{-1}), & \lambda\in G_+ \\
    0, & \lambda\in G_- \\
\end{array}%
\right.    \\[14pt]
\qquad\quad=\displaystyle-\frac{1}{2\pi
i}\int\limits_{{C}}\frac{\Theta_{*2+}^-(\overline{\zeta})\Theta_{*1-}^-(\overline{\zeta})}
{\overline{\zeta}-\overline{\lambda}}\,\,d\overline{\zeta}\,+\left\{%
\begin{array}{ll}
    \Theta_{*2+}^-(\overline{\lambda})(\Theta_{*1+}^-(\overline{\lambda})-\Theta_{*1}^+(\overline{\lambda})^{-1}), & \lambda\in G_+ \\
    0, & \lambda\in G_- \\
\end{array}%
\right.    \\[14pt]
\end{array}
\]
\vskip 4pt \noindent Hence,\\
\[
\begin{array}{l}
M_{1}(T_{1}-\lambda)^{-1}N_{1}^{21}=\\[14pt]
\qquad\quad=\displaystyle\frac{1}{2\pi
i}\int\limits_{{C}}\frac{\Theta_{*1-}^-(\overline{\zeta})^*\Theta_{*2+}^-(\overline{\zeta})^*}{\zeta-{\lambda}}\,
\,d\zeta\,+\left\{%
\begin{array}{ll}
    (\Theta_{*1+}^-(\overline{\lambda})^*-(\Theta_{*1}^+(\overline{\lambda})^*)^{-1})\,\Theta_{*2+}^-(\overline{\lambda})^*, & \lambda\in G_+ \\
    0, & \lambda\in G_- \\
\end{array}%
\right.    \\[14pt]
\qquad\quad=\displaystyle\frac{1}{2\pi
i}\int\limits_{{C}}\frac{\Theta_{1-}^-(\zeta)\Theta_{2+}^-(\zeta)}{\zeta-{\lambda}}\,
\,d\zeta\,+\left\{%
\begin{array}{ll}
    (\Theta_{1+}^-({\lambda})-\Theta_{1}^+({\lambda})^{-1})\,\Theta_{2+}^-({\lambda}), & \lambda\in G_+ \\
    0, & \lambda\in G_- \\
\end{array}%
\right.\,.    \\[14pt]
\end{array}
\]
\vskip 4pt \noindent Consider the case when $\,\lambda\in G_-\,$. Then\\
\[
\begin{array}{l}
M_{21}(T_{21}-\lambda)^{-1}N_{21}=\\[14pt]
\qquad=\displaystyle\frac{1}{2\pi
i}\int\limits_{{C}}\frac{\Theta_{1-}^-(\zeta)\Theta_{2+}^-(\zeta)}{\zeta-{\lambda}}\,
\,d\zeta\,-\,\Theta_{1-}^-(\lambda)\Theta_{2-}^-(\lambda)\,+\,\frac{1}{2\pi
i}\int\limits_{{C}}\frac{\Theta_{1+}^-(\zeta)\Theta_{2-}^-(\zeta)}{\zeta-{\lambda}}\,
\,d\zeta+0\,\\[14pt]
\qquad=\displaystyle\frac{1}{2\pi
i}\int\limits_{{C}}\frac{\Theta_{1-}^-(\zeta)\Theta_{2+}^-(\zeta)}{\zeta-{\lambda}}\,
\,d\zeta\,+\,\frac{1}{2\pi
i}\int\limits_{{C}}\frac{\Theta_{1-}^-(\zeta)\Theta_{2-}^-(\zeta)}{\zeta-{\lambda}}\,
\,d\zeta\,\\[14pt]
\qquad+\displaystyle\frac{1}{2\pi
i}\int\limits_{{C}}\frac{\Theta_{1+}^-(\zeta)\Theta_{2-}^-(\zeta)}{\zeta-{\lambda}}\,
\,d\zeta\,+\,\frac{1}{2\pi
i}\int\limits_{{C}}\frac{\Theta_{1+}^-(\zeta)\Theta_{2+}^-(\zeta)}{\zeta-{\lambda}}\,
\,d\zeta\,\\[14pt]
\qquad=\displaystyle\frac{1}{2\pi
i}\int\limits_{{C}}\frac{(\Theta_{1+}^-(\zeta)+\Theta_{1-}^-(\zeta))
(\Theta_{2+}^-(\zeta)+\Theta_{2-}^-(\zeta))}{\zeta-{\lambda}}\,
\,d\zeta\\[14pt]
\qquad=\displaystyle\frac{1}{2\pi
i}\int\limits_{{C}}\frac{\Theta_{1}^-(\zeta)\Theta_{2}^-(\zeta)}{\zeta-{\lambda}}\,
\,d\zeta =\frac{1}{2\pi
i}\int\limits_{{C}}\frac{\Theta_{21}^-(\zeta)}{\zeta-{\lambda}}\,
\,d\zeta=-\Theta_{21}^-(\lambda)\;.\\[14pt]
\end{array}
\]
\vskip 4pt \noindent Consider the case when $\,\lambda\in G_+\,$. Then\\
\[
\begin{array}{l}
M_{21}(T_{21}-\lambda)^{-1}N_{21}=\\[14pt]
\qquad=\displaystyle\frac{1}{2\pi
i}\int\limits_{{C}}\frac{\Theta_{1-}^-(\zeta)\Theta_{2+}^-(\zeta)}{\zeta-{\lambda}}\,
\,d\zeta\,+\,(\Theta_{1+}^-({\lambda})-\Theta_{1}^+({\lambda})^{-1})\,\Theta_{2+}^-({\lambda})\\[16pt]
\qquad-\;(\Theta_{1+}^-({\lambda})-\Theta_{1}^+({\lambda})^{-1})
(\Theta_{2+}^-(\lambda)-\Theta_{2}^+(\lambda)^{-1})\\[10pt]
\qquad\,+\,\displaystyle\frac{1}{2\pi
i}\int\limits_{{C}}\frac{\Theta_{1+}^-(\zeta)\Theta_{2-}^-(\zeta)}{\zeta-{\lambda}}\,
\,d\zeta+\Theta_{1+}^-(\lambda)(\Theta_{2+}^-(\lambda)-\Theta_{2}^+(\lambda)^{-1})\,\\[14pt]
\qquad=\displaystyle\frac{1}{2\pi
i}\int\limits_{{C}}\frac{\Theta_{1-}^-(\zeta)\Theta_{2+}^-(\zeta)+\Theta_{1+}^-(\zeta)\Theta_{2-}^-(\zeta)}
{\zeta-{\lambda}}\,\,d\zeta+\Theta_{1+}^-(\lambda)\Theta_{2+}^-(\lambda)-
\Theta_{1}^+(\lambda)^{-1}\Theta_{2}^+(\lambda)^{-1}\,\\[14pt]
\qquad=\displaystyle\frac{1}{2\pi
i}\int\limits_{{C}}\frac{(\Theta_{1}^-(\zeta)-\Theta_{1+}^-(\zeta))\,\Theta_{2+}^-(\zeta)+
(\Theta_{1}^-(\zeta)-\Theta_{1-}^-(\zeta))\,\Theta_{2-}^-(\zeta)}
{\zeta-{\lambda}}\,\,d\zeta\\[16pt]
\qquad+\;\Theta_{1+}^-(\lambda)\Theta_{2+}^-(\lambda)-
\Theta_{1}^+(\lambda)^{-1}\Theta_{2}^+(\lambda)^{-1}\,\\[10pt]
\qquad=\displaystyle\frac{1}{2\pi
i}\int\limits_{{C}}\frac{\Theta_{1}^-(\zeta)\Theta_{2}^-(\zeta)}
{\zeta-{\lambda}}\,\,d\zeta-\Theta_{1+}^-(\lambda)\Theta_{2+}^-(\lambda)
+\Theta_{1+}^-(\lambda)\Theta_{2+}^-(\lambda)-\Theta_{1}^+(\lambda)^{-1}\Theta_{2}^+(\lambda)^{-1}\,\\[14pt]
\qquad=\displaystyle\frac{1}{2\pi i}\int\limits_{{C}}\frac{\Theta_{21}^-(\zeta)}
{\zeta-{\lambda}}\,\,d\zeta-\Theta_{21}^+(\lambda)^{-1}
=\Theta_{21}^-(\lambda)-\Theta_{21}^+(\lambda)^{-1}\,.\\[14pt]
\end{array}
\]
\vskip 4pt \noindent Thus, for $\,\lambda\in (G_+\cup G_-)\cap\rho(T_{21})\,$, we obtain\\
\[
M_{21}(T_{21}-\lambda)^{-1}N_{21}=\Upsilon_{21}(\lambda)=
\left\{%
\begin{array}{ll}
    \Theta_{21+}^-(\lambda)-\Theta_{21}^+(\lambda)^{-1}, & \lambda\in G_+ \cap\rho(T_{21})\\
    -\Theta_{21-}^-(\lambda), & \lambda\in G_- \\
\end{array}%
\right.\,.
\]
\vskip 4pt \noindent That is
$\;\Upsilon_{21}=\mathcal{F}_{tc}(\Theta_{2}\Theta_{1})\,$. \proofend

Thus we have obtained important properties of product of systems. But the main
question  whether the product $\,\Sigma_2\cdot\Sigma_1\,$ of conservative curved
systems $\,\Sigma_1\,,\,\Sigma_2\,$ is a conservative curved system too leaves
unexplained. The following Proposition answers this question. It also answers a
question about author's motivation of the definition (Prod): in fact, the connection
between the product of systems and the product of models established in the
Proposition sheds genuine light on our definition (Prod).

\begin{prop}\label{M_and_S_con}
Suppose $\,\Pi_1,\Pi_2\in\Mod\,$, $\,\Pi=\Pi_2\cdot\Pi_1\,$,
$\,\Sigma_1=\mathcal{F}_{sm}(\Pi_1)\,$, $\,\Sigma_2=\mathcal{F}_{sm}(\Pi_2)\,$,
$\,\Sigma_{21}=\Sigma_{2}\cdot\Sigma_{1}\,$, and $\,\widehat{\Sigma}=\mathcal{F}_{sm}(\Pi)\,$. Then
$\;\;\Sigma_{21}\sim\widehat{\Sigma}\,$.
\end{prop}
\noindent We hope that it will cause no confusion if we use the same symbol
$\,\mathcal{F}_{ms}\,$ for the transformations $\,\mathcal{F}_{ms} : \Mod \to \Sys\,$
and $\,\mathcal{F}_{ms} : \Mod_n \to \Sys\,$: the latter one is defined by (MtoS) as
well (with $\,\pi_+=\pi_1\,$ and $\,\pi_-=\pi_n\,$).

\vskip 4pt \myproof\, Let
\[
\Sigma_1=(T_1, M_1, N_1)\,,\quad \Sigma_2=(T_2, M_2, N_2)\,,\quad \Sigma_{21}=(T_{21}, M_{21},
N_{21})
\]
Let also $\,\Pi=(\pi_1, \pi_2, \pi_3)\,$, $\,\widehat{\Sigma}=(\widehat{T}, \widehat{M},
\widehat{N})=\mathcal{F}_{sm}(\Pi)\,$ and
\[
\widehat{\Sigma}_1=(\widehat{T}_1, \widehat{M}_1, \widehat{N}_1)=\mathcal{F}_{sm}(\pi_1,
\pi_2)\,,\quad \widehat{\Sigma}_2=(\widehat{T}_2,
\widehat{M}_2,\widehat{N}_2)=\mathcal{F}_{sm}(\pi_2, \pi_3)\,.
\]
It is obvious that the systems $\,\widehat{\Sigma}_k\,$ and $\,{\Sigma}_k\,$, $\,k=1,2\,$ are
unitarily equivalent.

\vskip 4pt Since there are no simple and convenient expressions for operators
$\,\widehat{T}^*, \widehat{M}^*, \widehat{N}^*\,$ in terms of the model $\,\Pi\,$, we
need to employ the dual model $\,\Pi_*=(\pi_{*+},\pi_{*-})\,,\;\pi_{*\mp}\in {\mathcal
L}(L^2(\bar{C},{\mathfrak N_{\mp}}),{\mathcal H})\,$ defined by the conditions
$\,(f,\pi_{*\mp}v)_{\mathcal{H}}=<\pi_{\pm}^{\dag}f,v>_C\,, \;f\in\mathcal{H}\,,\;v\in
L^2(\bar{C},\mathfrak{N}_{\mp})\,$, where
\[
<u,v>_C:=\frac{1}{2\pi i}\int_C\;(u(z),v({\bar z)})_{\mathfrak N}\,dz\;,\;u\in L^2(C,{\mathfrak
N}),\;v\in L^2({\bar C},{\mathfrak N})\,.
\]
Then we can define the dual objects $\,\widehat{T}_*, \widehat{M}_*, \widehat{N}_*\,$ corresponding
to the subspace $\,\mathcal{K}_{*\Theta}=\Ran P_{*\Theta}\subset{\mathcal H}\,$. Note that
$\,P_{*\Theta}=P_{\Theta}^*\,$ and $\,(\widehat{T}_*, \widehat{M}_*,
\widehat{N}_*)\sim(\widehat{T}^*, \widehat{N}^*, \widehat{M}^*)\,$.

Since $\,P_{\Theta}^*\ne P_{\Theta}\,$, we have $\,\mathcal{K}_{\Theta}\ne
\mathcal{K}_{*\Theta}\,$. Besides, as is known from Section 1,
$\,\mathcal{K}_{\Theta}=\mathcal{K}_{(31)}\ne
\mathcal{K}_{(32)}\dot+\mathcal{K}_{(21)}\,$ and therefore the main challenge of the
Proposition is to handle all these subspaces coordinately. In~\cite{T1,T2}, the author
noticed that it was convenient to use the the pair of operators
$\,W,\,W_{*}\in\mathcal{L} (H,\,\mathcal{H})\,$ for a model and the dual one
simultaneously. We extend this construction to 3-models. By~\cite{T1,T2}, there exist
operators $\,W_{k},\,W_{*k}\in\mathcal{L} (H_{k},\,\mathcal{H})\,,\;k=1,2\,$ such that
$\,W_{*k}^*W_{k}=I\,$, $\,W_{k}W_{*k}^*=P_k\,$, and
\[
\begin{array}{lll}
{\widehat T}_{k}W_{k}=W_{k}T_{k}\,,\qquad &{\widehat M}_{k}W_{k}=M_{k}\,,\qquad &
{\widehat N}_{k}=W_{k}N_{k}\,,\\[2pt]
{\widehat T}_{*k}W_{*k}=W_{*k}T_{k}^*\,,\qquad &{\widehat M}_{*k}W_{*k}=N_{k}^*\,,\qquad &
{\widehat N}_{*k}=W_{*k}M_{k}^*\,,
\end{array}
\]
where $\,P_1=P_{(21)},\,P_2=P_{(32)}\,$ are projections related to the 3-model $\,\Pi=(\pi_1,
\pi_2, \pi_3)\,$. Define $\,W_{21}:=(W_{1},W_{2})\,$ and $\,W_{*21}:=(W_{*1},W_{*2})\,$. By
Lemma~\ref{lemm_mod2}, $\,P_{(32)}P_{(21)}=P_{(21)}P_{(32)}=0\,$. This implies
\[
W_{*21}^*W_{21}=\,{\rm diag}\;(I,I)\qquad\mbox{and}\qquad W_{21}W_{*21}^*=P_{(21)}+P_{(32)}\,.
\]
\noindent We put
\[
\begin{array}{lll}
{\widehat T}_{21}'=W_{21}T_{21}W_{*21}^*\,,\qquad &{\widehat M}_{21}'=M_{21}W_{*21}^*\,,\qquad &
{\widehat N}_{21}'=W_{21}N_{21}\,,\\[2pt]
{\widehat T}_{*21}'=W_{*21}T_{21}^*W_{21}^*\,,\qquad &{\widehat M}_{*21}'=N_{21}^*W_{21}^*\,,\qquad
& {\widehat N}_{*21}'=W_{*21}M_{21}^*\,,
\end{array}
\]
and (see Remark after Prop.~\ref{dir_sum})
\[
\begin{array}{lll}
{\widehat T}_{21}=(P_{(21)}+P_{(32)})\,\widehat{T}P_{(31)}\,,\, &{\widehat
M}_{21}=\widehat{M}P_{(31)}\,,\,&
{\widehat N}_{21}=(P_{(21)}+P_{(32)})\,\widehat{N}\,,\\[2pt]
{\widehat T}_{*21}=(P_{(21)}^*+P_{(32)}^*)\,\widehat{T}_{*}P_{(31)}^*\,,\, &{\widehat
M}_{*21}=\widehat{M}_{*}P_{(31)}^*\,,\,& {\widehat
N}_{21}=(P_{(21)}^*+P_{(32)}^*)\,\widehat{N}_{*}\,.
\end{array}
\]
Our aim is to show that $\,(\widehat{T}_{21}', \widehat{M}_{21}',
\widehat{N}_{21}')=(\widehat{T}_{21}, \widehat{M}_{21}, \widehat{N}_{21})\,$. If this identity is
hold, we get
\[
{\widehat T}W=WT_{21}\,,\qquad {\widehat M}W=M_{21}\,,\qquad {\widehat N}=WN_{21}\,,
\]
where $\;W=P_{(31)}W_{21}\,,\,W_*=P_{(31)}^*W_{*21}\,$. Thus,
$\;\Sigma_{21}\,{\stackrel{W}{\sim}}\,\widehat{\Sigma}\,$ and the Proposition is proved. Note also
that $\,W_{*}^*W=I\,$ and $\,WW_{*}^*=P_{(31)}\,$.

We check the desired identities by computations within the functional model. The identities
\[
{\widehat T}_{21}'=\left(%
\begin{array}{cc}
  {\widehat T}_{1} & {\widehat N}_{1}{\widehat M}_{2} \\
  0 & {\widehat T}_{2} \\
\end{array}%
\right),\; {\,\widehat M}_{21}'=(\widehat {M}_1,\displaystyle\frac{-1}{2\pi
i}\int\limits_C\Theta_1^-(z)[\widehat{M}_2(\widehat{T}_2-\zeta)^{-1}\,f]_-(z)\,dz\,)
\]
can be obtained  by a straightforward calculation. Indeed, we have\\
\[
\begin{array}{ll}
{\widehat T}_{21}'&=\; W_{21}T_{21}W_{*21}^*= \left(
\begin{array}{cc}
  W_1\,, & W_2
\end{array}
\right)\left(
\begin{array}{cc}
  T_1 & N_1M_2 \\
  0 & T_2 \\
\end{array}
\right)\left(
\begin{array}{cc}
  W_{*1}^*  \\
  W_{*2}^*  \\
\end{array}
\right)\\[24pt]
&=\left(
\begin{array}{cc}
  W_1T_1W_{*1}^* & W_1N_1M_2W_{*2}^* \\
  0 & W_2T_2W_{*2}^* \\
\end{array}
\right)=\left(%
\begin{array}{cc}
  {\widehat T}_{1} & {\widehat N}_{1}{\widehat M}_{2} \\
  0 & {\widehat T}_{2} \\
\end{array}%
\right)
\end{array}
\]
and\\
\[
\begin{array}{ll}
{\widehat M}_{21}'f&=\; M_{21}W_{*21}^*f= \left(
\begin{array}{cc}
  M_1\,, & M_{2}^{21}
\end{array}
\right)\left(
\begin{array}{cc}
  W_{*1}^*  \\
  W_{*2}^*  \\
\end{array}
\right)f\\[16pt]
&=M_1W_{*1}^*f_{21}+M_{2}^{21}W_{*2}^*f_{32}
\\[12pt]
&=\widehat{M}_{1}f_{21}-\displaystyle\frac{1}{2\pi
i}\int\limits_C\Theta_1^-(z)\,[{M}_2({T}_2-\cdot)^{-1}\,W_{*2}^*f_{32}]_-(z)\,dz\\[12pt]
&=\widehat{M}_{1}f_{21}-\displaystyle\frac{1}{2\pi
i}\int\limits_C\Theta_1^-(z)\,[\widehat{M}_2(\widehat{T}_2-\cdot)^{-1}\,f_{32}]_-(z)\,dz\,,
\end{array}
\]
\vskip 8pt\noindent where $\,f=f_{21}+f_{32}\in\mathcal{K}_{(21)}\dot+\mathcal{K}_{(32)}\,$.\\

On the other hand, using Lemma~\ref{lemm_mod2}, Prop.~\ref{proj_prop}, and the
inclusions $\,\mathcal{U}\mathcal{D}_{1+}\subset\mathcal{D}_{1+}\,$,
$\,\mathcal{U}\mathcal{H}_{31+}\subset\mathcal{H}_{31+}\,$ , we get
\[
\begin{array}{lcl}
{\widehat T}_{21}f&=&(P_{(21)}+P_{(32)})\,\widehat{T}P_{(31)}f=
(P_{(21)}+P_{(32)})P_{(31)}\mathcal{U}P_{(31)}f=\\[4pt]
&=& (P_{(21)}+P_{(32)})P_{(31)}\mathcal{U}f=(P_{(21)}+P_{(32)})\mathcal{U}f=\\[4pt]
&=&
P_{(21)}\mathcal{U}f_{21}+P_{(32)}\mathcal{U}f_{21}+P_{(21)}\mathcal{U}f_{32}+P_{(32)}\mathcal{U}f_{32}=\\[4pt]
&=& \widehat{T}_{1}f_{21}+0+P_{(21)}\mathcal{U}f_{32}+\widehat{T}_{2}f_{32}=\\[4pt]
&=& \widehat{T}_{1}f_{21}+P_{(21)}(I-P_{(32)})\mathcal{U}f_{32}+\widehat{T}_{2}f_{32}=\\[4pt]
&=& \widehat{T}_{1}f_{21}+P_{(21)}(\mathcal{U}f_{32}-\mathcal{U}f_{32}+\pi_2\widehat{M}_{2}f_{32})+\widehat{T}_{2}f_{32}=\\[4pt]
&=& \widehat{T}_{1}f_{21}+\widehat{N}_{1}\widehat{M}_{2}f_{32}+\widehat{T}_{2}f_{32}\,,
\end{array}
\]
where $\,f=f_{21}+f_{32}\in\mathcal{K}_{(21)}\dot+\mathcal{K}_{(32)}\,$. Thus we have $\,{\widehat
T}_{21}'={\widehat T}_{21}\,$.  Further, if we recall Lemma~\ref{lemm_mod1}, we obtain
\[
\begin{array}{lcl}
{\widehat
M}_{21}f&=&\widehat{M}P_{(31)}f=\widehat{M}f_{21}+\widehat{M}(I-\pi_1P_+\pi_1^{\dag})f_{32}=\\[4pt]
&=&\displaystyle\frac{1}{2\pi i}\int\limits_{C}(\pi_{1}^{\dag}f_{21})(z)\,dz+
\displaystyle\frac{1}{2\pi i}\int\limits_{C}[\pi_{1}^{\dag}(I-\pi_1P_+\pi_1^{\dag})f_{32}](z)\,dz=\\[4pt]
&=&\widehat{M}_{1}f_{21}+\displaystyle\frac{1}{2\pi
i}\int\limits_{C}(P_{-}\pi_{1}^{\dag}f_{32})(z)\,dz=
\widehat{M}_{1}f_{21}+\displaystyle\frac{1}{2\pi i}\int\limits_{C}(\pi_{1}^{\dag}f_{32})(z)\,dz=\\[4pt]
&=&\widehat{M}_{1}f_{21}+\displaystyle\frac{1}{2\pi
i}\int\limits_{C}(\pi_{1}^{\dag}\pi_{2}\pi_{2}^{\dag}f_{32})(z)\,dz+ \displaystyle\frac{1}{2\pi
i}\int\limits_{C}(\pi_{1}^{\dag}(I-\pi_{2}\pi_{2}^{\dag})f_{32})(z)\,dz=\\[4pt]
&=&\widehat{M}_{1}f_{21}+\displaystyle\frac{1}{2\pi
i}\int\limits_{C}\Theta_1^-(z)(\pi_{2}^{\dag}f_{32})(z)\,dz+0=\\[4pt]
&=& \widehat{M}_{1}f_{21}-\displaystyle\frac{1}{2\pi
i}\int\limits_C\Theta_1^-(z)\,[\widehat{M}_2(\widehat{T}_2-\zeta)^{-1}\,f_{32}]_-(z)\,dz\,.
\end{array}
\]
Therefore, $\,{\widehat M}_{21}'={\widehat M}_{21}\,$. Similarly, $\,{\widehat
M}_{*21}'={\widehat M}_{*21}\,$. We can obtain the residuary identity $\,{\widehat
N}_{21}'={\widehat N}_{21}\,$ if we make use of the duality relations
\[
({\widehat M}_{21}'f',n)=(f',{\widehat N}_{*21}'n)\,,\qquad ({\widehat
N}_{21}'m,g')=(m,{\widehat M}_{*21}'g')
\]
and
\[
({\widehat M}_{21}f,n)=(f,{\widehat N}_{*21}n)\,,\qquad ({\widehat
N}_{21}m,g)=(m,{\widehat M}_{*21}g)\,,\\[4pt]
\]
where $\,f'\in\mathcal{K}_{[31]}=\mathcal{K}_{(32)}\dot+\mathcal{K}_{(21)}\,$,
$\,g'\in\mathcal{K}_{*[31]}=\mathcal{K}_{*(32)}\dot+\mathcal{K}_{*(21)}\,$,
$\,f\in\mathcal{K}_{(31)}\,$, $\,g\in\mathcal{K}_{*(31)}\,$, $\,n\in \mathfrak{N}_1
\,$, and $\,m\in \mathfrak{N}_3 \,$. Therefore, using relations of duality, we have
\[
({\widehat N}_{21}'m,g)=(m,{\widehat M}_{*21}'g)=(m,{\widehat M}_{*21}g)=({\widehat
N}_{21}m,g)\,.\qquad\square
\]

\vskip 8pt\myrem\,\, Note that we do not claim that
$\mathcal{F}_{sm}(\Pi_2)\mathcal{F}_{sm}(\Pi_1)=\mathcal{F}_{sm}(\Pi_2\Pi_1)$. The statement and
proof of Prop.~\ref{M_and_S_con} is a good illustration to our remark that the linear similarity
(but not
unitary equivalence) is the natural kind of equivalence for conservative curved systems.\\

\noindent The following theorem is a direct consequence of Prop.~\ref{M_and_S_con}. We
shall use the notation $\,\mathcal{F}_{sc}:=\mathcal{F}_{sm}\circ\mathcal{F}_{mc}\,$.
\begin{theorema}
Let $\,\widehat{\Sigma}_1=\mathcal{F}_{sc}(\Theta_1)\,$,
$\,\widehat{\Sigma}_2=\mathcal{F}_{sc}(\Theta_2)\,$ and
$\,\widehat{\Sigma}_{21}=\mathcal{F}_{sc}(\Theta_{21})\,$, where
$\,\Theta_1,\Theta_2;\;\Theta_{21}=\Theta_2\Theta_1 \in \Cfn\,$. Suppose that
$\,\Sigma_1\sim(\widehat{\Sigma}_1\oplus\Sigma_{1u})\,$ and
$\,\Sigma_2\sim(\widehat{\Sigma}_2\oplus\Sigma_{2u})\,$, where the systems $\,\Sigma_{1u}\,$ and
$\,\Sigma_{2u}\,$ are ``purely normal'' systems. Then there exists a ``purely normal'' system
$\,\Sigma_u=(T_u,0,0,0)\,$ such that $\;\Sigma_2\cdot\Sigma_1\sim
(\widehat{\Sigma}_{21}\oplus\Sigma_u)\,$.
\end{theorema}
\proof By Prop.~\ref{M_and_S_con},
$\,\widehat{\Sigma}_2\cdot\widehat{\Sigma}_1\sim\mathcal{F}_{sc}(\Theta_2\cdot\Theta_1)
=\widehat{\Sigma}_{21}\oplus\widehat{\Sigma}_{u}\,$. Then, using Prop.~\ref{prod_sim},
we have
\[
\Sigma_2\cdot\Sigma_1\sim
(\widehat{\Sigma}_2\oplus\Sigma_{2u})\cdot(\widehat{\Sigma}_1\oplus\Sigma_{1u})\sim
(\widehat{\Sigma}_2\cdot\widehat{\Sigma}_1)\oplus\Sigma_{1u}\oplus\Sigma_{2u}
\]
and therefore $\,\Sigma_2\cdot\Sigma_1\sim(\widehat{\Sigma}_{21}\oplus\Sigma_u)\,$,
where $\,\Sigma_u=\widehat{\Sigma}_{u}\oplus\Sigma_{1u}\oplus\Sigma_{2u}\,$. \proofend

\noindent Thus we see that the definition of product of conservative curved systems
(Prod) is tightly linked to functional model though we do not refer to it explicitly.
On the other hand, its formal independence from functional model characterizes the
comparative autonomy of conservative curved system well enough. Moreover, we have
explicit formulas for  $\,\Sigma_2\cdot\Sigma_1\,$ of systems and the product depends
only on the factors $\,\Sigma_2,\,\Sigma_1\,$ and their characteristic functions
(theoretically, the dependence on characteristic functions is undesirable, but, in
author's opinion, we cannot count on having more than we have).

\vskip 2pt Now we turn to the associativity of multiplication of systems.
\begin{prop}\label{pr_ass}
One has $\;\Sigma_3\cdot(\Sigma_2\cdot\Sigma_1)\sim(\Sigma_3\cdot\Sigma_2)\cdot\Sigma_1\,$, where
$\,\Sigma_k\in\Sys\,,\,k=\overline{1,3}\,$.
\end{prop}
\myproof\, Let $\,\Sigma_k=\mathcal{F}_{sc}(\Theta_k)\,$,
$\,\Pi_k=\mathcal{F}_{mc}(\Theta_k)\,$, $\,\Pi=\Pi_3\cdot\Pi_2\cdot\Pi_1\,$ and
$\,\Sigma=\mathcal{F}_{sm}(\Pi)\,$. For the functional model $\,\Pi\,$, we consider
the following subspaces $\,\mathcal{K}_{1}=\mathcal{K}_{(21)}\,$,
$\,\mathcal{K}_{2}=\mathcal{K}_{(32)}\,$, $\,\mathcal{K}_{3}=\mathcal{K}_{(43)}\,$,
$\,\mathcal{K}_{21}=\mathcal{K}_{(31)}\,$, $\,\mathcal{K}_{32}=\mathcal{K}_{(42)}\,$,
$\,\mathcal{K}_{321}=\mathcal{K}_{(41)}\,$, and dual to them. We will denote by
$\,W_{1}\colon H_1\to\mathcal{K}_{1}\,$, $\,W_{2}\colon H_2\to\mathcal{K}_{2}\,$,
$\,W_{3}\colon H_3\to\mathcal{K}_{3}\,$ the operators that realize similarities of the
systems $\,\Sigma_{k}\,,\;k=1,2,3\,$ with the corresponding systems
$\,\widehat{\Sigma}_{k}\,$ in the model  $\,\Pi\,$ (see the proof of
Prop.~\ref{M_and_S_con}). Denote by $\,W_{*k}\,,\;k=1,2,3\,$ the dual operators. As in
the proof of Prop.~\ref{M_and_S_con}, we get that the operator
$\,W_{21}=P_{(31)}(W_1,W_2)\,$ realizes  similarity
$\,\Sigma_2\cdot\Sigma_1\sim\widehat{\Sigma}_{21}\,$. Similarly, the operator
$\,W_{32}=P_{(42)}(W_2,W_3)\,$ realizes  similarity
$\,\Sigma_3\cdot\Sigma_2\sim\widehat{\Sigma}_{32}\,$. By the same argument, we get
that the operator $\,W_{3(21)}=P_{(41)}(W_{21},W_3)\,$ realizes  similarity
$\,\Sigma_3\cdot(\Sigma_2\cdot\Sigma_1)\sim\widehat{\Sigma}_{321}\,$ and the operator
$\,W_{(32)1}=P_{(41)}(W_{1},W_{32})\,$ realizes  similarity
$\,(\Sigma_3\cdot\Sigma_2)\cdot\Sigma_1\sim\widehat{\Sigma}_{321}\,$. Thus, the
operators
\[
W_{3(21)}=P_{(41)}(P_{(31)}(W_1,W_2),W_3)\,,\quad
W_{(32)1}=P_{(41)}(W_1,P_{(42)}(W_2,W_3))\,
\]
realize the similarities
$\;\Sigma_3\cdot(\Sigma_2\cdot\Sigma_1)\sim\widehat{\Sigma}_{321}\,$ and
$\,(\Sigma_3\cdot\Sigma_2)\cdot\Sigma_1\sim\widehat{\Sigma}_{321}\,$, respectively.
Therefore,
$\;\Sigma_3\cdot(\Sigma_2\cdot\Sigma_1)\sim(\Sigma_3\cdot\Sigma_2)\cdot\Sigma_1\,$.
\proofend

Recall that the operator
$\,P_{[41]}=P_{(21)}(I-P_{(32)})(I-P_{(43)})+P_{(23)}(I-P_{(32)})+P_{(43)}
=P_{(21)}(I-P_{(43)})+P_{(32)}+P_{(43)}\,$ is a projection in $\,\mathcal{H}\,$ onto
the subspace $\,\mathcal{K}_{(21)}\dot+\mathcal{K}_{(32)}\dot+\mathcal{K}_{(43)}\,$
and its components $\,P_{(21)}(I-P_{(43)})\,$, $\,P_{(32)}\,$, $\,P_{(43)}\,$ are
commuting projections onto the subspaces $\,\mathcal{K}_{(21)}\,$,
$\,\mathcal{K}_{(32)}\,$, $\,\mathcal{K}_{(43)}\,$, respectively. Then we have
\[
\begin{array}{ll}
W_{3(21)}&=P_{(41)}(P_{(31)}(W_1,W_2),W_3)\\[4pt]
&=P_{(41)}(P_{(31)}(P_{(21)}(I-P_{(43)})+P_{(32)})+P_{(43)})(W_1,W_2,W_3)
\end{array}
\]
and
\[
\begin{array}{ll}
W_{(32)1}&=P_{(41)}(W_1,P_{(42)}(W_2,W_3))\\[4pt]
&=P_{(41)}(P_{(21)}(I-P_{(43)})+P_{(42)}(P_{(32)})+P_{(43)}))(W_1,W_2,W_3)\,.
\end{array}
\]
Thus, $\,W_{3(21)}=Y\,(W_1,W_2,W_3)\,$ and $\,W_{(32)1}=Z\,(W_1,W_2,W_3)\,$, where
$\,Y,Z :
\mathcal{K}_{(21)}\dot+\mathcal{K}_{(32)}\dot+\mathcal{K}_{(43)}\to\mathcal{K}_{(41)}\,$
and
\[
\begin{array}{l}
Y\;=\;P_{(41)}(P_{(31)}(P_{(21)}(I-P_{(43)})+P_{(32)})+P_{(43)})\,,\\[4pt]
Z\;=\;P_{(41)}(P_{(21)}(I-P_{(43)})+P_{(42)}(P_{(32)}+P_{(43)}))\,.
\end{array}
\]
Let us show that
$\,Z^{-1}\;=[P_{(21)}+(P_{(32)}+P_{(43)})P_{(42)}]|\mathcal{K}_{(41)}\;\,$. Indeed,
using Lemma~\ref{lemm_mod2} and Corollary of Lemma~\ref{lem_iii}, we have
\[
\begin{array}{l}
(P_{(21)}+(P_{(32)}+P_{(43)})P_{(42)})Z\\[6pt]
\qquad\quad=(P_{(21)}+(P_{(32)}+P_{(43)})P_{(42)})
P_{(41)}(P_{(21)}(I-P_{(43)})+P_{(42)}(P_{(32)}+P_{(43)}))\\[6pt]
\qquad\quad=(P_{(21)}P_{(41)}+(P_{(32)}+P_{(43)})P_{(42)})
(P_{(21)}(I-P_{(43)})+P_{(42)}(P_{(32)}+P_{(43)}))\\[6pt]
\qquad\quad=P_{(21)}\underline{P_{(41)}P_{(21)}}(I-P_{(43)})+P_{(21)}P_{(41)}P_{(42)}(P_{(32)}+P_{(43)})\\[6pt]
\qquad\quad+(P_{(32)}+P_{(43)})\underline{P_{(42)}P_{(21)}}(I-P_{(43)})+
\underline{(P_{(32)}+P_{(43)})P_{(42)}(P_{(32)}+P_{(43)})}\\[6pt]
\qquad\quad=P_{(21)}(I-P_{(43)})+P_{(21)}P_{(41)}P_{(42)}(P_{(32)}+P_{(43)})+(P_{(32)}+P_{(43)})\\[6pt]
\qquad\quad=I+P_{(21)}P_{(41)}P_{(42)}(P_{(32)}+P_{(43)})\\[6pt]
\end{array}
\]
Since
\[
\begin{array}{ll}
P_{(21)}P_{(41)}P_{(42)}&=(P_{(21)}+P_{(42)})P_{(41)}P_{(42)}-\underline{P_{(42)}P_{(41)}P_{(42)}}\\[6pt]
&=(P_{(21)}+P_{(42)})P_{(41)}P_{(42)}-P_{(42)}\\[6pt]
&=\underline{(P_{(21)}+P_{(42)})P_{(41)}(P_{(21)}+P_{(42)})}-(P_{(21)}+P_{(42)})P_{(41)}P_{(21)}-P_{(42)}\\[6pt]
&=(P_{(21)}+P_{(42)})-(P_{(21)}+P_{(42)})\underline{P_{(41)}P_{(21)}}-P_{(42)}\\[6pt]
&=P_{(21)}+P_{(42)}-P_{(21)}-P_{(42)}=0\,,\\[6pt]
\end{array}
\]
we obtain
$\,(P_{(21)}+(P_{(32)}+P_{(43)})P_{(42)})Z=I|\mathcal{K}_{(21)}\dot+\mathcal{K}_{(32)}\dot+\mathcal{K}_{(43)}\,.
\,$ On the other hand, we have
\[
\begin{array}{l}
Z(P_{(21)}+(P_{(32)}+P_{(43)})P_{(42)})\\[6pt]
\qquad\quad=P_{(41)}(P_{(21)}(I-P_{(43)})+P_{(42)}(P_{(32)}+P_{(43)}))
(P_{(21)}+(P_{(32)}+P_{(43)})P_{(42)})\\[6pt]
\qquad\quad=P_{(41)}(P_{(21)}+\underline{P_{(42)}(P_{(32)}+P_{(43)})P_{(42)}})
=P_{(41)}(P_{(21)}+P_{(42)})|\mathcal{K}_{(41)}\\[6pt]
\qquad\quad=\underline{P_{(41)}(P_{(21)}+P_{(42)})P_{(41)}}|\mathcal{K}_{(41)}
=P_{(41)}|\mathcal{K}_{(41)}=I|\mathcal{K}_{(41)}\,.\\[6pt]
\end{array}
\]
We need to compute the operator $\,Z^{-1}Y\,$ :
\[
\begin{array}{ll}
Z^{-1}Y&=(P_{(21)}+(P_{(32)}+P_{(43)})P_{(42)})P_{(41)}(P_{(31)}(P_{(21)}(I-P_{(43)})+P_{(32)})+P_{(43)})\\[6pt]
&=(P_{(21)}+(P_{(32)}+P_{(43)})P_{(42)})(P_{(31)}(P_{(21)}(I-P_{(43)})+P_{(32)})+P_{(41)}P_{(43)})\\[6pt]
&=P_{(21)}P_{(31)}(P_{(21)}(I-P_{(43)})+P_{(32)})+P_{(21)}P_{(41)}P_{(43)}\\[6pt]
&+(P_{(32)}+P_{(43)})P_{(42)}P_{(31)}(P_{(21)}(I-P_{(43)})+P_{(32)})+(P_{(32)}+P_{(43)})P_{(42)}P_{(41)}P_{(43)}\\[6pt]
&=P_{(21)}(I-P_{(43)})+P_{(21)}P_{(31)}P_{(32)}+P_{(21)}P_{(41)}P_{(43)}\\[6pt]
&+P_{(32)}P_{(42)}P_{(31)}(P_{(21)}(I-P_{(43)})+P_{(32)})\\[6pt]
&+\underline{P_{(43)}P_{(31)}}(P_{(21)}(I-P_{(43)})+P_{(32)})+(P_{(32)}+P_{(43)})P_{(42)}P_{(43)}\\[6pt]
\end{array}
\]
\[
\begin{array}{ll}
&=P_{(21)}(I-P_{(43)})+\underline{P_{(21)}P_{(31)}P_{(32)}}+P_{(21)}P_{(41)}P_{(43)}\\[6pt]
&+P_{(32)}\underline{P_{(42)}P_{(21)}}(I-P_{(43)})+P_{(32)}P_{(42)}P_{(31)}P_{(32)}
+\underline{P_{(32)}P_{(42)}P_{(43)}}+P_{(43)}\,.\\[6pt]
\end{array}
\]
Since
\[
\begin{array}{ll}
P_{(21)}P_{(31)}P_{(32)}&=\underline{(P_{(21)}+P_{(32)})P_{(31)}(P_{(21)}+P_{(32)})}
-P_{(21)}\underline{P_{(31)}P_{(21)}}\\[6pt]
&-\underline{P_{(32)}P_{(31)}}P_{(21)}-\underline{P_{(32)}P_{(31)}}P_{(32)}\\[6pt]
&=P_{(21)}+P_{(32)}-P_{(21)}-\underline{P_{(32)}P_{(21)}}-P_{(32)}=0
\end{array}
\]
and
\[
\begin{array}{ll}
P_{(32)}P_{(42)}P_{(43)}&=\underline{(P_{(32)}+P_{(43)})P_{(42)}(P_{(32)}+P_{(43)})}
-P_{(32)}\underline{P_{(42)}P_{(32)}}\\[6pt]
&-\underline{P_{(43)}P_{(42)}}P_{(32)}-\underline{P_{(43)}P_{(42)}}P_{(43)}\\[6pt]
&=P_{(32)}+P_{(43)}-P_{(32)}-\underline{P_{(43)}P_{(32)}}-P_{(43)}=0\,,
\end{array}
\]
we obtain
\[
\begin{array}{ll}
Z^{-1}Y&=P_{(21)}(I-P_{(43)})+P_{(21)}P_{(41)}P_{(43)}+P_{(32)}P_{(42)}P_{(31)}P_{(32)}+P_{(43)}\,.
\end{array}
\]
Taking into account that $\,\Ran P_{(43)}\subset\mathcal{H}_{41+}\,$ and using
Lemma~\ref{lem_iii}, we have
\[
\begin{array}{ll}
P_{(21)}P_{(41)}P_{(43)}&=\underline{(P_{(21)}+P_{(42)})P_{(41)}}P_{(43)}-\underline{P_{(42)}P_{(41)}}P_{(43)}\\[6pt]
&=(P_{(21)}+P_{(42)})P_{(43)}-P_{(42)}P_{(43)}=P_{(21)}P_{(43)}\,.
\end{array}
\]
Likewise, taking into account that $\,\Ran P_{(32)}\subset\mathcal{H}_{41+}\,$ and
using Lemma~\ref{lem_iii}, we have
\[
\begin{array}{ll}
P_{(32)}P_{(42)}P_{(31)}P_{(32)}&=P_{(32)}\underline{(P_{(21)}+P_{(42)})(P_{(31)}+P_{(43}))}P_{(32)}\\[6pt]
&=P_{(32)}(P_{(21)}+P_{(42)})P_{(32)}\\[6pt]
&=\underline{P_{(32)}P_{(21)}}P_{(32)}+P_{(32)}\underline{P_{(42)}P_{(32)}}=P_{(32)}P_{(32)}=P_{(32)}\,.
\end{array}
\]
Therefore,
\[
\begin{array}{ll}
Z^{-1}Y&=P_{(21)}(I-P_{(43)})+P_{(21)}P_{(43)}+P_{(32)}+P_{(43)}\\[6pt]
&=(P_{(21)}(I-P_{(43)})+P_{(32)}+P_{(43)})+P_{(21)}P_{(43)}\\[6pt]
&=(I+P_{(21)}P_{(43)})|\mathcal{K}_{(21)}\dot+\mathcal{K}_{(32)}\dot+\mathcal{K}_{(43)}\,.
\end{array}
\]
Note also that
$\,Y^{-1}Z=(I-P_{(21)}P_{(43)})|\mathcal{K}_{(21)}\dot+\mathcal{K}_{(32)}\dot+\mathcal{K}_{(43)}\,$.
Indeed, we have $\,P_{(21)}\underline{P_{(43)}P_{(21)}}P_{(43)}=0\,$.\\

Further, for the operator $\,(W_1,\,W_2,\,W_3)\, \colon H_1 \oplus H_2 \oplus H_3 \to
\mathcal{K}_{(21)}\dot+\mathcal{K}_{(32)}\dot+\mathcal{K}_{(43)}\,$, it can easily be
checked that $\,(W_1,\,W_2,\,W_3)^{-1}=\left(%
\begin{array}{ccc}
  W_{*1}^*P_{(21)}(I-P_{(43)}) \\
  W_{*2}^*P_{(32)} \\
  W_{*3}^*P_{(43)} \\
\end{array}%
\right)\,$. Thus we obtain

\[
\begin{array}{ll}
W_{(32)1}^{-1}W_{3(21)}&=(W_1,\,W_2,\,W_3)^{-1}Z^{-1}Y(W_1,\,W_2,\,W_3)\\[8pt]
&=\left(%
\begin{array}{ccc}
  I & 0 & W_{*1}^*P_{(21)}P_{(43)}W_3 \\
  0 & I & 0 \\
  0 & 0 & I \\
\end{array}%
\right)
\end{array}%
\]
and therefore $\,W_{(32)1}^{-1}W_{3(21)}\ne I\,$ because, in general, we have no the property
$\,P_{(21)}P_{(43)}=0\,$. Thus the identity
$\;\Sigma_3\cdot(\Sigma_2\cdot\Sigma_1)=(\Sigma_3\cdot\Sigma_2)\cdot\Sigma_1\,$ does not hold.\\

\exmpl\,\, We continue the example from Section 1. Consider the systems
\[
\Sigma_1=\Sigma_2=\mathcal{F}_{sc}(\theta)=((0),(1),(1)),\;\Sigma_3=\mathcal{F}_{sc}(\theta^2)=
\left(\left(%
\begin{array}{cc}
  0 & 1 \\
  0 & 0 \\
\end{array}%
\right),(1, 0),
\left(%
\begin{array}{c}
  2\varepsilon \\
  1 \\
\end{array}%
\right)\right).
\]
Then we can easily calculate that
\[
\Sigma_3\cdot(\Sigma_2\cdot\Sigma_1)=\left(\left(%
\begin{array}{cccc}
  0 & 1 & \varepsilon & 0 \\
  0 & 0 & 1 & 0 \\
  0 & 0 & 0 & 1 \\
  0 & 0 & 0 & 0 \\
\end{array}%
\right)\,, (1,\varepsilon,-\varepsilon^2,2\varepsilon^3)\,,
\left(%
\begin{array}{c}
  \varepsilon^3 \\
  -\varepsilon^2 \\
  2\varepsilon \\
  1 \\
\end{array}%
\right)
\right)
\]
and
\[
(\Sigma_3\cdot\Sigma_2)\cdot\Sigma_1=\left(\left(%
\begin{array}{cccc}
  0 & 1 & \varepsilon & -\varepsilon^2 \\
  0 & 0 & 1 & 0 \\
  0 & 0 & 0 & 1 \\
  0 & 0 & 0 & 0 \\
\end{array}%
\right)\,, (1,\varepsilon,0,0)\,,
\left(%
\begin{array}{c}
  \varepsilon^3 \\
  -\varepsilon^2 \\
  2\varepsilon \\
  1 \\
\end{array}%
\right) \right)\;.
\]
Thus, $\;\Sigma_3\cdot(\Sigma_2\cdot\Sigma_1)\ne(\Sigma_3\cdot\Sigma_2)\cdot\Sigma_1\,$. The matrix
\[
X=\left(%
\begin{array}{cccc}
  1 & 0 & -\varepsilon^2 & 2\varepsilon^3 \\
  0 & 1 & 0 & 0 \\
  0 & 0 & 1 & 0 \\
  0 & 0 & 0 & 1 \\
\end{array}%
\right)\qquad\qquad\quad
\]
realizes the similarity
$\;\Sigma_3\cdot(\Sigma_2\cdot\Sigma_1)\sim(\Sigma_3\cdot\Sigma_2)\cdot\Sigma_1\,$.

\section{Regular factorizations}

We start with extension of the notion of regularity to the class of n-characteristic functions.
\begin{defin}
We shall say that a n-characteristic function $\,\Theta\,$ is regular (and write
$\,\Theta\in\Cfn_n^{reg}\,$) if $\;\forall \; i\ge j \ge k\,$
\[
\;\;\Ran\,(I-\Theta^{\dag}_{ij}(z)\Theta_{ij}(z))^{1/2}\cap
\Ran\,(I-\Theta_{jk}(z)\Theta^{\dag}_{jk}(z))^{1/2}=\{0\}\,,\;\;a.e.\;\;z\in C\,.
\]
\end{defin}
\noindent Note that it suffices to check these conditions for
$\;k=1,\;i=n\;,j=\overline{1,n}\,$ (it follows  from~\cite{NF}, Lemma VII.4.1).

Taking into account the fact that we identify n-characteristic function with the factorization of
Schur class function, we obtain the definition of regularity for factorization of Schur class
function in the case when $\,n=3\,$. If additionally $\,\Xi_k\equiv 1,\;k=\overline{1,3}\,$, we
arrive at the standard definition ~\cite{NF,Br} (see the Introduction).

\vskip 6pt In the context of functional models the corresponding notion is the
following.
\begin{defin} Let
$\, \Mod_n^{reg}:=\{ \Pi\in\Mod_n : \Ran\pi_1\vee\Ran\pi_n=\mathcal{H} \}\,$. We shall
say that an n-model $\,\Pi\in\Mod_n\,$ is regular if $\,\Pi\in\Mod_n^{reg}\,$.
\end{defin}
We are going to show that these two notions of regularity (for n-characteristic functions and for
n-models) agree. With that end in mind we employ the construction of Prop.~\ref{CtoM}. It is easy
to show that for any two contractive operators $\,A_{21} \colon \mathcal{N}_1\to\mathcal{N}_2\,$
and $\,A_{32} \colon \mathcal{N}_2\to\mathcal{N}_3\,$ there exist three isometries $\,V_1 \colon
\mathcal{N}_1\to\mathcal{H}\,$, $\,V_2 \colon \mathcal{N}_2\to\mathcal{H}\,$, and  $\,V_3 \colon
\mathcal{N}_3\to\mathcal{H}\,$ such that
\[
A_{21}=V_2^*V_1\,,\qquad A_{32}=V_3^*V_2\,,\qquad A_{32}A_{21}=V_3^*V_1\,.
\]
Note that we need no to assume as in Prop.~\ref{CtoM} that the operators $\,A_{21}\,$
and $\,A_{32}\,$ are operator valued functions of weighted Schur class: it suffices to
assume that they are merely contractive operators. Evolving this approach, we obtain
the following Lemmas~\ref{lem_fact} and~\ref{lem_reg}. Let $\,V_1,\,V_2,\,V_3\,$ be
isometries; $\,A_{21}=V_2^*V_1\,,\,A_{32}=V_3^*V_2\,,\,A_{31}=V_3^*V_1\,$;
$\,\mathcal{E}_1=\Ran V_1\,$, $\,\mathcal{E}_2=\Ran V_2\,$, and $\,\mathcal{E}_3=\Ran
V_3\,$.
\begin{lemma}\label{lem_fact}
The following conditions are equivalent:\\[3pt]
1)\,\,$\,A_{31}=A_{32}A_{21}\,$; \,\,\,2)\,\,$\,V_3^*V_1-V_3^*V_2V_2^*V_1=0\,$;\\[2pt]
\,\,\,3)\,\,$\,((\mathcal{E}_{1}\vee\mathcal{E}_{2})\ominus\mathcal{E}_{2})\,\bot\,
((\mathcal{E}_{3}\vee\mathcal{E}_{2})\ominus\mathcal{E}_{2})\,$.
\end{lemma}
\proof \,\,1)\,$\,\Longleftrightarrow\,$ 2)\, is obvious. To prove
\,\,2)\,$\,\Longleftrightarrow\,$ 3)\, we need the following
\begin{lemm}\label{lem_iv}
One has \,\,$\;\clos\Ran (I-V_2V_2^*)V_1=(\Ran V_1\vee\Ran V_2)\ominus\Ran V_2\,$.
\end{lemm}
\proof Let $\,f\in \Ran (I-V_2V_2^*)V_1\,$. Then $\,f=(I-V_2V_2^*)V_1u\in\Ran
V_1\vee\Ran V_2\,$. On the other hand,
$\,V_2^*f=V_2^*(I-V_2V_2^*)V_1u=(I-V_2^*V_2)V_2^*V_1u=0\cdot V_2^*V_1u=0\,$, that is,
$\,f\bot\Ran V_2\,$ and therefore $\,\Ran (I-V_2V_2^*)V_1\subset(\Ran V_1\vee\Ran
V_2)\ominus\Ran V_2\,$.

Conversely, let $\,f\in(\Ran V_1\vee\Ran V_2)\ominus\Ran V_2\,$. Then we have
$\,V_2^*f=0\,$ and $\,f=\lim_{n\to\infty}(V_1u_{1n}+V_2u_{2n})\,$. Hence,
\[
\begin{array}{ll}
f&=(I-V_2V_2^*)f=\lim_{n\to\infty}((I-V_2V_2^*)V_1u_{1n}+(I-V_2V_2^*)V_2u_{2n})\\[6pt]
&=(I-V_2V_2^*)f=\lim_{n\to\infty}(I-V_2V_2^*)V_1u_{1n} \in \clos\Ran
(I-V_2V_2^*)V_1\,.\quad \square
\end{array}
\]
To complete the proof of Lemma~\ref{lem_fact} we need only to make use of the
following observation
\[ \Ran (I-V_2V_2^*)V_1 \;\bot\; \Ran
(I-V_2V_2^*)V_3 \,\Longleftrightarrow\, V_3^*(I-V_2V_2^*)V_1=0\,.\quad \square
\]
\begin{lemm}\label{lem_v}
Assume that \,\,\,$\,V_3^*V_1-V_3^*V_2V_2^*V_1=0\,$. Then
\[
(\mathcal{E}_{1}\vee\mathcal{E}_{3})\ominus\mathcal{E}_{3}\subset
((\mathcal{E}_{1}\vee\mathcal{E}_{2})\ominus\mathcal{E}_{2})\,\oplus\,
((\mathcal{E}_{3}\vee\mathcal{E}_{2})\ominus\mathcal{E}_{3})\,.
\]
\end{lemm}
\proof Using Lemma~\ref{lem_iv} and the obvious identity
\[
(I-V_3V_3^*)V_1=(I-V_2V_2^*)V_1+(I-V_3V_3^*)V_2V_2^*V_1\,,
\]
we get
\[
(\mathcal{E}_{1}\vee\mathcal{E}_{3})\ominus\mathcal{E}_{3}\subset
((\mathcal{E}_{1}\vee\mathcal{E}_{2})\ominus\mathcal{E}_{2})\,\vee\,
((\mathcal{E}_{3}\vee\mathcal{E}_{2})\ominus\mathcal{E}_{3})\,.
\]By Lemma~\ref{lem_fact},
$\,((\mathcal{E}_{1}\vee\mathcal{E}_{2})\ominus\mathcal{E}_{2})\,\bot\,
((\mathcal{E}_{3}\vee\mathcal{E}_{2})\ominus\mathcal{E}_{2})\,$. Then
\[
((\mathcal{E}_{1}\vee\mathcal{E}_{2})\ominus\mathcal{E}_{2})\,\bot\,
((\mathcal{E}_{3}\vee\mathcal{E}_{2})\ominus\mathcal{E}_{2})\oplus\mathcal{E}_{2}
=\mathcal{E}_{3}\vee\mathcal{E}_{2}
\]
and therefore $\,((\mathcal{E}_{1}\vee\mathcal{E}_{2})\ominus\mathcal{E}_{2})\,\bot\,
((\mathcal{E}_{3}\vee\mathcal{E}_{2})\ominus\mathcal{E}_{3})\,$. \proofend

\myrem\,\, If we define the isometries $\,\tau_{jij} \colon
\clos\Ran\,(I-V_j^*V_iV_i^*V_j)^{1/2} \to \mathcal{H} \,$ by the formula $\,
\tau_{jij}(I-V_j^*V_iV_i^*V_j)^{1/2}=(I-V_iV_i^*)V_j\,$, we can rewrite the identity
\[
(I-V_3V_3^*)V_1=(I-V_2V_2^*)V_1+(I-V_3V_3^*)V_2V_2^*V_1
\]
in the form
\[
\tau_{131}(I-A^*_{31}A_{31})^{1/2}=\tau_{121}(I-A^*_{21}A_{21})^{1/2}+
\tau_{232}(I-A^*_{32}A_{32})^{1/2}A_{21}\,.
\]
Note that $\,\clos\Ran\tau_{jij}=\clos\Ran(I-V_iV_i^*)V_j\,$. Then, by
Lemma~\ref{lem_fact}, the condition $\,V_3^*(I-V_2V_2^*)V_1=0\,$ means
$\,\tau_{232}^*\tau_{121}=0\,$ and $\,\tau_{121}^*\tau_{232}=0\,$. Therefore we have
\[
Z\,(I-A^*_{31}A_{31})^{1/2}=\left(%
\begin{array}{c}
  (I-A^*_{21}A_{21})^{1/2} \\
  (I-A^*_{32}A_{32})^{1/2}A_{21} \\
\end{array}%
\right)\eqno{\rm{(Z)}}
\]
and the operator $\,Z=(\tau_{121}^*+\tau_{232}^*)\tau_{131}\,$ is an isometry. We need
the following Lemma established in~\cite{Br}.
\begin{lemm}\label{lem_vi}
The following conditions are equivalent:\\[3pt]
1)\, $\,\Ran\,(I-A^*_{32}A_{32})^{1/2}\cap \Ran\,(I-A_{21}A^*_{21})^{1/2}=\{0\}\,$;\\[2pt]
2)\, $\,A_{21}^*(I-A^*_{32}A_{32})^{1/2}m +(I-A^*_{21}A_{21})^{1/2}n=0\,$, $\;m\in
\clos\Ran\,(I-V_2^*V_3V_3^*V_2)^{1/2}\,$ and
$\;n\in \clos\Ran\,(I-V_1^*V_2V_2^*V_1)^{1/2}\,$ $\,\Longrightarrow\,$ $\,m=0,\,n=0\,$;\\[2pt]
\,\,\,\,\,\, 3)\, The operator $\,Z\,$ defined by the condition {\rm (Z)} is an
unitary operator.
\end{lemm}
\noindent Now we have prepared to prove the Lemma~\ref{lem_reg}, which allow us to
translate factorization problems into geometrical language and now we can point out
the purely geometrical nature of the notion of regularity. Note that this fact is the
underlying basis of the generalization of Sz.-Nagy-Foia\c{s}'s regularity criterion
in~\cite{BD}, where the authors drop the condition of analyticity.
\begin{lemma}\label{lem_reg}
Assume that \,\,\,$\,V_3^*V_1-V_3^*V_2V_2^*V_1=0\,$.\\[2pt]
Then the following conditions are equivalent:\\[3pt]
1)\, $\,\Ran\,(I-A^*_{32}A_{32})^{1/2}\cap \Ran\,(I-A_{21}A^*_{21})^{1/2}=\{0\}\,$;\\[2pt]
2)\, $\,\clos\Ran\,(I-V_3V_3^*)V_1=\clos\Ran\,(I-V_2V_2^*)V_1\oplus
\clos\Ran\,(I-V_3V_3^*)V_2\,$;\\[2pt]
\,\,\,\,\,\, 3)\, $\,\mathcal{E}_2\subset\mathcal{E}_1\vee\mathcal{E}_3\,$.
\end{lemma}
\proof \,\,1)\,$\,\Longleftrightarrow\,$ 2)\, By Lemma~\ref{lem_vi}, condition 1) is
equivalent to the condition that the operator $\,Z\,$ is unitary. Since under our
assumptions $\,Z\,$ is always isometrical, we can check only that
$\,Z^*=\tau_{131}^*(\tau_{121}+\tau_{232})\,$ is an isometrical operator. The latter
is equivalent to the condition
$\,\Ran\tau_{121}\oplus\Ran\tau_{232}\subset\Ran\tau_{131}\,$. The inverse inclusion
is Lemma~\ref{lem_v}.

\noindent \,\,2)\,$\,\Longrightarrow\,$ 3)\, Since we have
$\,(\mathcal{E}_{1}\vee\mathcal{E}_{3})\ominus\mathcal{E}_{3}=
((\mathcal{E}_{1}\vee\mathcal{E}_{2})\ominus\mathcal{E}_{2})\,\oplus\,
((\mathcal{E}_{3}\vee\mathcal{E}_{2})\ominus\mathcal{E}_{3})\,$, we obtain
\[
\begin{array}{ll}
\mathcal{E}_{1}\vee\mathcal{E}_{3}&=
\mathcal{E}_{3}\oplus((\mathcal{E}_{1}\vee\mathcal{E}_{3})\ominus\mathcal{E}_{3})
=\mathcal{E}_{3}\oplus((\mathcal{E}_{3}\vee\mathcal{E}_{2})\ominus\mathcal{E}_{3})
\,\oplus\,((\mathcal{E}_{1}\vee\mathcal{E}_{2})\ominus\mathcal{E}_{2})\\[6pt]
&=(\mathcal{E}_{3}\vee\mathcal{E}_{2})
\,\oplus\,((\mathcal{E}_{1}\vee\mathcal{E}_{2})\ominus\mathcal{E}_{2}).
\end{array}
\]
Hence, $\,\mathcal{E}_2\subset\mathcal{E}_1\vee\mathcal{E}_3\,$.

\noindent \,\,3)\,$\,\Longrightarrow\,$ 2)\, We have
$\,\mathcal{E}_2\subset\mathcal{E}_1\vee\mathcal{E}_3\,$. Then
\[
\mathcal{E}_{3}\oplus((\mathcal{E}_{2}\vee\mathcal{E}_{3})\ominus\mathcal{E}_{3})
=\mathcal{E}_{2}\vee\mathcal{E}_{3}\subset\mathcal{E}_{1}\vee\mathcal{E}_{3}
=\mathcal{E}_{3}\oplus((\mathcal{E}_{1}\vee\mathcal{E}_{3})\ominus\mathcal{E}_{3})
\]
and therefore $\,(\mathcal{E}_{2}\vee\mathcal{E}_{3})\ominus\mathcal{E}_{3}
\subset(\mathcal{E}_{1}\vee\mathcal{E}_{3})\ominus\mathcal{E}_{3}\,$.

On the other hand, we have
$\,((\mathcal{E}_{1}\vee\mathcal{E}_{2})\ominus\mathcal{E}_{2})\,\bot\,\mathcal{E}_{2})\,$
and $\,((\mathcal{E}_{1}\vee\mathcal{E}_{2})\ominus\mathcal{E}_{2})\,\bot\,
((\mathcal{E}_{3}\vee\mathcal{E}_{2})\ominus\mathcal{E}_{2})\,$. Hence,
$\,((\mathcal{E}_{1}\vee\mathcal{E}_{2})\ominus\mathcal{E}_{2})
\,\bot\,\mathcal{E}_{3}\vee\mathcal{E}_{2}\,$. Then, we get
\[
\mathcal{E}_{3}\oplus((\mathcal{E}_{1}\vee\mathcal{E}_{2})\ominus\mathcal{E}_{2})
\subset \mathcal{E}_{1}\vee\mathcal{E}_{2}\vee\mathcal{E}_{3} \subset
\mathcal{E}_{1}\vee\mathcal{E}_{3}
=\mathcal{E}_{3}\oplus((\mathcal{E}_{1}\vee\mathcal{E}_{3})\ominus\mathcal{E}_{3})
\]
and therefore $\,(\mathcal{E}_{1}\vee\mathcal{E}_{2})\ominus\mathcal{E}_{2} \subset
(\mathcal{E}_{1}\vee\mathcal{E}_{3})\ominus\mathcal{E}_{3}\,$. Thus, we obtain
\[
((\mathcal{E}_{1}\vee\mathcal{E}_{2})\ominus\mathcal{E}_{2})
\oplus((\mathcal{E}_{2}\vee\mathcal{E}_{3})\ominus\mathcal{E}_{3})
\subset(\mathcal{E}_{1}\vee\mathcal{E}_{3})\ominus\mathcal{E}_{3}\,.
\]
The inverse inclusion is Lemma~\ref{lem_v}. \proofend

\noindent The following assertion is a straightforward consequence of the lemmas.
\begin{prop}\label{reg_M_C}
One has
$\;\,\Pi=\mathcal{F}_{mc}(\Theta)\in\Mod_n^{reg}\;\Longleftrightarrow\;\Theta\in\Cfn_n^{\,reg}\,$.
\end{prop}

\vskip 3pt We have defined the notions of regularity for $\Cfn_n$ and $\Mod_n$. Now we
pass over to curved conservative systems looking for a counterpart of the regularity
in this new context. In the Introduction we have defined the notion of \textit{simple}
curved conservative systems. For them, we have
\begin{prop}\label{simpl}
1)\,\,Let $\;\widehat{\Sigma}=\mathcal{F}_{sm}(\Pi)\,,\;\Pi\in\Mod\,$ and
$\,\rho(\widehat{T})\cap G_+\ne\emptyset\;$. Then the system $\,\widehat{\Sigma}\,$ is
simple;\,\,\,2)\,\, If $\,\Sigma\sim\Sigma'\,$, then $\,\Sigma\,$ is simple
$\,\Leftrightarrow\,$ $\,\Sigma'\,$ is simple;\,\,\,3)\,\, If a system $\,\Sigma\,$ is
simple, then the system $\,\Sigma^*\,$ is also simple;\,\,\,4)\,\, If
$\,\Sigma\;{\stackrel{X}{\sim}}\;\Sigma'\,$,
$\,\Sigma\;{\stackrel{X'}{\sim}}\;\Sigma'\,$ and the system $\,\Sigma\,$ is simple,
then $\,X=X'\,$;\,\,\,5)\,\, If the system $\,\Sigma=\Sigma_2\cdot\Sigma_1\,$ is
simple and $\,\rho(T)\cap\rho(T_1)\cap G_+\ne\emptyset\;$, then the systems
$\,\Sigma_1\,$ and $\,\Sigma_2\,$ are simple.
\end{prop}
\proof \,1)\, By~\cite{T3}, we have
\[
-\widehat{M}(\widehat{T}-z)^{-1}f=\left\{\begin{array}{ll}
\Theta^+(z)^{-1}(\pi_-^{\dag}f)(z)&,\,z\in \rho(\widehat{T})\cap G_+\\
(\pi_+^{\dag}f)(z)&,\,z\in G_-
\end{array}\right.
\]
and therefore
\[
\begin{array}{ll}
\bigcap\limits_{z\in\rho(T)}\,\Ker
\widehat{M}(\widehat{T}-z)^{-1}&=\{f\in\mathcal{K}\, :
\,\pi_+^{\dag}f=0\,,\;\pi_-^{\dag}f=0\,\}\\[6pt]
&=\{f\in\mathcal{K}\, : \,f\,\bot\,\Ran\pi_+\,,\;f\,\bot\,\Ran\pi_-\,\}=\{0\}\,.
\end{array}
\]
\noindent 2)\, Let $\,\Sigma\;{\stackrel{X}{\sim}}\;\Sigma'\,$. Then
$\,M(T-z)^{-1}=M'(T'-z)^{-1}X\,$ and the property follows straightforwardly from this
identity.

\noindent 3)\, is a direct consequence of properties 1) and 2).

\noindent 4)\, It is sufficient to check that
$\,\Sigma\;{\stackrel{X}{\sim}}\;\Sigma\,\Rightarrow\, X=I\,$. We have
\[
\begin{array}{c}
M(T-z)^{-1}Xf=MX(T-z)^{-1}f=M(T-z)^{-1}f\,\Longrightarrow\,\\[8pt]
Xf-f\in \bigcap\limits_{z\in\rho(T)}\,\Ker M(T-z)^{-1}=\{0\}\,\Longrightarrow\,Xf=f\,.
\end{array}
\]
\noindent 5)\, It can easily be checked that
$\,\rho(T)\cap\rho(T_1)\subset\rho(T_2)\,$. Then we have
\[
(T_{21}-z)^{-1}=\left(
\begin{array}{cc}
  (T_1-z)^{-1} & -(T_1-z)^{-1}N_1M_2(T_2-z)^{-1} \\
  0 & (T_2-z)^{-1} \\
\end{array}
\right).
\]
and therefore $\,\forall\,f_1\in H_1\;M_{21}(T_{21}-z)^{-1}f_1=M_1(T_1-z)^{-1}f_1\,$.
Hence,
\[
\bigcap\limits_{z\in\rho(T)}\,\Ker M_1(T_1-z)^{-1} \subset
\bigcap\limits_{z\in\rho(T)}\,\Ker M_{21}(T_{21}-z)^{-1}=\{0\}\,,
\]
that is, the system $\,\Sigma_1\,$ is simple.

Further, by property 3), the system $\,\Sigma^*\,$ is simple. Then, using the same
arguments as above, it follows that the system $\,\Sigma_2^*\,$ is simple. Hence the
system $\,\Sigma_2\,$ is simple too. \proofend

\begin{defin} The product of systems $\,\Sigma_{21}=\Sigma_{2}\cdot\Sigma_{1}\,$ is called
regular if the system $\,\Sigma_{21}\,$ is simple.
\end{defin}
\begin{prop}\label{reg_M_S}
Let $\,\widehat{\Sigma}_1=\mathcal{F}_{sm}(\Pi_{1})\,$,
$\,\widehat{\Sigma}_2=\mathcal{F}_{sm}(\Pi_{2})\,$. Suppose $\,\Sigma_1\sim\widehat{\Sigma}_1\,$,
$\,\Sigma_2\sim\widehat{\Sigma}_2\,$, $\,\Sigma_{21}=\Sigma_2\cdot\Sigma_1\,$, and
$\,\rho(T_{21})\cap G_+\ne\emptyset\;$. Then the product $\,\Sigma_2\cdot\Sigma_1\,$ is regular
$\,\Leftrightarrow\,$ the product $\,\Pi_{2}\cdot\Pi_{1}\,$ is regular.
\end{prop}
\myproof\,\,Without loss of generality (see Prop.~\ref{M_and_S_con} and
Prop.~\ref{simpl}) it can be assumed that
$\,\Sigma_{21}=\widehat{\Sigma}=\mathcal{F}_{sm}(\Pi_2\cdot\Pi_1)\,$ and
$\,T_{21}=\widehat{T}\,$. As above (see the proof of Prop.~\ref{simpl}), we get
\[
\mathcal{K}_{u}=\bigcap\limits_{z\in\rho(T)}\,\Ker
\widehat{M}(\widehat{T}-z)^{-1}=\{f\in\mathcal{K}_{(31)}\, :
\,\pi_+^{\dag}f=0\,,\;\pi_-^{\dag}f=0\,\}
\]
and therefore $\,\mathcal{K}_{u}\subset(\Ran\pi_+\vee\Ran\pi_-)^{\bot}\,$. On the
hand, if $\,f\in(\Ran\pi_+\vee\Ran\pi_-)^{\bot}\,$, then $P_{(31)}f=f$ and
$\,f\in\mathcal{K}_{(31)}\,$. Thus,
$\,\mathcal{K}_{u}=(\Ran\pi_+\vee\Ran\pi_-)^{\bot}\,$. It remains to note that the
product $\,\Sigma_2\cdot\Sigma_1\,$ is regular iff $\,\mathcal{K}_{u}=\{0\}\,$ and the
product $\,\Pi_{2}\cdot\Pi_{1}\,$ is regular iff
$\,\Ran\pi_+\vee\Ran\pi_-=\mathcal{H}\,$ (recall that $\,\pi_+=\pi_3\,$ and
$\,\pi_-=\pi_1\,$).\proofend

\noindent Combining Prop.~\ref{reg_M_C} and Prop.~\ref{reg_M_S}, we arrive at
\begin{prop}[Criterion of regularity]\label{reg_C_S}
Let $\,\widehat{\Sigma}_1=\mathcal{F}_{sc}(\Theta_{1})\,$,
$\,\widehat{\Sigma}_2=\mathcal{F}_{sc}(\Theta_{2})\,$. Suppose
$\,\Sigma_1\sim\widehat{\Sigma}_1\,$, $\,\Sigma_2\sim\widehat{\Sigma}_2\,$,
$\,\Sigma_{21}=\Sigma_2\cdot\Sigma_1\,$, and $\,\rho(T_{21})\cap G_+\ne\emptyset\;$. Then the
product $\,\Sigma_2\cdot\Sigma_1\,$ is regular $\,\Leftrightarrow\,$ the factorization
$\,\Theta_{2}\cdot\Theta_{1}\,$ is regular.
\end{prop}
\noindent Thus we obtain the correspondence between regular factorizations of characteristic
functions and regular products of systems.

\vskip 3pt \myrem\,\, It can easily be shown that the inner-outer
factorization~\cite{NF} of Schur class functions is regular (see~\cite{NF}). Hence,
using the criterion of regularity, one can prove that the product of colligations with
$C_{11}$ and $C_{00}$ contractions is regular. It is possible to extend this result to
the case of weighted Schur functions employing the generalization of regularity
criterion (Prop.~\ref{reg_C_S}). Note that, for J-contractive analytic operator
functions, J-inner-outer factorization is regular too~\cite{T4}. However, since in
this situation we lose such a geometrical functional model  as is the
Sz.-Nagy-Foia\c{s} model for contractions (and such a geometrical description of
regularity), we have to establish at first the regularity of the product of
``absolutely continuous'' and ``singular'' colligations (analogous of $C_{11}$ and
$C_{00}$ contractions) and then to obtain the uniqueness of J-inner-outer
factorization~\cite{T4}.

\section{Factorizations and invariant subspaces}

The most remarkable feature of the product of systems is its connection with invariant
subspaces. We see that the subspace $H_1$ in the definition (Prod) is invariant under
the operator $\,T_{21}\,$ (and under its resolvent $\;(T_{21}-z)^{-1}\,,\;z\in G_-$).
In the context of functional model this implies that the subspace
$\,\mathcal{K}_{(21)}\,$ is invariant under the operator $\,\widehat{T}\,$ (see
Prop.~\ref{M_and_S_con}). Following B.Sz.-Nagy and C.Foia\c{s}, we shall work within
the functional model and use the model as a tool for studying the correspondence
``factorizations $\,\leftrightarrow\,$ invariant subspaces''. Let
$\,\Theta\in\Mod_3\,,\Pi=\mathcal{F}_{mc}(\Theta)=(\pi_1, \pi_2, \pi_3)\in\Mod_3\,$.
We define the transformation $\,L=\mathcal{F}_{ic}(\Theta)\,$ as a mapping that takes
each 3-characteristic function $\,\Theta\,$ (which we identify with factorization of
Schur class function) to the invariant subspace $\,L:=\mathcal{K}_{(21)}=\Ran
P_{(21)}\,$. To study the transformation $\,\mathcal{F}_{ic}\,$ (and its ingenuous
extension to n-characteris\-tic functions), we need to make some preliminary work.

Let $\,\Pi\in\Mod_n\,$. Consider the chain of subspaces
$\;\mathcal{H}_{11+}\subset\ldots\subset\mathcal{H}_{n1+}\,$ (see the definition of
$\,\mathcal{H}_{ij+}\,$ after Lemma~\ref{lemm_mod2}). These subspaces are invariant under the
resolvent $\;(\mathcal{U}-z)^{-1}\,,\;z\in G_-$. The inverse is also true accurate up to the
``normal'' part of the chain.

\begin{prop}\label{u_chain}
Suppose $\,\mathcal{U}\in\mathcal{L}(\mathcal{H})\,$ is a normal operator,
$\,\sigma(\mathcal{U})\subset C$, and $\,\mathcal{H}_{1+}\subset\ldots\subset\mathcal{H}_{n+}\,$ is
a chain of invariant under $\;(\mathcal{U}-z)^{-1}\,,\;z\in G_-$ subspaces. Then there exists an
n-model $\,\Pi\in\Mod_n\,$ such that
$\,\mathcal{H}_{k1+}\subset\mathcal{H}_{k+}\,,\;k=\overline{1,n}\,$ and the subspaces
$\,\mathcal{H}_{uk}:=\mathcal{H}_{k+}\ominus\mathcal{H}_{k1+}\,$ reduce the operator
$\,\mathcal{U}\,$. If an n-model $\,\Pi'\in\Mod_n\,$ satisfies the same conditions, then
$\,\mathcal{H}_{k+}'=\mathcal{H}_{k+}\,$ and $\,\exists\,\psi_k\,$ such that
$\,\psi_k,\psi_k^{-1}\in H^{\infty}(G_+,\mathcal{L}(\mathfrak{N}_k))\,$ and
$\,\pi_k'=\pi_k\psi_k\,$. Besides, we have
$\,\mathcal{H}_{u1}\subset\ldots\subset\mathcal{H}_{un}\,$.
\end{prop}
\proof Consider the Wold type decomposition
$\,\mathcal{H}_{k+}=\mathcal{H}_{k+}^{pur}\oplus\mathcal{H}_{k+}^{nor}\,$  with
respect to the normal operator $\,\mathcal{U}\,,\,\sigma(\mathcal{U})\subset C\,$
(see~\cite{AD}). The operator $\,\mathcal{U}|\mathcal{H}_{k+}^{pur}\,$ is the pure
subnormal part of $\,\mathcal{U}|\mathcal{H}_{k+}\,$ and
$\,\mathcal{U}|\mathcal{H}_{k+}^{nor}\,$ is a normal operator. This decomposition is
unique. We set
\[
\mathcal{E}_{k+}=\mathcal{H}_{k+}^{pur},\qquad\mathcal{E}_{k}=\vee_{z\notin
C}(\mathcal{U}-z)^{-1}\mathcal{E}_{k+},\qquad\mathcal{E}_{k-}=\mathcal{E}_{k}\ominus\mathcal{E}_{k+}
\]
Obviously, $\,\mathcal{E}_{k-}\subset\mathcal{H}_{k+}^{\bot}\,$ and
$\,\mathcal{U}^*|\mathcal{E}_{k-}\,$ is the pure subnormal part of
$\,\mathcal{U}^*|\mathcal{H}_{k+}^{\bot}\,$. For $\,i\ge j\ge k\,$, we have
$\,\mathcal{E}_{i-}\bot\mathcal{E}_{j+}\,$, $\,\mathcal{E}_{k+}\subset\mathcal{H}_{j+}\,$, and
$\,\mathcal{E}_{i-}\subset\mathcal{H}_{j+}^{\bot}\,$. Hence,
\[
\mathcal{E}_{k}\subset\vee_{z\notin C}(\mathcal{U}-z)^{-1}\mathcal{H}_{j+}\qquad\mbox{and}\qquad
\mathcal{E}_{i}\subset\vee_{z\notin C}(\mathcal{U}^*-\bar{z})^{-1}\mathcal{H}_{j+}^{\bot}\,.
\]
This implies that
\[
\mathcal{E}_{j}\oplus((\mathcal{E}_{k}\vee\mathcal{E}_{j})\ominus\mathcal{E}_{j})
=\mathcal{E}_{k}\vee\mathcal{E}_{j}\subset \mathcal{E}_{j}\vee\mathcal{H}_{j+}
=\mathcal{E}_{j}\oplus\mathcal{H}_{j+}^{nor}\,.
\]
Therefore we get
\[
(\mathcal{E}_{k}\vee\mathcal{E}_{j})\ominus\mathcal{E}_{j}\subset
\mathcal{H}_{j}^{nor}\qquad\mbox{and}\qquad
\mathcal{E}_{j+}\oplus((\mathcal{E}_{k}\vee\mathcal{E}_{j})\ominus\mathcal{E}_{j})\subset
\mathcal{H}_{j}
\]
In the same way,
$\,\mathcal{E}_{j-}\oplus((\mathcal{E}_{i}\vee\mathcal{E}_{j})\ominus\mathcal{E}_{j})\subset
\mathcal{H}_{j}^{\bot}\,$. And finally,
\[
((\mathcal{E}_{i}\vee\mathcal{E}_{j})\ominus\mathcal{E}_{j})\,\bot\,
((\mathcal{E}_{k}\vee\mathcal{E}_{j})\ominus\mathcal{E}_{j}\,.
\]
We need to make use of the following lemma.
\begin{lem}
Suppose $\,\mathcal{U}\in\mathcal{L}(\mathcal{H})\,$ is a normal operator,
$\,\sigma(\mathcal{U})\subset C$, $\,\mathcal{E}_+\subset\mathcal{H}\,$ and
$\,\mathcal{U}|\mathcal{E}_+\,$ is a pure subnormal operator. Then there exists an
operator $\pi\in {\mathcal L}(L^2(C,{\mathfrak N}),{\mathcal H})$ such that
$\,\Ran\pi=\vee_{\lambda\notin C}(\mathcal{U}-\lambda)^{-1}\mathcal{E}_{+}\,$,
$\,\Ker\pi=\{0\}\,$, $\,\pi E^2(G_+,\mathfrak{N})=\mathcal{E}_+$ and
$\,\mathcal{U}\pi=\pi z\,$.
\end{lem}
\proof Without loss of generality we can assume that
$\,\mathcal{H}=\vee_{\lambda\notin C}(\mathcal{U}-\lambda)^{-1}\mathcal{E}_{+}\,$.
By~\cite{AD}, there exists an unitary operator $\,Y_0\in {\mathcal
L}(E_{\alpha}^2(G_+,{\mathfrak N}),{\mathcal E}_+)\,$ such that
$\,\mathcal{U}Y_0=Y_0z\,$, where $\,E_{\alpha}^2(G_+,{\mathfrak N})\,$ is the Smirnov
space of character-automorphic functions (see the comments between
Prop.\ref{proj_prop} and Prop.\ref{dir_sum}). By Mlak's lifting theorem~\cite{Ml}, the
operator $\,Y_0\,$ can be extended to the space $\,L^2(C,{\mathfrak N})\,$ lifting the
intertwining condition. This extension will be denoted by $\,\pi_0\in {\mathcal
L}(L^2(C,{\mathfrak N}),{\mathcal H})\,$. So, we have $\,\mathcal{U}\pi_0=\pi_0z\,$.
Similarly, there exists an extension $\,X_0\in {\mathcal L}({\mathcal
H},L^2(C,{\mathfrak N}))\,$ of the operator $\,Y_0^{-1}\,$ such that
$\,X_0\mathcal{U}=zX_0\,$. Thus, $\,X_0\pi_0|E_{\alpha}^2(G_+,{\mathfrak
N})=I|E_{\alpha}^2(G_+,{\mathfrak N})\,$. Since $\,L^2(C,{\mathfrak
N}))=\vee_{\lambda\notin C}(z-\lambda)^{-1}E_{\alpha}^2(G_+,{\mathfrak N})\,$, we get
$X_0\pi_0=I$. Likewise, since $\,\pi_0X_0|\mathcal{E}_+=I|\mathcal{E}_+\,$ and
$\,\mathcal{H}=\vee_{\lambda\notin C}(\mathcal{U}-\lambda)^{-1}\mathcal{E}_{+}\,$, we
get $\pi_0X_0=I$ and therefore $\,\pi_0^{-1}=X_0\in {\mathcal L}({\mathcal
H},L^2(C,{\mathfrak N}))\,$.

According to~\cite{AD}, the ``bundle'' shift $\,z|E_{\alpha}^2(G_+,\mathfrak{N})\,$ is
similar to the trivial shift $\,z|E^2(G_+,\mathfrak{N})\,$. The similarity is realized
by operator valued function $\,\chi\in L^{\infty}(C,\mathcal{L}(\mathfrak{N}))\,$ such
that $\,\chi^{-1}\in L^{\infty}(C,\mathcal{L}(\mathfrak{N}))\,$ and $\,\chi
E^2(G_+,\mathfrak{N})=E_{\alpha}^2(G_+,\mathfrak{N})\,$. Then we put
$\,\pi:=\pi_0\chi\,$. \proofend

Since $\,\mathcal{U}|\mathcal{E}_{j+}\,$ is the pure subnormal part of
$\,\mathcal{U}|\mathcal{H}_{j}\,$, there exists operators $\pi_{j}\in {\mathcal L}(L^2(C,{\mathfrak
N_{j}}),{\mathcal H})$ such that $\,\Ran\pi_j=\mathcal{E}_{j}\,$, $\,\pi_j
E^2(G_+)=\mathcal{E}_{j+}\,$, and $\,\mathcal{U}\pi_j=\pi_j z\,$. In terms of operators $\,\pi_j\,$
we rewrite the relations obtained earlier. The relation $\,\mathcal{E}_{i-}\bot\mathcal{E}_{j+}\,$
implies $\,P_{-}(\pi_{i}^{\dagger}\pi_{j})P_+=0\,$ and the orthogonality
$\,((\mathcal{E}_{i}\vee\mathcal{E}_{j})\ominus\mathcal{E}_{j})\bot
((\mathcal{E}_{k}\vee\mathcal{E}_{j})\ominus\mathcal{E}_{j})\,$ means that
$\,\Ran(I-\pi_j\pi_j^{\dag})\pi_i \,\bot\, \Ran(I-\pi_j\pi_j^{\dag})\pi_k\,$. Hence,
$\,\pi_i^{\dag}(I-\pi_j\pi_j^{\dag})\pi_k=0\,$ and
$\,\pi_i^{\dag}\pi_k=\pi_i^{\dag}\pi_j\pi_j^{\dag}\pi_k\,$. Thus, the n-tuple
$\,\Pi=(\pi_1,\ldots,\pi_n)\,$ is an n-model.

\noindent We put $\,\mathcal{H}_{j1+}=\mathcal{H}_{\pi_j\vee\dots\vee\pi_1}\cap\Ker
\pi_jP_-\pi_j^{\dag}\,$. Then,
\[
\begin{array}{lcl}
\mathcal{H}_{j1+}&=&\Ran P_{\pi_j\vee\dots\vee\pi_1}\cap\,\Ran\,(I-\pi_jP_-\pi_j^{\dag}) =\Ran
P_{\pi_j\vee\dots\vee\pi_1}(I-\pi_jP_-\pi_j^{\dag})\\[4pt]
&=&\Ran P_{\pi_j\vee\dots\vee\pi_1}((I-\pi_j\pi_j^{\dag})+\pi_jP_+\pi_j^{\dag})=
\mathcal{E}_{j+}\oplus\Ran(I-\pi_j\pi_j^{\dag})P_{\pi_j\vee\dots\vee\pi_1}\\[4pt]
&=&\mathcal{E}_{j+}\oplus\vee_{k=1}^{j-1}\clos\Ran(I-\pi_j\pi_j^{\dag})\pi_k=
\mathcal{E}_{j+}\oplus(\vee_{k=1}^{j-1}((\mathcal{E}_{k}\vee\mathcal{E}_{j})\ominus\mathcal{E}_{j}))\,.
\end{array}
\]
Hence we get
$\,\mathcal{H}_{j1+}\subset\mathcal{H}_{j+}=\mathcal{E}_{j+}\oplus\mathcal{H}_{j+}^{nor}\,$ and
\[
\mathcal{H}_{uj}=\mathcal{H}_{j+}\ominus\mathcal{H}_{j1+}=\mathcal{H}_{j+}^{nor}\ominus\
(\vee_{k=1}^{j-1}((\mathcal{E}_{k}\vee\mathcal{E}_{j})\ominus\mathcal{E}_{j}))\,.
\]
It is obvious that the subspace $\,\mathcal{H}_{uj}\,$ reduces the operator $\,\mathcal{U}\,$.

Assume that $\,\mathcal{H}_{j1+}'=\mathcal{H}_{\pi_j'\vee\dots\vee\pi_1'}\cap\Ker
\pi_j'P_-{\pi_j'}^{\dag}\,$, $\,\mathcal{H}_{j1+}'\subset\mathcal{H}_{j+}\,$ and the subspace
$\,\mathcal{H}_{j+}\ominus\mathcal{H}_{j1+}'\,$ reduces the operator $\,\mathcal{U}\,$, where
$\,\Pi'=(\pi_1',\ldots,\pi_n')\in\Mod_n\,$. Then we have the generalized Wold
decomposition~\cite{AD} $\,\mathcal{H}_{j+}=\mathcal{E}_{j+}'\oplus(\vee_{k=1}^{j-1}
((\mathcal{E}_{k}'\vee\mathcal{E}_{j}')\ominus\mathcal{E}_{j}'))
\oplus(\mathcal{H}_{j+}\ominus\mathcal{H}_{j1+}')\,$. Since this decomposition is unique, we obtain
$\,\mathcal{E}_{j+}'=\mathcal{E}_{j+}\,$, $\,\mathcal{E}_{j}'=\mathcal{E}_{j}\,$ and, by induction,
$\,\mathcal{H}_{j1+}'=\mathcal{H}_{j1+}\,$. Then, $\,\pi_j'=\pi_j\psi_j\,$, where
$\,\psi_j=\pi_j^{\dag}\pi_j'\, ,\;\psi_j^{-1}={\pi_j'}^{\dag}\pi_j\in
H^{\infty}(G_+,\mathcal{L}(\mathfrak{N}_j))\,$.

Since $\,\mathcal{H}_{uj}\,\bot\,(\mathcal{E}_{j+}\oplus
((\mathcal{E}_{k}\vee\mathcal{E}_{j})\ominus\mathcal{E}_{j}))\,$ and
$\,\mathcal{H}_{uj}\subset\mathcal{H}_{j+}\,$, we get
$\,\mathcal{H}_{uj}\,\bot\,(\mathcal{E}_{k}\vee\mathcal{E}_{j})\,$. For $i>j\,$, we have
$\,\mathcal{E}_{i-}\subset\mathcal{H}_{j+}^{\bot}\subset\mathcal{H}_{uj}^{\bot}\,$. Hence,
$\,\mathcal{H}_{uj}\,\bot\,\mathcal{E}_{i}\,$ and
$\,\mathcal{H}_{uj}\,\bot\,\mathcal{H}_{\pi_n\vee\dots\vee\pi_1}\,$. Since
$\,\mathcal{H}_{j1+}\subset\mathcal{H}_{\pi_n\vee\dots\vee\pi_1}\,$,
$\,\mathcal{H}_{uj}\,\subset\,\mathcal{H}_{\pi_n\vee\dots\vee\pi_1}^{\bot}\,$ and
$\,\mathcal{H}_{j+}=\mathcal{H}_{uj}\oplus\mathcal{H}_{j1+}\,$, we have
$\,\mathcal{H}_{uj}=\mathcal{H}_{\pi_n\vee\dots\vee\pi_1}^{\bot}\cap\mathcal{H}_{j+}\,$ and
therefore $\,\mathcal{H}_{u1}\subset\ldots\subset\mathcal{H}_{un}\,$.\proofend

Let $\,\theta\in\Cfn\,$. We fix $\,\theta\,$ and define $\,\Mod_n^{\theta} := \{\Pi\in\Mod_n :\,
\pi_n^{\dag}\pi_1=\theta  \}\,$. Then we can consider the chain of subspaces
$\;\,\mathcal{F}_{im}^{\theta}(\Pi):=(\mathcal{K}_{(11)}\subset\mathcal{K}_{(21)}\subset\ldots\subset
\mathcal{K}_{(n1)})\,$, where $\,\mathcal{K}_{(k1)}=\Ran P_{(k1)}\,$. The subspaces
$\,\mathcal{K}_{(k1)}\,$ are invariant under the operator $\,\widehat{T}\,$ and this observation
motivates the following definition.

Let $\,\theta=\pi_-^{\dag}\pi_+\,$, where the operators
$\,\pi_{\pm}\in\mathcal{L}(L^2(\Xi_{\pm}),\mathcal{H})\,$ are isometries. Let
$\,\mathcal{U}\in\mathcal{L}(\mathcal{H})\,$ be a normal operator such that
$\,\mathcal{U}\pi_{\pm}=\pi_{\pm} z\,$ and $\,\sigma(\mathcal{U})\subset C\,$. Let also
$\,\mathcal{K}=\Ran P\,,\,P=(I-\pi_+P_+\pi_+^{\dag})(I-\pi_-P_-\pi_-^{\dag})\,$, and
$\,T=P\mathcal{U}|\mathcal{K}\,$.
\begin{defin}
A chain of subspaces $\,L=(L_{1}\subset L_{2}\subset\ldots\subset L_{n})\,$ is called
n-invariant if $\,L_n\subset\mathcal{K}\,$, $\,(T-z)^{-1}L_k\subset L_k\,,\;z\in
G_-\,,\;k=\overline{1,n}\,$, and the subspaces $\,L_1\,$, $\,\mathcal{K}\ominus L_n\,$
reduce the operator $\,\mathcal{U}\,$. We will denote the class of all n-invariant
chains by~$\,\Inv_n^{\theta}\,$.
\end{defin}
\noindent In fact, we have already defined the transformation $\,\mathcal{F}_{im}^{\theta} \colon
\Mod_n^{\theta}\to\Inv_n^{\theta}\,$, which takes each $\,\Pi\in\Mod_n^{\theta}\,$ to the
n-invariant chain of subspaces $\;\,(\mathcal{K}_{(11)}\subset
\mathcal{K}_{(21)}\subset\ldots\subset \mathcal{K}_{(n1)})\in\Inv_n^{\theta}\,$. This
transformation is surjective accurate up to the ``normal'' part of the chain.
\begin{prop}\label{t_chain}
Suppose a chain $\,L\,$ is n-invariant. Then there exists an n-model $\,\Pi\in\Mod_n^{\theta}\,$
such that $\,\mathcal{K}_{(k1)}\subset L_{k}\,,\;k=\overline{1,n}\,$ and the subspaces
$\,L_{uk}:=L_{k}\ominus\mathcal{K}_{(k1)}\,$ reduce the operator $\,\mathcal{U}\,$. If an n-model
$\,\Pi'\in\Mod_n\,$ satisfies the same conditions, then
$\,\mathcal{K}_{(k1)}'=\mathcal{K}_{(k1)}\,$ and $\,\exists\,\psi_k\,$ such that
$\,\psi_k,\psi_k^{-1}\in H^{\infty}(G_+,\mathcal{L}(\mathfrak{N}_k))\,$ and
$\,\pi_k'=\pi_k\psi_k\,$.  Besides, we have $\,L_{u1}\subset\ldots\subset L_{un}\,$.
\end{prop}
\myproof We put $\,\mathcal{H}_{k+}=L_{k}\dot+\mathcal{D}_{+}\,$, where $\,\mathcal{D}_{+}=\Ran
q_+\,$, $\,q_+=\pi_+P_+\pi_+^{\dag}\,$. Then, for $\;z\in G_-$, we get
$\;\,(\mathcal{U}-z)^{-1}\mathcal{D}_{+}\subset\mathcal{D}_{+}\subset\mathcal{H}_{k}\;\,$ and
\[
(\mathcal{U}-z)^{-1}L_{k} \subset P(\mathcal{U}-z)^{-1}L_{k}\dot+ q_{+}(\mathcal{U}-z)^{-1}L_{k}
\subset (T-z)^{-1}L_{k}\dot+\mathcal{D}_{+}\subset\mathcal{H}_{k}\,.
\]
Therefore the chain $\,\mathcal{H}_{1+}\subset\ldots\subset\mathcal{H}_{n+}\,$ is
invariant under $\;(\mathcal{U}-z)^{-1}\,,\;z\in G_-$. By Prop.~\ref{u_chain}, there
exists an n-model $\,\Pi\in\Mod_n\,$ such that
$\,\mathcal{H}_{k1+}\subset\mathcal{H}_{k+}\,,\;k=\overline{1,n}\,$ and the subspaces
$\,\mathcal{H}_{uk}=\mathcal{H}_{k+}\ominus\mathcal{H}_{k1+}\,$ reduce the operator
$\,\mathcal{U}\,$. Since $L_1$ reduces $\mathcal{U}$, we have that
$\,\mathcal{H}_{1+}=\mathcal{D}_{+}\dot+L_{1}\,$ is the generalized Wold decomposition
of $\,\mathcal{H}_{1+}\,$. Taking into account the uniqueness of Wold decomposition,
we obtain $\,\pi_+=\pi_1\psi_1\,$. Comparing  the Wold decompositions of the equal
subspaces $\,\mathcal{K}\dot+\mathcal{D}_{+}\,$ and $\,(\mathcal{K}\ominus
L_n)\dot+\mathcal{H}_{n+}\,$, we obtain $\,\pi_-=\pi_n\psi_n\,$. Thus we can assume
without loss of generality (see the proof of Prop.~\ref{u_chain}) that
$\,\pi_1=\pi_+\,$ and $\,\pi_n=\pi_-\,$, i.e., $\,\Pi\in\Mod_n^{\theta}\,$.

Since $\,L_k\subset\mathcal{K}\,$, we have
$\,L_k=(I-\pi_+P_+\pi_+^{\dag})\mathcal{H}_{k+}\,$. Taking into account that
$\,\mathcal{K}_{(k1)}=(I-\pi_+P_+\pi_+^{\dag})\mathcal{H}_{k1+}\,$ and
$\,\mathcal{H}_{k1+}\subset\mathcal{H}_{k+}\,$, we get $\,\mathcal{K}_{(k1)}\subset
L_k\,$. Since
$\,\mathcal{H}_{uk}=\mathcal{H}_{\pi_n\vee\dots\vee\pi_1}^{\bot}\cap\mathcal{H}_{k+}\,$,
we have $\,\mathcal{H}_{uk}=P\mathcal{H}_{uk}\subset P\mathcal{H}_{k+}=L_k\,$ and
therefore $\,\mathcal{H}_{uk}\oplus\mathcal{K}_{(k1)}\subset L_k\,$. In fact, these
two spaces are equal. Consider the operator $\,q_{+}'=\pi_+'P_{+}'\pi_+'^{*}\,$, which
is the orthogonal counterpart to $\,q_{+}=\pi_+P_+\pi_+^{\dag}\,$ (see the comments
between Prop.\ref{proj_prop} and Prop.\ref{dir_sum}). Put $\,L_k':=(I-q_{+}')L_k\,$
and $\,\mathcal{K}_{(k1)}':=(I-q_{+}')\mathcal{K}_{(k1)}\,$. By Corollary of
Lemma~\ref{lem_iii}, $\,L_k=(I-q_{+})L_k'\,$ and
$\,\mathcal{K}_{(k1)}=(I-q_{+})\mathcal{K}_{(k1)}'\,$. Further, we have
\[
L_k'\ominus\mathcal{K}_{(k1)}'=(L_k'\oplus\mathcal{D}_{+})
\ominus(\mathcal{K}_{(k1)}'\oplus\mathcal{D}_{+})=\mathcal{H}_{k+}\ominus\mathcal{H}_{k1+}=
\mathcal{H}_{uk}\,.
\]
Then
\[
L_k=(I-q_{+})L_k'=(I-q_{+})(\mathcal{K}_{(k1)}'\oplus\mathcal{H}_{uk})
=\mathcal{K}_{(k1)}\dot+\mathcal{H}_{uk}=\mathcal{K}_{(k1)}\oplus\mathcal{H}_{uk}
\]
and therefore $\,L_k=\mathcal{K}_{(k1)}\oplus\mathcal{H}_{uk}\,$. Hence,
$\,L_{uk}=L_k\ominus\mathcal{K}_{(k1)}=\mathcal{H}_{uk}\,$. Then, by
Prop.~\ref{u_chain}, we have $\,L_{u1}\subset\ldots\subset L_{un}\,$.

Let $\,\Pi'\in\Mod_n^{\theta}\,$ be a n-model such that $\,\mathcal{K}_{(k1)}'\subset
L_{k}\,,\;k=\overline{1,n}\,$ and the subspaces $\,L_{uk}'=L_{k}\ominus\mathcal{K}_{(k1)}'\,$
reduce the operator $\,\mathcal{U}\,$. Then
$\,\mathcal{H}_{k1+}'=\mathcal{K}_{(k1)}'\dot+\mathcal{D}_{+}\subset
L_k\dot+\mathcal{D}_{+}=\mathcal{H}_{k+}\,$ and the subspaces
$\,\mathcal{H}_{k+}\ominus\mathcal{H}_{k1+}'=L_{uk}'\,$ reduce the operator $\,\mathcal{U}\,$. By
Prop.~\ref{u_chain}, we get $\,\mathcal{H}_{k1+}'=\mathcal{H}_{k1+}\,$. Hence,
$\,\mathcal{K}_{(k1)}'=(I-\pi_+P_+\pi_+^{\dag})\mathcal{H}_{k1+}'=
(I-\pi_+P_+\pi_+^{\dag})\mathcal{H}_{k1+}=\mathcal{K}_{(k1)}\,$.\proofend

\myrem\,\, In the case $\,n=2\,$ this proposition is an analogue of the well-known decomposition of
a contraction $T$ into the orthogonal sum $\,T=T_{cnu}\oplus T_u\,$ of the completely non-unitary
part $\,T_{cnu}\,$ and the unitary part $\,T_u\,$ (see~\cite{NF}). In this connection, we will use
the notation
\[
\Inv_n^{\theta\,cnu}:=\{(L_{1}\subset L_{2}\subset\ldots\subset L_{n})\in\Inv_n^{\theta} :
L_{un}=\{0\}\}
\]
In this notation Prop.~\ref{t_chain} means merely that
$\,\Ran\mathcal{F}_{im}^{\theta}=\Inv_n^{\theta\,cnu}\,$. Note also that the condition
$\,L_{un}=\{0\}\,$ is equivalent to the condition
\[
\vee_{k=1}^{n}[\vee_{z\notin C}(\mathcal{U}-z)^{-1}([\vee_{z\notin
C}(\mathcal{U}-z)^{-1}(L_{k}\dot+\mathcal{D}_{+})]\ominus(L_{k}\dot+\mathcal{D}_{+}))]=\mathcal{H}\,.
\]

Let us now return to the transformation $\,\mathcal{F}_{ic}\,$. Fix $\,\theta\in\Cfn\,$ and define
$\,\Mod_n^{\theta}:=\{\Theta\in\Mod \colon \Theta_{n1}=\theta\}\,$. Then we can consider the
restriction $\,\mathcal{F}_{ic}|\Mod_3^{\theta}\,$, which takes each 3-characteristic function
$\,\Theta\in\Mod_3^{\theta}\,$ to the invariant subspace
$\,L:=\mathcal{K}_{(21)}\subset\mathcal{H}\,$. The main difficulty to handle effectively
factorizations of the function $\,\theta\,$ is the fact that the space $\,\mathcal{H}\,$ is
variable and we cannot compare invariant subspaces when we run over factorizations of $\,\theta\,$.
To avoid this effect we shall restrict ourselves to models for which
$\,\mathcal{H}=\mathcal{H}_{\pi_+\vee\pi_-}=\Ran\pi_+\vee\Ran\pi_-\,$ and
$\,\Pi=(\pi_+,\pi_2,\pi_-)\,$, where $\,\pi_{\pm}\in\mathcal{L}(L^2(\Xi_{\pm}),\mathcal{H})\,$ are
some fixed isometries such that $\,\theta=\pi_-^{\dag}\pi_+\,$. Then we obviously have $\,\Ran\pi_2
\subset \Ran\pi_1\vee\Ran\pi_3\,$ and therefore $\,\Pi=(\pi_1,\pi_2,\pi_3)\in\Mod_3^{reg}\,$.

In this connection we define the subclasses
\[
\Cfn_n^{\theta\,reg}:=\Cfn_n^{\theta}\cap\Cfn_n^{reg}\,,\qquad
\Mod_n^{\theta\,reg}:=\Mod_n^{\theta}\cap\Mod_n^{reg}
\]
and
\[
\Inv_n^{\theta\,reg}:=\{ L\in\Inv_n^{\theta} : \Ran\pi_+\vee\Ran\pi_-=\mathcal{H}\}\,.
\]
By Prop.\ref{reg_M_C}, it can easily be shown that
\[
\mathcal{F}_{im}^{\theta}(\Pi)\in\Inv_n^{\theta\,reg}\,\Longleftrightarrow\,\Pi\in\Mod_n^{\theta\,reg}\,.
\]
Besides, it is clear that $\,\Inv_n^{\theta\,reg}\subset\Inv_n^{\theta\,cnu}\,$.

\vskip 3pt Finally, we define the transformation $\,\mathcal{F}_{ic}^{\theta} :
\Cfn_n^{\theta\,reg}\to\Inv_n^{\theta\,reg}\,$ by the following procedure. Let
$\,\Theta\in\Cfn_n^{\theta\,reg}\,$ and $\,\Pi=\mathcal{F}_{mc}(\Theta)\in\Mod_n\,$
(in fact, by Prop.~\ref{reg_M_C}, $\,\Pi\in\Mod_n^{reg}\,$). Then,
$\,\mathcal{H}=\mathcal{H}_{\pi_1\vee\pi_n}=\Ran\pi_1\vee\Ran\pi_n\,$ and
$\,\theta=\pi_n^{\dag}\pi_1\,$. By Prop.~\ref{CtoM}, there exists an unique unitary
operator $\,X \colon \mathcal{H}_{\pi_1\vee\pi_n} \to \mathcal{H}_{\pi_+\vee\pi_-}\,$
such that $\,\pi_+=X\pi_1\,$ and $\,\pi_-=X\pi_n\,$. Then we put
$\,\mathcal{F}_{ic}^{\theta}(\Theta):=\mathcal{F}_{im}^{\theta}(X\Pi)\in\Inv_n^{\theta\,reg}\,$,
where $\,X\Pi=(X\pi_1, X\pi_2, \ldots, X\pi_n)\in\Mod_n^{\theta\,reg}\,$. This
definition of the fundamental transformation $\,\mathcal{F}_{ic}^{\theta}\,$ is rather
indirect. As justification of it we note that in the case of the unit disk the known
approaches~\cite{NF,Br,NV} are not simpler than our procedure. The following
Proposition is a straightforward consequence of Prop.~\ref{t_chain}.
\begin{prop}\label{F_cm}
One has \,\,\,\,1)\, $\,\Ran\mathcal{F}_{ic}^{\theta}=\Inv_n^{\theta\,reg}\,$;\,\,\, 2)\, If\,
$\,\mathcal{F}_{ic}^{\theta}(\Theta')=\mathcal{F}_{ic}^{\theta}(\Theta)\,$,
$\,\Theta,\Theta'\in\Cfn_n^{\theta\,reg}\,$, then $\,\Theta'\sim\Theta\,$, where $\,\sim\,$ is
equivalence relation: $\,\Theta'\sim\Theta\,$ if $\,\exists\,\psi_k\,,\;k=\overline{2,n-1}\,$ such
that $\,\psi_k,\,\psi_k^{-1}\in H^{\infty}(G_+,\mathcal{L}(\mathfrak{N}_k))\,$,
$\,\Theta_{ij}'=\psi_i^{-1}\Theta_{ij}\psi_j\,$, and $\,\Xi_k'=\psi_k^*\Xi_k\psi_k\,$;
$\,\psi_1=I,\,\psi_n=I\,$.
\end{prop}
\noindent Thus, one can consider the quotient space
$\,\Cfn_n^{\theta\,reg\sim}:=\Cfn_n^{\theta\,reg}/_\sim\,$ and the corresponding one-to-one
transformation $\,\mathcal{F}_{ic}^{\theta\sim} :
\Cfn_n^{\theta\,reg\sim}\to\Inv_n^{\theta\,reg}\,$. Note that the functions $\,\psi_k\,$ can be
regarded as $\Xi$-unitary constants, i.e., $\,\psi_k^{\dag}=\psi_k^{-1}\in H^\infty(G_+,{\mathcal
L}({\mathfrak N}_k))\,$, where $\,\psi_k^{\dag}\,$ are adjoint to $\,\psi_k\colon L^2(\Xi_k') \to
L^2(\Xi_k)\,$.

Let us consider particular cases. In the case of $\,n=3\,$ we obtain that the transformation
$\,\mathcal{F}_{ic}^{\theta\sim} : \Cfn_3^{\theta\,reg\sim}\to\Inv_3^{\theta\,reg}\,$ is an
one-to-one correspondence between regular factorizations of a characteristic function and invariant
subspaces of the corresponding model operator.

Consider the case $\,n=4\,$. Let $\,L=(L_1,L_2,L_3,L_4)\in\Inv_4^{\theta\,reg}\,$. By
Prop.~\ref{F_cm}, there exists $\,\Theta\in\Cfn_4^{\theta\,reg}\,$ such that
$\,L=\mathcal{F}_{ic}^{\theta}(\Theta)\,$. If we rename $\,L'=L_2,\,L''=L_3\,$ (recall that
$\,L_1=\{0\},\,L_4=\mathcal{K}_{\theta}\,$) and
$\,\theta=\Theta_{41},\;\theta_1'=\Theta_{21},\;\theta_2'=\Theta_{42},
\;\theta_1''=\Theta_{31},\;\theta_2''=\Theta_{43}\,$, $\,\Xi_+=\Xi_1\,$, $\,\Xi'=\Xi_2\,$,
$\,\Xi''=\Xi_3\,$, $\,\Xi_-=\Xi_4\,$, then we have
\[
\theta=\theta_2'\theta_1' =\theta_2''\theta_1''\quad\mbox{and}\quad\exists\:\vartheta\in
S_{\Xi}\quad\mbox{such that}\quad
\theta_1''=\vartheta\theta_1'\,,\quad\theta_2'=\theta_2''\vartheta\,. \eqno{(\prec)}
\]
Certainly, $\,\vartheta=\Theta_{32}\,$ and $\,\Xi=(\Xi',\Xi'')\,$. We shall say that the
factorization $\,\theta=\theta_2'\theta_1'\,$ precedes the factorization
$\,\theta=\theta_2''\theta_1''\,$ (and write $\,\theta_2'\theta_1'\prec\theta_2''\theta_1''\,$) if
the condition ($\prec$) is satisfied. Thus, $\,L'\subset
L''\,\Longrightarrow\,\theta_2'\theta_1'\prec\theta_2''\theta_1''\,$.

Conversely, suppose that factorizations $\,\theta_2'\theta_1'=\theta_2''\theta_1''\,$ are regular
and $\,\theta_2'\theta_1'\prec\theta_2''\theta_1''\,$. After backward renaming we have
$\,\Theta\in\Cfn_4^{\theta\,reg}\,$. Let $\,L=\mathcal{F}_{ic}^{\theta}(\Theta)\,$,
$\,L'=\mathcal{F}_{ic}^{\theta}(\theta_2'\theta_1')\,$, and
$\,L''=\mathcal{F}_{ic}^{\theta}(\theta_2''\theta_1'')\,$. Since the factorizations are regular, we
have $\,L'=L_2, L''=L_3\,$. Therefore,
$\,\theta_2'\theta_1'\prec\theta_2''\theta_1''\,\Longrightarrow\,L'\subset L''\,$. Finally, we have
\[
\mathcal{F}_{ic}^{reg}(\Theta_{42}\Theta_{21})\subset\mathcal{F}_{ic}^{reg}(\Theta_{43}\Theta_{31})\,
\,\Longleftrightarrow\,\Theta_{42}\Theta_{21}\prec\Theta_{43}\Theta_{31}\,.
\]
It is easy to check that
$\,\theta_2'\theta_1'\prec\theta_2''\theta_1''\,,\theta_2'\theta_1'\sim\vartheta_2'\vartheta_1'\,,
\theta_2''\theta_1''\sim\vartheta_2''\vartheta_1''\,\Longrightarrow\,
\vartheta_2'\vartheta_1'\prec\vartheta_2''\vartheta_1''\,$, i.e., the order relation $\,\prec\,$ is
well defined on the quotient space $\,\Cfn_3^{\theta\,reg\sim}\,$. Taking all this into account, we
arrive at the main result of the Section.
\begin{theoremb}
There is an order preserving one-to-one correspondence $\,\mathcal{F}_{ic}^{reg}\,$
between regular factorizations of a characteristic function (up to the equivalence
relation) and invariant subspaces of the resolvent $\,(\widehat{T}-z)^{-1}\,,\;z\in
G_-\,$ of the corresponding model operator.
\end{theoremb}
\noindent This Theorem is an extension of the fundamental result from~\cite{NF} (Theorems VII.1.1
and VII.4.3;\, see also~\cite{K} for some refinement).
\begin{cor}
Suppose that factorizations $\,\theta_2'\theta_1',\,\theta_2''\theta_1''\,$ are regular,
$\,\theta_2'\theta_1'\prec\theta_2''\theta_1''\,$ and
$\,\theta_2''\theta_1''\prec\theta_2'\theta_1'\,$. Then
$\,\theta_2'\theta_1'\sim\theta_2''\theta_1''\,$.
\end{cor}
\proof Let $\,L'=\mathcal{F}_{ic}^{reg}(\theta_2'\theta_1')\,$ and
$\,L''=\mathcal{F}_{ic}^{reg}(\theta_2''\theta_1'')\,$. By Theorem B, we get $\,L'
\subset L'' \subset L'\,$ and therefore $\,L'=L''\,$. Then, by Prop.~\ref{F_cm}, we
have $\,\theta_2'\theta_1'\sim\theta_2''\theta_1''\,$. \proofend

\noindent Note that the Corollary can be proved independently from Theorem B. The
corresponding argumentation make use of Lemmas~\ref{lem_fact} and ~\ref{lem_reg} and
therefore we can drop the assumptions that
$\,\theta_2',\,\theta_2',\,\theta_1'',\,\theta_2''\,$ are operator valued functions.

\begin{prp}
Let $\,A_{21},A_{42},A_{31},A_{43}\,$ be contractions. Suppose that factorizations
$\,A_{42}\cdot A_{21},\; A_{43}\cdot A_{31}\,$ are regular, $\,A_{42}\cdot A_{21}\prec
A_{43}\cdot A_{31}\,$ and $\,A_{43}\cdot A_{31}\prec A_{42}\cdot A_{21}\,$. Then there
exists an unitary operator $\,U\,$ such that $\,A_{31}=UA_{21}\,$ and
$\,A_{43}=A_{42}U^{-1}\,$.
\end{prp}
\proof We shall make use of the following two lemmas.
\begin{lemm}\label{lem_vii}
Suppose that $\,||A||\le 1\,$ and $\,A|H_1=I|H_1\,$. Then $\,A^*|H_1=I|H_1\,$.
\end{lemm}
\proof We have $\,A=\left(
\begin{array}{cc}
  I & a_{12} \\
  0 & a_{22} \\
\end{array}\right)\,$. Then
\[
0\le ((I-A^*A)\,\left(
\begin{array}{c}
  f_1  \\
  0  \\
\end{array}\right)\,,\,\left(
\begin{array}{c}
  f_1  \\
  0  \\
\end{array}\right))=-(a_{12}^*f_1,a_{12}^*f_1)\le 0\,.
\]
Therefore, $\,a_{12}=0\,$ and
$\,A=\left(
\begin{array}{cc}
  I & 0 \\
  0 & a_{22} \\
\end{array}\right)
\,$. \proofend

\begin{lemm}\label{lem_viii}
Let $\,A_{21},A_{32}\,$ be contractions. Suppose that factorization $\,A_{32}\cdot
A_{21}\,$ is regular. Then $\,(\Ran A_{32}^*\vee\Ran A_{21})^{\bot}=\{0\}\,$.
\end{lemm}
\proof Let $\,f\,\bot\,(\Ran A_{32}^*\vee\Ran A_{21})\,$. Then $\,f\in\Ker A_{32}\,$
and $\,f\in\Ker A_{21}^*\,$. Hence, $\,(I-A_{21}A^*_{21})f=f\,$ and therefore
$\,(I-A_{21}A^*_{21})^{1/2}f=f\,$. Similarly, we have
$\,(I-A^*_{32}A_{32})^{1/2}f=f\,$. Then $\,f\in \Ran\,(I-A^*_{32}A_{32})^{1/2}\cap
\Ran\,(I-A_{21}A^*_{21})^{1/2}=\{0\}\,$. \proofend From the definition of the order
relation $\,\prec\,$ we get that there exists contractions $\,A_{32},A_{23}\,$ such
that $\,A_{42}=A_{43}A_{32}\,$, $\,A_{31}=A_{32}A_{21}\,$, $\,A_{43}=A_{42}A_{23}\,$,
and $\,A_{32}=A_{23}A_{31}\,$. Let $\,A=A_{23}A_{32}\,$. Then we have
$\,A_{21}=AA_{21}\,$ and $\,A_{42}=A_{42}A\,$ and therefore $\,A|\Ran A_{21}=I|\Ran
A_{21}\,$ and $\,A^*|\Ran A_{42}=I|\Ran A_{42}\,$. By Lemma~\ref{lem_vii}, $\,A|\Ran
A_{42}=I|\Ran A_{42}\,$. Finally, by Lemma~\ref{lem_viii}, we get $\,A=A|(\Ran
A_{21}\vee\Ran A_{42})=I|(\Ran A_{21}\vee\Ran A_{42})=I\,$, that is,
$\,A_{23}A_{32}=I\,$. Likewise, we get $\,A_{32}A_{23}=I\,$. Since $\,A_{23}\,$ and
$\,A_{32}\,$ are contraction, they are unitary operators. It remains to put
$\,U=A_{32}\,$. \proofend

\vskip 5pt In conclusion we again consider curved conservative systems. The following
assertion is just a translation of Prop.~\ref{t_chain} into the language of systems.
\begin{prop}\label{sys_inv}
Suppose $\,\Sigma=(T, M, N)\in\Sys\,$ and a subspace $\,L\,$ is invariant under the
resolvent $\,(T-z)^{-1}\,,\;z\in G_-\,$. Then there exist systems
$\,\Sigma_1,\,\Sigma_2\in\Sys\,$ and an operator $\,X : H_1\oplus H_2 \to H\,$ such
that $\;\Sigma\;{\stackrel{X}{\sim}}\;\Sigma_2\cdot\Sigma_1\,$ and $\,L=XH_1\,$.
\end{prop}
\proof Let
$\,\Sigma\;{\stackrel{Y}{\sim}}\;\widehat{\Sigma}=\mathcal{F}_{sc}(\theta)\,$ and
$\,\Pi=(\pi_1,\pi_3)=\mathcal{F}_{mc}(\theta)\,$. Then $\,\widehat{L}=YL\,$ is an
invariant subspace for the model operator. By Theerem~B, there exists an regular
factorization $\,\theta=\theta_2\cdot\theta_1\,$ such that
$\,\widehat{L}=\mathcal{F}_{ic}^{reg}(\theta_2\cdot\theta_1)=\Ran P_{(21)}\,$.
Besides, $\,\theta_1=\pi_2^\dag\pi_1\,$ and $\,\theta_2=\pi_3^\dag\pi_2\,$. We put
$\,\Pi_1=(\pi_1,\pi_2)\,$, $\,\Pi_2=(\pi_2,\pi_3)\,$,
$\,\widehat{\Sigma}_1=\mathcal{F}_{sc}(\Pi_1)\,$, and
$\,\widehat{\Sigma}_2=\mathcal{F}_{sc}(\Pi_2)\,$. Let
$\,\Sigma_1\;{\stackrel{Y_1}{\sim}}\;\widehat{\Sigma}_1\,$ and
$\,\Sigma_2\;{\stackrel{Y_2}{\sim}}\;\widehat{\Sigma}_2\,$. By
Prop.~\ref{M_and_S_con}, we get  $\,\Sigma_1\cdot\Sigma_2\sim\widehat{\Sigma}\,$ with
the operator $\,P_{(31)}(Y_1,Y_2)\,$ realizing the similarity. It can easily be
checked that $\,\widehat{L}=P_{(31)}(Y_1,Y_2)H_1\,$. Then, for
$\,X=Y^{-1}P_{(31)}(Y_1,Y_2)\,$, we get
$\;\Sigma\;{\stackrel{X}{\sim}}\;\Sigma_2\cdot\Sigma_1\,$ and $\,L=XH_1\,$. \proofend

\noindent Besides, we have the following assertion.
\begin{prop}
Suppose the system $\,\Sigma_2\cdot\Sigma_1=\Sigma_2'\cdot\Sigma_1'\,$ is simple, $\,H_1=H_1'\,$,
and $\,\Theta_2\Theta_1=\Theta_2'\Theta_1'\,$. Then there exists $\,\psi\,$ such that
$\,\psi,\,\psi^{-1}\in H^{\infty}(G_+,\mathcal{L}(\mathfrak{N}))\,$ and
$\,\Sigma_1'\sim\Sigma_1''=(T_1, M_1, N_1'')\,$, where $\; {N_1''}^*f_1=\displaystyle-\frac{1}{2\pi
i}\int_{\overline{C}}\psi(z)^*\,[N_1^*(T_1^*-\cdot)^{-1}f_1]_-(z)\,dz\,$, $\,f_1\in H_1\,$.
\end{prop}
\proof Let
$\,\Sigma=\Sigma_2\cdot\Sigma_1\;{\stackrel{Y}{\sim}}\;\widehat{\Sigma}=\mathcal{F}_{sc}(\theta)\,$
and $\,\Pi=(\pi_1,\pi_3)=\mathcal{F}_{mc}(\theta)\,$. Using the same notation as in
the proof of Prop.~\ref{sys_inv}, we obtain that the operators $\,P_{(31)}(Y_1,Y_2)\,$
and $\,P_{(31)}(Y_1',Y_2')\,$ realize the similarities
$\,\Sigma_2\cdot\Sigma_1\sim\widehat{\Sigma}\,$ and
$\,\Sigma_2'\cdot\Sigma_1'\sim\widehat{\Sigma}\,$, respectively. Since the system
$\,\widehat{\Sigma}\sim{\Sigma}\,$ is simple, by Prop.~\ref{simpl}(4), we get
$\,P_{(31)}(Y_1,Y_2)=P_{(31)}(Y_1',Y_2')\,$ and therefore
$\,P_{(31)}(Y_1,Y_2)H_1=P_{(31)}(Y_1',Y_2')H_1'\,$. Then, by Prop.~\ref{F_cm}, there
exists an operator valued function $\,\psi\,$ such that $\,\psi,\,\psi^{-1}\in
H^{\infty}(G_+,\mathcal{L}(\mathfrak{N}_2))\,$, $\,\theta_1'=\psi^{-1}\theta_1\,$, and
$\,\theta_2'=\theta_2\psi\,$. According to~\cite{T2},
$\,\Sigma_1''\sim\widehat{\Sigma}''=\mathcal{F}_{sc}(\psi^{-1}\theta_1)\,$. Since
$\,\Sigma_1'\sim\widehat{\Sigma}'=\mathcal{F}_{sc}(\theta_1')\,$, we get
$\,\Sigma_1'\sim\Sigma_1''\,$. \proofend

Further, we shall say that a system $\,\Sigma\in\Sys\,$ possesses the property of
uniqueness of characteristic function if there exists an unique characteristic
function $\,\Theta\in\Cfn\,$ such that $\,\Sigma=\mathcal{F}_{cs}(\Theta)\,$. Recall
(see the Introduction) the sufficient condition for this property: the transfer
function $\,\Upsilon(z)\,$ of the system $\,\Sigma\,$ is an operator valued function
of Nevanlinna class. For products of systems we have the following (non-trivial) fact:
\textit{suppose that a system $\,\Sigma=\Sigma_2\Sigma_1\,$ is simple, possesses the
property of uniqueness, and $\,\rho(T_1)\cap G_+\ne\emptyset\;$; then the system
$\,\Sigma_1\,$ possesses the same property too}.

\begin{prop}
Suppose the system $\,\Sigma_2\cdot\Sigma_1=\Sigma_2'\cdot\Sigma_1\,$ is simple and possesses the
property of uniqueness. Suppose also $\,\rho(T_1)\cap G_+\ne\emptyset\;$. Then
$\,\Sigma_2=\Sigma_2'\,$.
\end{prop}
\proof Let $\,\Sigma=\Sigma_2\cdot\Sigma_1=\Sigma_2'\cdot\Sigma_1\,$ and
$\,\theta=\theta_2\theta_1=\theta_2'\theta_1'\,$ be the corresponding factorizations.
Then $\,\theta_1=\theta_1'\,$ (see the comments before the Proposition). Since
$\,\forall\;\lambda\in\rho(T_1)\cap G_+\ne\emptyset\;$
$\,\exists\;\theta_1(\lambda)^{-1}\,$, we get $\,\theta_2=\theta_2'\,$. Then
$\,\Sigma_2\sim\mathcal{F}_{sc}(\theta_2)\,$,
$\,\Sigma_2'\sim\mathcal{F}_{sc}(\theta_2)\,$ and therefore
$\,\Sigma_2\;{\stackrel{X_2}{\sim}}\;\Sigma_2'\,$. Taking this into account, we have
$\,\Sigma\;{\stackrel{I}{\sim}}\;\Sigma\,$ and $\,\Sigma\;{\stackrel{I\oplus
X_2}{\sim}}\;\Sigma\,$. By Prop.~\ref{simpl}(4), we get $\,X_2=I\,$. \proofend

\end{document}